\pdfoutput=1
\documentclass[preprint,3p,times]{elsarticle}

\usepackage{graphicx}
\usepackage{epsfig}

\usepackage{amssymb}
\usepackage{amsthm}
\usepackage{amsbsy}
\usepackage{amsmath}
\usepackage{amsfonts}
\usepackage{dsfont}
\usepackage{rotating}
\usepackage{color}

\usepackage{pgf,tikz}
\usepackage{mathrsfs}
\usepackage{amssymb,fancyhdr,txfonts,pxfonts}
\usetikzlibrary{arrows}

\setcitestyle{square}

\DeclareMathAlphabet{\mathpzc}{OT1}{pzc}{m}{it}

\journal{Computer Methods in Applied Mechanics and Engineering}

\begin{document}

\begin{frontmatter}



\title{Semi-implicit discontinuous Galerkin methods for the incompressible Navier-Stokes equations on staggered adaptive Cartesian grids}

\author[UNITN]{Francesco Fambri\corref{cor1}}
\ead{francesco.fambri@unitn.it}
\author[UNITN]{Michael Dumbser}
\ead{michael.dumbser@unitn.it}

\address[UNITN]{Laboratory of Applied Mathematics, Department of Civil, Environmental and Mechanical Engineering, University of Trento, Via Mesiano 77, I-38123 Trento, Italy}

\cortext[cor1]{Corresponding author}

\begin{abstract}
In this paper a new high order semi-implicit discontinuous Galerkin  method (SI-DG) is presented for the solution of the incompressible Navier-Stokes equations 
on \textit{staggered} space-time \textit{adaptive} Cartesian grids (AMR) in two and three space-dimensions. 
The pressure is written in the form of piecewise polynomials on the main grid, which is dynamically adapted within a \emph{cell-by-cell} AMR framework. According to the time dependent main grid, 
 different face-based spatially \emph{staggered} dual grids are defined for the piece-wise polynomials of the respective velocity components. 
Although the resulting adaptive staggered grids are more complex than classical uniform Cartesian meshes, the numerical scheme can still be written in a rather compact form. 
Thanks to the use of a tensor-product formulation for the definition of the nodal basis in the  $d$-dimensional  space ($d=2,3$), all the discrete operators can be efficiently 
written as a combination of linear one-dimensional operators acting in the $d$ space directions separately. 

Arbitrary high order of accuracy is achieved in space, while a very simple semi-implicit time discretization is obtained via 
an explicit discretization of the nonlinear convective terms, and an implicit discretization of the pressure gradient in the momentum equation 
and of the divergence of the velocity field in the continuity equation.  
The real advantages of the staggered grid arise in the solution of the Schur complement associated with the saddle point problem of the discretized incompressible Navier-Stokes equations, i.e. 
after substituting the discrete momentum equations into the discrete continuity equation. This leads to a linear system for only one unknown, the scalar pressure. Indeed, the 
resulting linear pressure system is shown to be symmetric and positive-definite. 
In order to avoid a quadratic stability condition for the parabolic terms given by the viscous stress tensor, an implicit discretization is also used for the diffusive terms in the momentum
equation. A particular feature of our staggered DG approach is that the viscous stress tensor is discretized on the \emph{dual} mesh. This corresponds to the use of a lifting operator, but on 
the staggered grid. Both linear systems for pressure and velocity are very efficiently 
solved by means of a classical matrix-free conjugate gradient method, for which fast convergence is observed. Moreover, it should be noticed that all test cases shown in this paper have 
been performed without the use of any preconditioner. 
Due to the explicit discretization of the nonlinear convective terms, the final algorithm is stable within the classical CFL-type time step restriction for explicit DG methods.  

The new space-time adaptive staggered DG scheme has been thoroughly verified for polynomial degrees up to $N=9$ for a large set of non-trivial test problems in two and three space 
dimensions, for which analytical, numerical or experimental reference solutions exist. The high order of accuracy of the scheme is shown by means of a numerical convergence 
table obtained by comparing the numerical solution against a smooth analytical solution. 

To the knowledge of the authors, this is the first staggered semi-implicit DG scheme for the incompressible Navier-Stokes equations on space-time adaptive meshes in two and three 
space dimensions. 

\end{abstract}


\end{frontmatter}

 
\newcommand{\CK}{Cauchy-Kovalewski }
\newcommand{\CPU}{\textnormal{CPU}}
\newcommand{\p}{\textnormal{P}}
\newcommand{\PM}{\mathbb{P}_M}
\newcommand{\PNPM}{\mathbb{P}_N\mathbb{P}_M}
\newcommand{\PNM}{\mathbb{P}_N\mathbb{P}_M}
\newcommand{\PMM}{\mathbb{P}_N\mathbb{P}_N}
\newcommand{\PzM}{\mathbb{P}_0\mathbb{P}_M} 
\newcommand{\PMPM}{\mathbb{P}_M\mathbb{P}_M}
\newcommand{\dt}{\frac{\partial}{\partial t}}
\newcommand{\dx}{\frac{\partial}{\partial x}}
\newcommand{\dy}{\frac{\partial}{\partial y}}
\newcommand{\dz}{\frac{\partial}{\partial z}}
\newcommand{\Path}{\mathbf{\Psi}}
\newcommand{\dtau}{\frac{\partial}{\partial \tau}}
\newcommand{\dxi}{\frac{\partial}{\partial \xi}}
\newcommand{\deta}{\frac{\partial}{\partial \eta}}
\newcommand{\dzeta}{\frac{\partial}{\partial \zeta}}
\newcommand{\tens}[1]{\underline{\underline{#1}}}
\newcommand{\K}{\tens{K}}
\newcommand{\M}{\tens{M}}
\newcommand{\Roe}{\textnormal{Roe}}
\newcommand{\HLL}{\textnormal{HLL}}
\newcommand{\U}{\mathcal{U}}
\newcommand{\Fd}{\mathcal{F}}
\newcommand{\Gd}{\mathcal{G}}
\newcommand{\Hd}{\mathcal{H}}
\newcommand{\Sd}{\mathcal{S}}
\newcommand{\tot}{\textnormal{tot}}
\newcommand{\CFL}{\textnormal{CFL}}
\newcommand{\vv}{\vec v^{\,2}}
\newcommand{\for}{\textnormal{for}}
\newcommand{\nn}{n+\frac{1}{2}}
\newcommand{\jj}{j+\frac{1}{2}}
\newcommand{\Qi}{\mathbf{Q}_i^n} 
\newcommand{\Qj}{\mathbf{Q}_j^n} 
\newcommand{\Qjj}{\mathbf{Q}_{j+\frac{1}{2}}^{n+\frac{1}{2}}} 
\newcommand{\Aip}{V_i^+} 
\newcommand{\Aim}{V_i^-} 
\newcommand{\Ajp}{V_j^+} 
\newcommand{\Ajm}{V_j^-} 
\newcommand{\Aom}{V_1^-} 
\newcommand{\Atm}{V_2^-} 
\newcommand{\Ahm}{V_3^-} 
\newcommand{\AT}{\left|T_i\right|}
\newcommand{\halb}{\frac{1}{2}}
\newcommand{\FQi}{\tens{\mathbf{F}}\left(\Qi\right)}
\newcommand{\FQj}{\tens{\mathbf{F}}\left(\Qj\right)}
\newcommand{\FQjj}{\tens{\mathbf{F}}\left(\Qjj\right)}
\newcommand{\nj}{\vec n_j}
\newcommand{\FORCE}{\textnormal{FORCE}}
\newcommand{\GFORCE}{\textnormal{GFORCEN}}
\newcommand{\LF}{\textnormal{LF}'}
\newcommand{\LW}{\textnormal{LW}'}
\newcommand{\WL}{\mathcal{W}_h^-}
\newcommand{\WR}{\mathcal{W}_h^+}
\newcommand{\nur}{\boldsymbol{\nu}^\textbf{r} }
\newcommand{\nuf}{\boldsymbol{\nu}^{\boldsymbol{\phi}} }
\newcommand{\nut}{\boldsymbol{\nu}^{\boldsymbol{\theta}} }
\newcommand{\ar}{\phi_1\rho_1}
\newcommand{\arr}{\phi_2\rho_2}
\newcommand{\ur}{u_1^r}
\newcommand{\uf}{u_1^{\phi}}
\newcommand{\ut}{u_1^{\theta}}
\newcommand{\urr}{u_2^r}
\newcommand{\uff}{u_2^{\phi}}
\newcommand{\utt}{u_2^{\theta}}
\newcommand{\ub}{\textbf{u}_\textbf{1}}
\newcommand{\ubb}{\textbf{u}_\textbf{2}}
\newcommand{\RoeMat}{{\tilde A}_{\Path}^G} 
\newcommand{\be}{\begin{equation}}
\newcommand{\ee}{\end{equation}}
\newcommand{\bdm}{\begin{displaymath}}
\newcommand{\edm}{\end{displaymath}}

\newcommand{\apriori}{\textit{a priori} }
\newcommand{\aposteriori}{\textit{a posteriori} }

\newcommand{\Mass}[1]{\mathbf{M}^{\textsl{#1}} }
\newcommand{\iMass}[1]{\mathbf{M}^{\mathbf{-1}\,\textsl{#1}} }
\newcommand{\barr}[2]{\mathbf{#1}_{\text{\textbf{\textsl{p}}}}^{\textsl{#2}}  }
\newcommand{\tild}[2]{\mathbf{#1}_{\text{\textbf{\textsl{v}}}}^{\textsl{#2}}  }
\newcommand{\operator}[2]{\mathbf{#1}^{\textsl{#2}}  }
\newcommand{\Hoperator}[2]{\mathbb{#1}^{\textsl{#2}}  }
\newcommand{\dof}[1]{\widehat{\mathbf{#1}}}
\newcommand{\ndof}{N_{\text{dof}}}
\newcommand{\SIDG}{\text{SI-DG}}
\newcommand{\err}{\mathfrak{r}}
\newcommand{\Nedg}[1]{N_{{#1}-\text{faces}}}
\newcommand{\Nel}{N_{\text{elem}}}
\newcommand{\real}{{\rm I\!R}}
\newcommand{\nat}{{\rm I\!N}}
\newcommand{\dualk}[1]{*_{(#1)}}
\newcommand{\kapped}{\left._{(k)} \right.}
\newcommand{\ixed}[1]{\left._{(#1)} \right.}
\newcommand{\US}{\text{u.s.}}
\newcommand{\SP}[2]{\left\langle #1 , #2 \right\rangle} 
\newcommand{\unit}{\mathbb{I}} 
\newcommand{\basis}{\mathpzc{B}}
\newcommand{\virg}[1]{``#1''} 
\newcommand{\ov}[1]{|_{#1}}
\newcommand{\emm}{}  

\newcommand{\vb}{\mathbf{v}}
\def\Circlearrowleft{\ensuremath{%
  \rotatebox[origin=c]{180}{$\circlearrowleft$}}}
\def\Circlearrowright{\ensuremath{%
  \rotatebox[origin=c]{180}{$\circlearrowright$}}}
\def\CircleArrowleft{\ensuremath{%
  \reflectbox{\rotatebox[origin=c]{180}{$\circlearrowleft$}}}}
\def\CircleArrowright{\ensuremath{%
  \reflectbox{\rotatebox[origin=c]{180}{$\circlearrowright$}}}}
	
	\newcommand{\cyclic}[3]{\left\{\left.^{#1}\CircleArrowleft\right.^{#2}_{#3}\right\}}

\definecolor{qqqqff}{rgb}{0.,0.,1.}
\definecolor{ffqqqq}{rgb}{1.,0.,0.}
\definecolor{qqqqffb}{rgb}{0.3333333333333333,0.3333333333333333,0.3333333333333333}
\definecolor{cqcqcq}{rgb}{0.7529411764705882,0.7529411764705882,0.7529411764705882} 
\definecolor{ffffff}{rgb}{1.,1.,1.}
\definecolor{color1}{rgb}{1.,1.,1.}
\definecolor{color2}{rgb}{0.8627450980392157,0.0784313725490196,0.23529411764705882}
\definecolor{color3}{rgb}{0.,0.8,0.}
\definecolor{color4}{rgb}{0.,1.,0.}
\definecolor{color5}{rgb}{1.,0.,0.}
\definecolor{color6}{rgb}{0.,0.,1.}

\newcommand{\x}{\mbf{x}}
\newcommand{\mbf}[1]{\mathbf{#1}}			%
\newcommand{\su}[1]{\text{\rotatebox{90}{#1}}}

\newcommand{\boxalign}[2][0.97\textwidth]{
 \par\noindent\tikzstyle{mybox} = [draw=black,inner sep=6pt]
 \begin{center}\begin{tikzpicture}
  \node [mybox] (box){%
   \begin{minipage}{#1}{\vspace{-5mm}#2}\end{minipage}
  };
 \end{tikzpicture}\end{center}
}




\section{Introduction} 
\label{Intro}

In this paper a novel semi-implicit discontinuous Galerkin method for solving the incompressible Navier-Stokes equations is derived within the framework of \emph{staggered adaptive} meshes 
(staggered AMR), in two and three space dimensions. The governing equations for an incompressible fluid read 
\begin{align}
&\frac{\partial  \mathbf{v}}{\partial t} + \nabla \cdot \mathbf{F} + \nabla p = 0, \label{eq:NSmom} \\
&\nabla \cdot \mathbf{v} = 0, \label{eq:NSinc}
\end{align}
where $\vb(\x,t)$ is the velocity vector field, $p(\x,t)$ is the pressure and $\mathbf{F}$ takes into account non-linear convection $\mathbf{F}_c= \vb \otimes \vb$ as well as the viscous 
stress tensor $\boldsymbol{\sigma} =-\nu \nabla \vb$, i.e. 
\begin{align}
\mathbf{F} = \mathbf{F}_c + \boldsymbol{\sigma} = \mathbf{v} \otimes \mathbf{v} - \nu \nabla \mathbf{v}, 
\label{eqn:fluxtensor} 
\end{align}
$\nu$ being the kinematic viscosity of the considered fluid. The governing partial differential equations (\ref{eq:NSmom}-\ref{eq:NSinc}) are of great interest in many different scientific areas 
ranging from geophysical flows like meteorological flows and free-surface hydrodynamics in oceans and rivers, over industrial flows in aerospace and mechanical engineering
to biological flows like those in the human cardiovascular and respiratory system. 

Resolving the smallest spatial scales appearing in the flow within large domains requires a higher-order accurate method with very low numerical diffusion and dispersion errors. 
%
Over the years, finite-difference (FD) and finite-volume  (FV) methods have been widely used for solving many different 
families of partial differential equations. Very high order of accuracy at low computational cost can be easily reached for FD schemes 
whenever regular structured grids are used. 
In contrast to FD methods, FV schemes are particularly suitable for general unstructured meshes, but the expensive and rather cumbersome 
\emph{recovery} or \emph{reconstruction} step that is needed for obtaining high order of accuracy in the FV framework constitutes 
a clear drawback. 
In all cases, higher order formulations of FD and FV methods require large stencils, leading to a deterioration of the parallel 
scalability of the algorithms. In this context, the excellent parallel scalability properties of discontinuous Galerkin (DG)
finite element methods make this class of finite element (FE) schemes well suited for large-scale simulations. 
 
In 1973 Reed and Hill derived a new class of methods that allowed the '\emph{flux to be discontinuous across triangle interfaces}' \cite{reedhill} 
for solving the neutron transport equation.  Around the 90'ies, the DG method has been extended to general nonlinear hyperbolic systems in a series 
of well-known papers by Cockburn and Shu and collaborators, see e.g. \cite{Cockburn1989b,Cockburn1990,CockburnShu98}. A review of DG finite element 
methods can be found in \cite{cockburn_2000_dg,cockburn_2001_rkd}. 
In the DG formulation, arbitrary high order of accuracy in space can be achieved by adding new independent degrees of freedom to the chosen polynomial 
basis. Typically, a stable higher order time discretization was reached by means of the method of lines (MOL) approach in combination with explicit 
TVD Runge Kutta schemes, leading to the well known family of RKDG schemes, although also other explicit time-discretizations are possible, see e.g. 
the family of ADER-DG and Lax-Wendroff DG schemes presented in \cite{dumbser_jsc,QiuDumbserShu}. 
On the other hand, one of the major drawbacks for any \emph{explicit} DG scheme is represented by the severe $\CFL$ stability condition that limits the time-step 
of the simulations to be proportional to $h/(2N+1)$ for hyperbolic PDE or even proportional to $h^2/(2N+1)^2$ for parabolic PDE, where $h$ is the characteristic mesh size and 
$N$ is the degree of the polynomial basis. 
Van der Vegt et al. have extended the DG method to an elegant space-time formalism  \cite{spacetimedg1,spacetimedg2,KlaijVanDerVegt}, providing a fully 
implicit and \emph{unconditionally stable} DG scheme for the compressible Euler and Navier-Stokes equations, and later also for the incompressible Navier-Stokes equations, 
see \cite{Rhebergen2013}. 
The drawback is that a fully implicit time-discretization leads to a highly coupled non-linear system for the complete set of the degrees of freedom of the physical variables 
to be solved at every time-step. In this case large scale simulations in two and three space-dimensions can become computationally very demanding. 
Alternative families of \textit{linearly implicit} time discretizations for DG schemes have been considered recently in \cite{Bassi2015}. 

In their pioneering work \cite{BassiRebay} Bassi and Rebay have presented the first DG scheme for the solution of the compressible Navier-Stokes equations, and shortly after 
Baumann and Oden \cite{BaumannOden1,BaumannOden2} have formulated a DG scheme based on penalty terms for the treatment of convection-diffusion equations. Indeed, 
whenever parabolic (second order) or higher order spatial derivatives appear in the considered PDE, obtaining a DG finite-element formulation is not straightforward,  
see \cite{CockburnShu1998,yan2002}. 
High order DG methods are actually a very active field of the ongoing research and several different formulation for the Navier-Stokes equations have been provided in the meantime, see 
\cite{Bassi2007,MunzDiffusionFlux,stedg2,DumbserNSE,HartmannHouston1,HartmannHouston2,Crivellini2013,KleinKummerOberlack2013} to mention a few.


In this work we adopt a specific technique for circumventing a direct approach for solving the saddle point problem of the incompressible Navier-Stokes equations, which goes  
back to  a  family of very efficient semi-implicit finite difference methods developed by Casulli et al. in the context of simulating hydrostatic and non-hydrostatic gravity-driven free-surface 
flows on staggered grids, see \cite{CasulliCompressible,CasulliWalters,Casulli2009,CasulliVOF}. A theoretical analysis of this approach has been provided in 
\cite{CasulliCattani,BrugnanoCasulli,BrugnanoCasulli2,CasulliZanolli2012}. In this family of methods, exact mass conservation is ensured via a conservative finite-volume formulation of 
the discrete continuity equation and the nonlinear convective terms are discretized \textit{explicitly}, in order to obtain a well-behaved pressure system that is at most mildly 
nonlinear (i.e. with nonlinearities only on the diagonal) and whose linear part is at least symmetric and positive semi-definite. 
The common point of the numerical methods mentioned above consist in the application of the Schur complement for the solution of the discrete saddle point problem that results 
after a semi-implicit discretization of the PDE within a staggered-mesh framework. 
Characterized by a high computational efficiency, these methods have been extended to other problems, e.g. blood flow in the human cardiovascular system 
\cite{CasulliDumbserToro,Blood3D2014}, compressible gas dynamics in compliant tubes \cite{DumbserIbenIoriatti} and the dynamics of compressible fluids 
with general equation of state  \cite{DumbserCasulli2016}. 

The first direct extension of staggered semi-implicit finite volume and finite difference schemes to the DG framework has been derived in \cite{DumbserCasulli2013,TavelliDumbser2014} 
for the shallow water equations on Cartesian and unstructured triangular 
grids. The resulting staggered semi-implicit DG method has furthermore been extended to the incompressible Navier-Stokes equations in two and three space dimensions on uniform Cartesian 
and conforming unstructured simplex meshes, see \cite{FambriDumbser,TavelliDumbser2014b,TavelliDumbser2015,TavelliDumbser2016}. 
Concerning the uniform Cartesian grid case \cite{FambriDumbser}, a rigorous theoretical analysis of the corresponding algebraic systems has been very recently presented 
in \cite{SIDG_analysis2016} by employing the theory of matrix-valued symbols and Generalized Locally Toeplitz (GLT) algebras, see \cite{serra1998,GSz,glt}. 

Actually, there exist several different alternative formulations that combine the stability properties of semi-implicit methods and the high order of accuracy of DG schemes 
within \emph{collocated grids}. Some important examples are provided by Dolejsi et al. \cite{Dolejsi2004,Dolejsi2007,Dolejsi2008} for compressible gas dynamics and 
convection-diffusion equations, as well as the work of \cite{GiraldoRestelli,Tumolo2013} for shallow water systems. Concerning staggered-meshes, relevant research has 
been carried out by Chung et al. in \cite{chung2012staggered,ChungNS} for \emph{edge-based}  staggered meshes, and by Liu et al. for the analysis of a DG finite element 
method based on the alternative \emph{vertex-based} staggering approach, see \cite{Liu2007,Liu2008}. 

In 1984, Berger and collaborators presented the adaptive mesh refinement (AMR) approach for finite difference and finite volume schemes for hyperbolic equations, 
see \cite{Berger-Oliger1984,Berger-Colella1989}. Their version of AMR was written in the form of \emph{nested, logically rectangular and refined meshes}, 
or \emph{patches} and the employed numerical schemes were at most second order accurate.  The first higher order path-based AMR method was provided by 
Baeza and Mulet in \cite{BaezaMulet2006}, using up to fifth order accurate WENO schemes. The so-called \emph{cell-by-cell} AMR approach, which has been adopted
in this paper, has first been introduced by Khokhlov in \cite{Khokhlov98} and was later also extended to high-order ADER-WENO finite volume schemes in 
\cite{AMR3DCL,AMR3DNC} for general conservative and non-conservative hyperbolic systems of PDE. 
Because of their great flexibility and since they directly allow the use of non-conforming meshes, DG methods have already been extensively implemented on adaptive meshes, 
commonly known as \emph{hp-adaptive} DG methods, see in particular \cite{BaumannOden1,BaumannOden2,Houston2002,Houston2000,Houston2002}. 
DG methods with AMR have been successfully extended also to the \emph{unstructured} and the anisotropic mesh case, see respectively \cite{Luo2008,Yu2011} for the Euler equations 
and \cite{Leicht2008} for the compressible Navier-Stokes equations. 
Concerning implicit time discretizations, Kopera and Giraldo \cite{KoperaGiraldo} presented an interesting implicit-explicit (IMEX) DG method on AMR meshes for the compressible 
Euler equations with application to atmospheric flow simulations. For further references see also \cite{Georgoulis2009,Lu2014}. 

In this paper, the family of \textit{spectral} semi-implicit DG methods for the solution of the two and three dimensional Navier-Stokes equations on 
\textit{edge-based staggered} Cartesian grids \cite{FambriDumbser} is extended to staggered grids with adaptive mesh refinement (AMR). 
The main novelty of the paper consists in the development of the first high order DG scheme on \emph{staggered AMR meshes}. 
Similar to the uniform Cartesian case, even within staggered AMR grids, an important achievement in terms of computational efficiency has been obtained by 
succeeding in writing all discrete operators as a combination of simple one-dimensional operators thanks to the use of tensor-products of one-dimensional 
operators. 
The method is tested on a large set of test problems in two and three space dimensions, employing polynomial degrees up to $N=9$. To the knowledge of the authors, 
this is the \emph{first time} that a high order accurate semi-implicit DG scheme is derived on \emph{staggered adaptive grids}. 
 
The rest of the  paper is structured as follows:  Section \ref{SI-DG} is devoted to a description of the presented numerical method; in particular 
Section \ref{sec:AMR} outlines the chosen \emph{cell-by-cell} AMR strategy on the main grid used to represent the discrete pressure; 
Section \ref{sec:sAMR} defines the corresponding \emph{staggered AMR mesh};
Section \ref{sec:pN} defines the solution space for the adopted discrete formulation of the governing equations within our staggered DG framework; 
Sections \ref{sec:AMRDGSI}-\ref{sec:FA} outline the details of our $\SIDG$ discretization procedure; 
Section \ref{sec:tests} is devoted to the validation of the scheme and providing comparisons with available analytical, numerical or experimental 
reference solution for several non-trivial test cases in two and three space-dimensions. 
Finally, Section \ref{sec:conclusions} provides some concluding remarks and an outlook to future work.

\section{Semi-implicit DG schemes on staggered AMR grids} 
\label{SI-DG}
In the following,  wherever not specified, the equations and illustrations are given for the three-dimensional case, i.e. $d=3$.

		\subsection{A cell-by-cell AMR framework}  
		\label{sec:AMR}
\begin{figure}
	\centering 
	\vspace{10pt}
	\resizebox{0.45\textwidth}{!}{		\begin{tikzpicture}[line cap=round,line join=round,>=triangle 45,x=0.7cm,y=0.7cm]
\clip(0.,0) rectangle (45.,29.);
\fill[line width=1.6pt,color=ffffff,fill=ffffff,fill opacity=1.0] (0.,9.) -- (0.,0.) -- (9.,0.) -- (9.,9.) -- cycle;
\fill[line width=1.6pt,color=ffffff,fill=ffffff,fill opacity=1.0] (9.,9.) -- (9.,0.) -- (18.,0.) -- (18.,9.) -- cycle;
\fill[line width=1.6pt,color=ffffff,fill=ffffff,fill opacity=1.0] (18.,9.) -- (18.,0.) -- (27.,0.) -- (27.,9.) -- cycle;
\fill[line width=1.6pt,fill=black,fill opacity=0.01] (0.,18.) -- (0.,9.) -- (9.,9.) -- (9.,18.) -- cycle;
\fill[line width=1.6pt,color=ffffff,fill=ffffff,fill opacity=1.0] (18.,18.) -- (18.,9.) -- (27.,9.) -- (27.,18.) -- cycle;
\fill[line width=1.6pt,color=ffffff,fill=ffffff,fill opacity=1.0] (0.,27.) -- (0.,18.) -- (9.,18.) -- (9.,27.) -- cycle;
\fill[line width=1.6pt,color=ffffff,fill=ffffff,fill opacity=1.0] (9.,27.) -- (9.,18.) -- (18.,18.) -- (18.,27.) -- cycle;
\fill[line width=1.6pt,color=ffffff,fill=ffffff,fill opacity=1.0] (18.,27.) -- (18.,18.) -- (27.,18.) -- (27.,27.) -- cycle;
\fill[line width=1.6pt,color=ffffff,fill=ffffff,fill opacity=1.0] (9.,12.) -- (9.,9.) -- (12.,9.) -- (12.,12.) -- cycle;
\fill[line width=1.6pt,color=ffffff,fill=ffffff,fill opacity=1.0] (12.,12.) -- (12.,9.) -- (15.,9.) -- (15.,12.) -- cycle;
\fill[line width=1.6pt,color=ffffff,fill=ffffff,fill opacity=1.0] (15.,12.) -- (15.,9.) -- (18.,9.) -- (18.,12.) -- cycle;
\fill[line width=1.6pt,dash pattern=on 10pt off 10pt,color=ffffff,fill=ffffff,fill opacity=1.0] (9.,15.) -- (9.,12.) -- (12.,12.) -- (12.,15.) -- cycle;
\fill[line width=1.6pt,color=ffffff,fill=ffffff,fill opacity=1.0] (9.,18.) -- (9.,15.) -- (12.,15.) -- (12.,18.) -- cycle;
\fill[line width=0.pt,dash pattern=on 10pt off 10pt,color=ffffff,fill=ffffff,fill opacity=1.0] (12.,18.) -- (12.,15.) -- (15.,15.) -- (15.,18.) -- cycle;
\fill[line width=1.6pt,dash pattern=on 10pt off 10pt,color=ffffff,fill=ffffff,fill opacity=1.0] (15.,18.) -- (15.,15.) -- (18.,15.) -- (18.,18.) -- cycle;
\fill[line width=1.6pt,color=ffffff,fill=ffffff,fill opacity=1.0] (15.,15.) -- (15.,12.) -- (18.,12.) -- (18.,15.) -- cycle;
\draw [line width=1.6pt] (0.,9.)-- (0.,0.);
\draw [line width=1.6pt] (0.,0.)-- (9.,0.);
\draw [line width=1.6pt] (9.,9.)-- (0.,9.);
\draw [line width=1.6pt] (9.,9.)-- (9.,0.);
\draw [line width=1.6pt] (9.,0.)-- (18.,0.);
\draw [line width=1.6pt] (18.,0.)-- (18.,9.);
\draw [line width=1.6pt] (18.,9.)-- (9.,9.);
\draw [line width=1.6pt] (18.,9.)-- (18.,0.);
\draw [line width=1.6pt] (18.,0.)-- (27.,0.);
\draw [line width=1.6pt] (27.,0.)-- (27.,9.);
\draw [line width=1.6pt] (27.,9.)-- (18.,9.);
\draw [line width=1.6pt] (0.,18.)-- (0.,9.);
\draw [line width=1.6pt] (9.,9.)-- (9.,18.);
\draw [line width=1.6pt] (9.,18.)-- (0.,18.);
\draw [line width=1.6pt,dash pattern=on 10pt off 10pt] (9.,18.)-- (9.,9.);
\draw [line width=1.6pt,dash pattern=on 10pt off 10pt] (9.,9.)-- (18.,9.);
\draw [line width=1.6pt,dash pattern=on 10pt off 10pt] (18.,9.)-- (18.,18.);
\draw [line width=1.6pt,dash pattern=on 10pt off 10pt,color=qqqqffb] (18.,18.)-- (9.,18.);
\draw [line width=1.6pt] (18.,18.)-- (18.,9.);
\draw [line width=1.6pt] (18.,9.)-- (27.,9.);
\draw [line width=1.6pt] (27.,9.)-- (27.,18.);
\draw [line width=1.6pt] (27.,18.)-- (18.,18.);
\draw [line width=1.6pt] (0.,27.)-- (0.,18.);
\draw [line width=1.6pt] (0.,18.)-- (9.,18.);
\draw [line width=1.6pt] (9.,18.)-- (9.,27.);
\draw [line width=1.6pt] (9.,27.)-- (0.,27.);
\draw [line width=1.6pt] (9.,27.)-- (9.,18.);
\draw [line width=1.6pt] (9.,18.)-- (18.,18.);
\draw [line width=1.6pt] (18.,18.)-- (18.,27.);
\draw [line width=1.6pt] (18.,27.)-- (9.,27.);
\draw [line width=1.6pt] (18.,27.)-- (18.,18.);
\draw [line width=1.6pt] (18.,18.)-- (27.,18.);
\draw [line width=1.6pt] (27.,18.)-- (27.,27.);
\draw [line width=1.6pt] (27.,27.)-- (18.,27.);
\draw [line width=1.6pt] (9.,12.)-- (9.,9.);
\draw [line width=1.6pt] (9.,9.)-- (12.,9.);
\draw [line width=1.6pt,dash pattern=on 10pt off 10pt] (12.,9.)-- (12.,12.);
\draw [line width=1.6pt] (12.,9.)-- (15.,9.);
\draw [line width=1.6pt,dash pattern=on 10pt off 10pt] (15.,9.)-- (15.,12.);
\draw [line width=1.6pt,dash pattern=on 10pt off 10pt] (15.,12.)-- (12.,12.);
\draw [line width=1.6pt] (15.,9.)-- (18.,9.);
\draw [line width=1.6pt] (18.,9.)-- (18.,12.);
\draw [line width=1.6pt] (9.,15.)-- (9.,12.);
\draw [line width=1.6pt,dash pattern=on 10pt off 10pt] (12.,15.)-- (9.,15.);
\draw [line width=1.6pt,dash pattern=on 10pt off 10pt] (15.,15.)-- (15.,12.);
\draw [line width=1.6pt,dash pattern=on 10pt off 10pt] (15.,12.)-- (18.,12.);
\draw [line width=1.6pt] (18.,12.)-- (18.,15.);
\draw [line width=1.6pt,dash pattern=on 10pt off 10pt] (18.,15.)-- (15.,15.);
\draw [line width=1.6pt] (9.,18.)-- (9.,15.);
\draw [line width=1.6pt,dash pattern=on 10pt off 10pt] (12.,15.)-- (12.,18.);
\draw [line width=1.6pt] (12.,18.)-- (9.,18.);
\draw [line width=1.6pt,dash pattern=on 10pt off 10pt] (15.,18.)-- (15.,15.);
\draw [line width=1.6pt,dash pattern=on 10pt off 10pt] (18.,15.)-- (18.,18.);
\draw [line width=1.6pt,dash pattern=on 10pt off 10pt] (18.,18.)-- (15.,18.);
\draw [line width=1.6pt] (18.,12.)-- (18.,15.);
\draw [line width=1.6pt] (21.,18.)-- (21.,9.);
\draw [line width=1.6pt] (24.,18.)-- (24.,9.);
\draw [line width=1.6pt,dash pattern=on 10pt off 10pt] (12.,15.)-- (12.,12.);
\draw [line width=1.6pt,dash pattern=on 10pt off 10pt] (12.,12.)-- (9.,12.);
\draw [line width=1.6pt] (27.,12.)-- (18.,12.);
\draw [line width=1.6pt] (18.,15.)-- (27.,15.);
\draw [line width=1.6pt,dash pattern=on 10pt off 10pt] (12.,15.)-- (15.,15.);
\draw [line width=1.6pt] (29.93225667479689,24.32531136108621)-- (35.68980308772662,24.32531136108621);
\draw (37.11445614620858,24.51363700222249) node[anchor=north west] {\Huge active mesh};
\draw [line width=1.6pt,dash pattern=on 10pt off 10pt] (30.021520960268656,22.40612922344297)-- (35.77906737319837,22.40612922344297);
\draw (37.170594049241,22.604948293150024) node[anchor=north west] {\Huge  virtual mesh};
\draw (9.88757317548866,29.17308296848763) node[anchor=north west] {\Huge  $\text{ref. level: }  \ell $};
\draw (18.58894814551255,20.191018455205434) node[anchor=north west] {\Huge  $\text{ref. level: }  \ell + 1$};
\draw (30,5) node[anchor=north west] {\resizebox{0.15\textwidth}{!}{\Huge  $\Omega_h$ }};
 
\end{tikzpicture} }\vline\; 
		\includegraphics[width=0.4\textwidth]{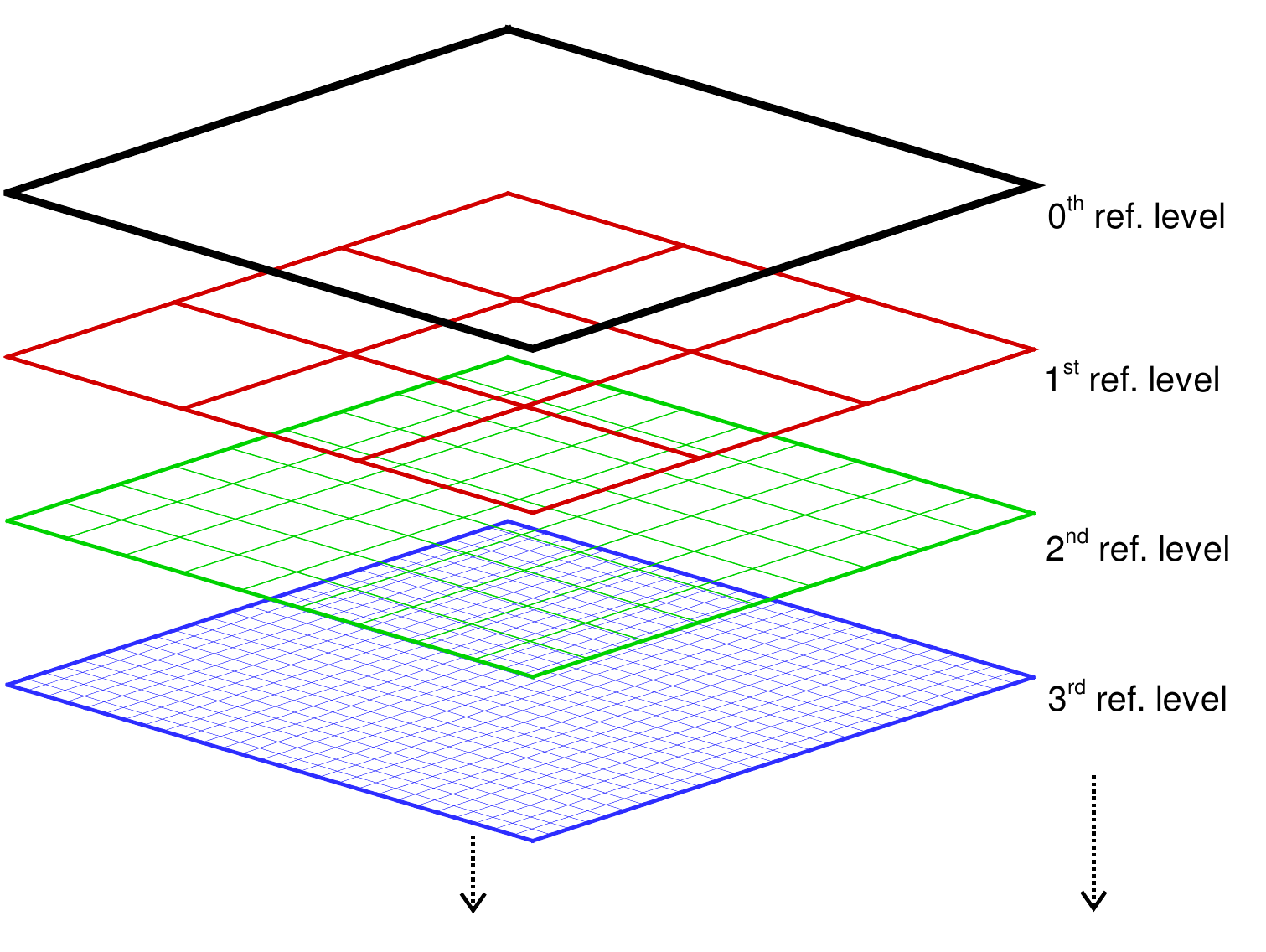}\\
		\vspace{10pt}
	\caption{At the left, a simple sketch of the AMR grid $\Omega_h$ with one single refinement level is shown. At the right, the \emph{tree-structure} of the refinement levels $\ell=0$,$1$, $\ldots$, $\ell_{\text{max}}$  for a single element at the coarsest level $T_i\in\Omega_h^0$ is shown. 
(See colored version on-line)}
	\label{fig:AMR}
\end{figure}
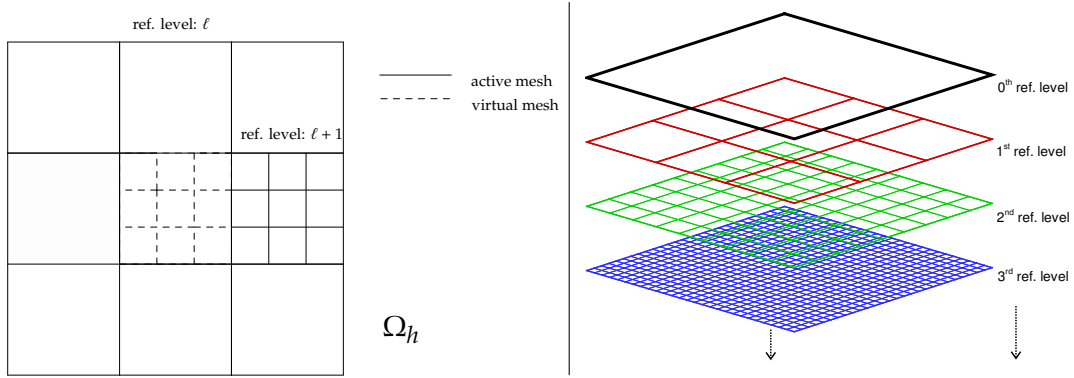

		   There  are mainly two strategies for defining an automatically refined grid, namely the so called \emph{patched} AMR and the \emph{cell by cell} method. 
			In this paper the second strategy is chosen, leading to a tree data structure by defining up to $\ell_{\text{max}}$ refinement levels 
			$\Omega_h^{\ell}$, $\ell=1,2,\ldots,\ell_{\text{max}}$, the coarser grid $\Omega^0_h$ given, covering the entire spatial domain $\Omega$, i.e. 
					\begin{align}
					 \Omega \equiv \bigcup_{T_i \in \Omega_h^\ell} T_i, \;\;  \forall \ell=0,1,\ldots,\ell_{\text{max}}.
					\end{align}
	An integer number  $\err$ for the spatial refinement ratio, or \emph{refinement factor}, between adjacent refinement levels drives the refinement scales. Figure \ref{fig:AMR} gives an illustrative sketch of the resulting adaptive mesh. In order to refine the grid only whenever and wherever necessary, a \emph{refinement-estimator function} ${\chi}_{\emm}$ can be chosen to be a function of the space-derivatives of a given indicator function $\Phi(\x,t)$, e.g. in the form of 
\begin{equation} \label{eq:chi_m}
	{\chi}_{\emm} (\Phi) =  \sqrt{ \frac{\sum_{k,l}{\left(\left. \partial ^2 \Phi \middle/ \partial x_k \partial x_l \right. \right)^2}}{ \sum_{k,l}{ \left[ \left. \Big( \left| \left. \partial \Phi \middle/ \partial x_k\right.\right|_{i+1} + \left| \left.\partial \Phi \middle/ \partial x_k \right. \right|_i \Big) \middle/ \Delta x_l \right. + \epsilon \left|\frac{\partial^2 }{\partial x_k\partial x_l} \right|\left| \Phi \right|    \right]^2}}  }\,,
\end{equation}
which considers up to the second space derivatives. 
The indicator function $\Phi$ can be any physical quantity of interest, e.g. the pressure, the vorticity, the kinetic energy or any other function of the flow quantities. 
$\chi$ will be evaluated periodically in time, according to the time-scales of the physical problem. Then, the space-elements $T_i$ will be \emph{refined} or \emph{recoarsened} 
whenever the prescribed upper and lower threshold values $\chi_{\text{ref}}$ and $\chi_{\text{rec}}$, respectively, are exceeded.  
The resulting \emph{active} mesh $\Omega_h$ is the set of space-elements belonging to the refinement levels $\Omega_h^{\ell}$, $\ell=1,2,\ldots,\ell_{\text{max}}$, such that $\Omega$ is spanned and the non-overlapping property is satisfied, i.e.\footnote{ $\circ$ denotes the \emph{interior} operator, i.e.  with $T^{\circ} = T \setminus \partial T$.}
\begin{align}
\Omega_h = \left\{  T_i \; \Big| \; T_i \in  \bigcup_{\ell=1}^{\ell_{\text{max}}} \Omega_h^{\ell}, \;\;\text{and}\;\; \Omega  = \bigcup_{i} T_i ,  \;\;\text{and}\;\;  \varnothing = \bigcup_{T_i\neq T_j } \left(  T^{\circ}_{i} \cap T^{\circ}_{j} \right), \;\;\text{with}\;\; i,j=1,2,\ldots,\Nel \; \right\}.\label{def:Omega_h}
\end{align} 
For practical purposes, it becomes useful to define the '\emph{status} $\beta$' of the complete set of spatial elements
\begin{align}
 &\forall T_i\in \bigcup_{\ell} \Omega_h^{\ell} & \beta_i = \left\{\begin{array}{rclc} 
-1, & & \text{for the so called \emph{virtual parent cells}, i.e.} &\exists T_j \in \Omega_h \big| T_j \supset T_i \\
0, & & \text{for \emph{active elements}, i.e.} & T_i \in \Omega_h  \\ 
1, & & \text{for the so called \emph{virtual children}, i.e.} & \exists T_j\in \Omega_h \big| T_j \subset T_i 
\end{array}\right..
\end{align}
In order to make things simpler, in this work, two \emph{neighbor and active} elements $T_i$ and $T_j$ are allowed to belong only to the same or to an adjacent refinement level, i.e. if $T_i\in \Omega_h^{\ell}$ then  $T_j \in \Omega_h^{\ell-1\leq \tilde{\ell} \leq \ell+1}$, or in a simpler notation $| \ell(T_j)-\ell(T_i)| \leq 1$, where $\ell(T)$ is the refinement level of a general space-element $T$. 

%

		\subsection{The $d$ spatially staggered AMR meshes}  
		\label{sec:sAMR}
		   
			Chosen a main Cartesian mesh $\Omega_h$, there are basically two different strategies for choosing a staggered Cartesian dual mesh $\Omega^*_h$: the 'node-based' and the 'edge-based'  staggering, i.e.  the \emph{B-}  and the \emph{C-grid}, respectively, according to the nomenclature of Arakawa \& Lamb \cite{Arakawa}.  In this work  the \emph{edge}-based staggering  has been selected to be the optimal one in terms of numerical efficiency for the resolution of the resulting system of discrete equations. A $d$-dimensional uniform Cartesian element $T_i \in \Omega_h$ is defined as the Cartesian product 

\begin{equation*}
T_i = \prod _{s=x,y,z} \Delta s(\x_i) ,\;\;\; \forall\; i=1,2,\ldots,\Nel, 
\end{equation*} 
   where $\Delta s(\x_i)$ are the spatial discretization-steps, centered in $\x_i =(x_i,y_i,z_i)$, i.e. the barycenter of $T_i$. Then, in order to treat hanging-nodes in the chosen AMR-mesh, one can define
$\Nedg{i}^{\kapped}$ being the total number of \emph{faces} (or \emph{edges} in two space dimensions) of $T_i$ that are oriented in the k-th space direction, i.e. $\Gamma_{i,l}^{\kapped}$ with $l=1,2,\ldots,\Nedg{i}^{\kapped}$, having  $\Nedg{i}^{\kapped}=2$ (backward and forward) for uniform grids. Then, the complete set for the faces of $T_i$ will be $\Sigma_i=\left\{\Gamma_{i,1}, \Gamma_{i,2}, \ldots, \Gamma_{i,\Nedg{i}} \right\}$, where $\Nedg{i}$ is the total number of faces of $T_i$ with
	\begin{align*}
	& \Nedg{i} =\sum_{k=1}^d \Nedg{i}^{\kapped} 
	&  \partial T_i   \equiv \Sigma_i  = \bigcup_{j=1}^{\Nedg{i}}\Gamma_{i,j} 
	\end{align*}
The resulting edge-based Cartesian \emph{staggered} (or '\emph{dual}', or '\emph{starred}')  elements are identified by $T^*_{i,j}$,   $j=1$, $2$, $\ldots$, $\Nedg{i}$ and $i =1$, $2$, $\ldots$, $\Nel$.
 It becomes useful to distinguish up to $d$ sets of non-overlapping Cartesian staggered-elements  $T^{\dualk{k}}_{i,j}$ $j=1$, $2$, $\ldots$, $\Nedg{i}^{\kapped}$, and $i =1$, $2$, $\ldots$, $\Nel$, 
	referring to the $k$-th space-direction of the staggering, $k=x$, $y$, $z$.
	Since every internal face $\Gamma$ is shared by two distinct \emph{neighbor} elements $T_{r(\Gamma)}$ and $T_{l(\Gamma)}$, the chosen two-index notation is \emph{surjective}, i.e. exists $j_1$ and $j_2$ such that $T_{l(\Gamma),j_1}^{*} \equiv T_{r(\Gamma),j_2}^{*}$, or equivalently $\Gamma_{l(\Gamma),j_1} \equiv \Gamma_{r(\Gamma),j_2} \equiv \Gamma$.
	Here, $r(\Gamma)$ and $l(\Gamma)$ are defined to be the integer indexes for the right and the left space element, with respect to the oriented face $\Gamma$, which orientation is well definite and no sign-function is needed, because the mesh is Cartesian.
	
	Then, in the aim of simplicity, a one-index \emph{injective} notation can be obtained after an adequate  surjective reordering  map $\rho$
	\begin{align*}
	&  
	\begin{array}{rcccll}
	   \nat \times \nat \supset & A &\xrightarrow{\hspace{0.5cm}\rho\hspace{0.5cm}} & B & = \left\{1,2, \ldots, N_{\text{faces}} \right\}\subset \nat &  \\ 
		& \su{$\in$} & & \su{$\in$}  & &\\
				&  (i,j)&  \xrightarrow{\hspace{0.5cm}\hspace{0.5cm}}&  m &  = \rho(i,j) ,& \forall \;i=1,2,\ldots,\Nel;\;j=1,2,\ldots,\Nedg{i};
		\end{array} 
	\end{align*} 
	so that exists  $m=\rho(l(\Gamma),j_1) =\rho(r(\Gamma),j_2)$ and holds
\begin{equation*} 
T_{l(\Gamma),j_1}^{*} \equiv T_{r(\Gamma),j_2}^{*}  \equiv   T_m^{*} = \prod _{s=x,y,z} \Delta s(\x_m^{*}) ,\;\;\;\forall\; m=1,2,\ldots,N_{\text{faces}} 
\end{equation*} 
	where  $\Delta s (\x_m^{*})$ are the spatial discretization-steps centered in $\x_m^{*} $, i.e.  the barycenter of face $\Gamma_m$; $N_{\text{faces}}$ is the total number of faces of $\Omega_h$. Then, it follows the injectivity property
	\begin{align*}
  T^*_{m_1} \neq T^*_{m_2},\;\;\; \forall m_1\neq m_2, \;\;\; m_1 , m_2\in \left\{1,2, \ldots, N_{\text{faces}} \right\}
	\end{align*}
	or alternatively 
	\begin{align*}
  T^{\dualk{k}}_{m_1} \neq T^{\dualk{k}}_{m_2},\;\;\; \forall m_1\neq m_2, \;\;\; m_1 , m_2\in \left\{1,2, \ldots, N^{\kapped}_{\text{faces}} \right\},\;\;\; k=x,y,z 
	\end{align*}
	$N^{\kapped}_{\text{faces}} $ being the total number of faces in the $k$-th space direction in the computational domain. In the following, the two-index notation will be used only whenever strictly necessary. 
	%
In this notation, one can distinguish the $d+1$ (the main $\Omega_h$ and the $d$ dual meshes $\Omega^{\dualk{k}}_h$, $k=x$, $y$or $z$ if $d=3$) spatially staggered non-overlapping meshes with the property
\begin{align}
&\Omega = \bigcup_i T_i =  \bigcup_i T^{\dualk{k}}_{i}, & \label{def:stagg1}\\   
&  \varnothing = \bigcup_{i\neq j} \left(  T^{\circ}_{i} \cap T^{\circ}_{j} \right) =  \bigcup_{i\neq l}  \left( T^{\dualk{k} \circ}_{i}\cap  T^{\dualk{k} \circ}_{j} \right), &  T \in \Omega_h,\; T^{\dualk{k}}\in \Omega_h^{\dualk{k}} \;k=x,y,z, 
\label{def:stagg2}
\end{align}

\begin{figure} 
\centering 
\begin{tabular}{cc}
			\resizebox{0.45\textwidth}{!}{\input{StaggEl2d}}& 
			\includegraphics[width=0.45\textwidth]{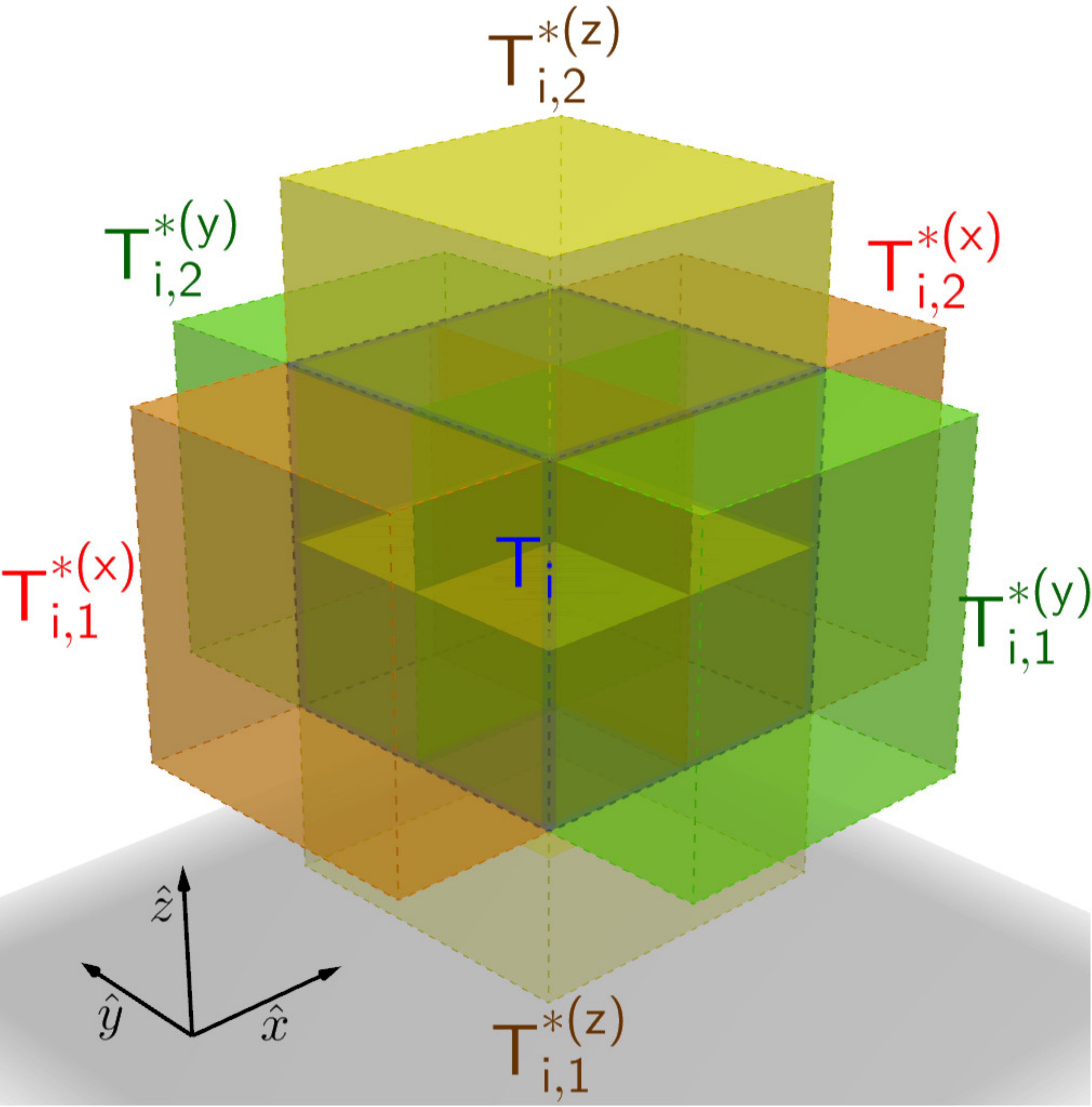}  
\end{tabular} 			
\caption{Mesh-staggering for the two dimensional  (left) and for the three-dimensional  (right) purely-Cartesian case, i.e. in our notation  $\ell_{\text{max}}=0$ and consequently $\Omega_h \equiv \Omega_h^{\ell=0}$ and $\Nedg{i}^{\kapped} = 2$, with $k=x,y,z$. The staggered elements $T^{\dualk{k}} \in \Omega_h^{\dualk{k}}$ corresponding to the faces of $T_i \in \Omega_h$ in the $x$-th, $y$-th (for d=2) and $z$-th (for d=3) direction are shown and highlighted in red, green and ocher, respectively. For the $2d$ case, the nearest neighbor elements of $T_i$ are shown and highlighted in blue. In our notation holds the following correspondence $(T^{\dualk{x}}_{i,1} , T^{\dualk{x}}_{i,2} , T^{\dualk{y}}_{i,1} , T^{\dualk{y}}_{i,2}, T^{\dualk{z}}_{i,1} , T^{\dualk{z}}_{i,2}) = (T^{*}_{i,1} , T^{*}_{i,2}  ,T^{*}_{i,3} , T^{*}_{i,4},T^{*}_{i,5} , T^{*}_{i,6})$ for the $d=3$ case. (See colored version online)}\label{fig:staggering}
\end{figure}

Figure \ref{fig:staggering} depicts the staggered meshes $\Omega_h^{\dualk{k}}$ next to the main grid $\Omega_h$ in two and three space-dimensions in the purely Cartesian case, i.e. $\ell_{\text{max}}=0$ and $\Nedg{i}^{\kapped} = 2$. 

\begin{figure} \centering 
    \includegraphics[width=0.7\textwidth]{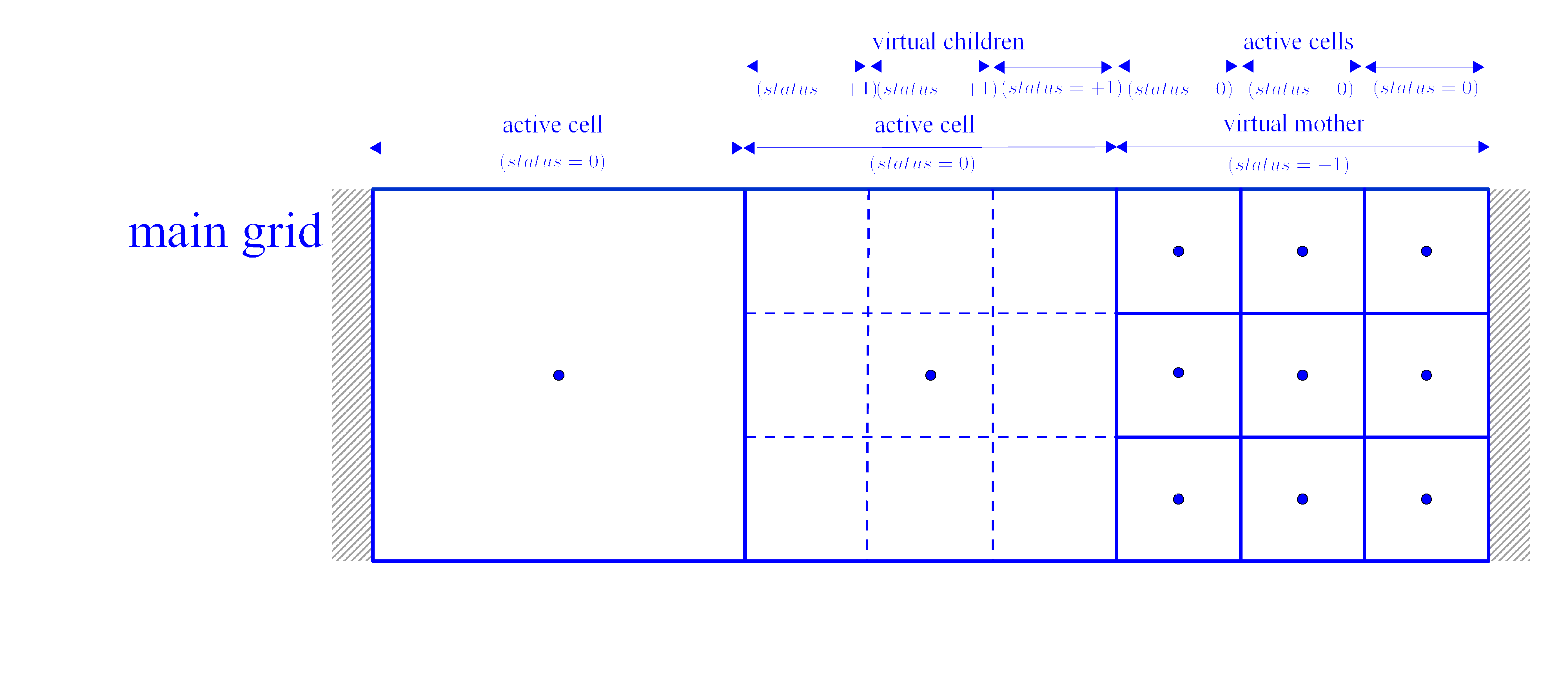} 
    \includegraphics[width=0.8\textwidth]{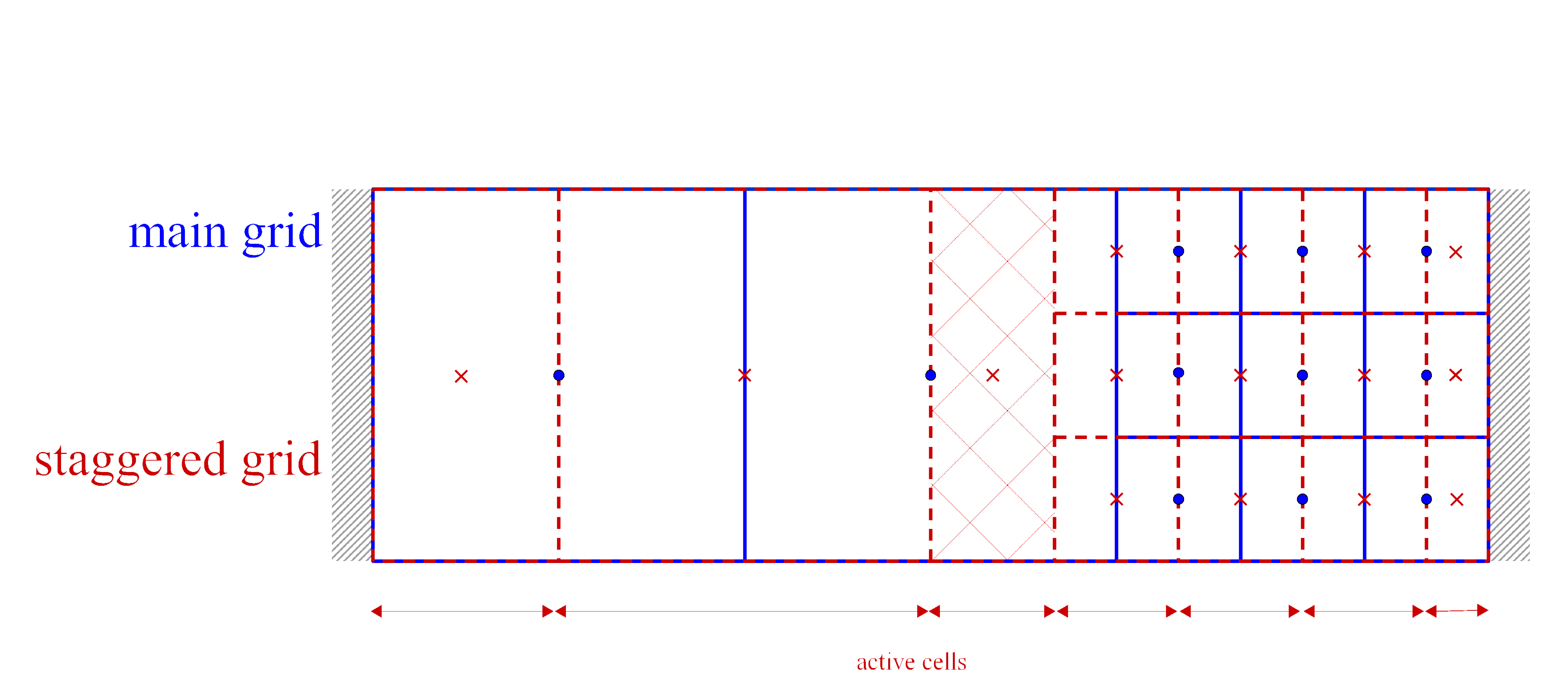} 
			\caption{Schematic two dimensional view of the adopted cell by cell AMR framework (top), and  the chosen spatially staggered  and adaptively refined grid, e.g., in the in $x$ direction $\Omega_h^{\dualk{x}}$ ( highlighted in red at the bottom). In particular, the \emph{anomalous} dual element, over which the pressure is continuous, is highlighted by a checkered texture (See colored version online)}\label{fig:sAMR}
\end{figure}

Whenever AMR is considered the starred notation for the refinement meshes is used, i.e. $\Omega_h^{\dualk{k}}$ for the dual-staggered active grid and $\Omega_h^{\ell \dualk{k}}$ for the dual-staggered grid of refinement level $\ell$. Furthermore, for $\Omega_h^{\dualk{k}}$ a simple trick becomes necessary in order to satisfy condition (\ref{def:stagg1}). Indeed, if uniform Cartesian elements are used for the dual elements $T_m^*$, centered in $\x^*_m$, then some \emph{new kind of elements}, characterized by a different mesh size ratio, need to be defined. An $\US$ element $T_u^*$ appears in the dual computational domain $\Omega_h^*$ wherever two \emph{neighbor} and \emph{active} elements $T_i$ and $T_j$ belong  to adjacent refinement levels, i.e. $T_i \in \Omega_h^{ell}$,  $T_j \in \Omega_h^{ell+1}$  
with $\beta_i = \beta_j=0$.
Figure \ref{fig:sAMR} gives an illustration of the resulting staggered auto-adapted grid $\Omega_h^{\dualk{x}}$  next to $\Omega_h$. The peculiarity of such \emph{unusual staggered} ($\US$) elements will arise in the definition of the discrete solution.  For the sake of simplicity, $N_{\text{faces}}$ will take account also of the number of $\US$  dual elements belonging to $\Omega_h^*$. An auxiliary value for labeling the refinement level of the $\US$ elements $T^*_{u}$ is given by the fractional intermediate index $\ell(T^*_{u}) = (\ell(T_i) + \ell(T_j))/2$. Finally, referring to an $\US$ element $T^*_{u}$, the corresponding face $\Gamma_{u}$ should not exist, but is can be imposed to be equal to the void set $\varnothing$. Similarly the right and left element of $\Gamma_u$ are chosen to be $l(u) = r(u) \equiv i$, where for  $T_i$  holds $T^*_{u} \subset T_i$.

Similarly, also the faces of the dual elements, i.e. the \emph{dual faces} $\Gamma_i^{\dualk{k}} \in \partial T_{m}^{\dualk{k}}$, can be defined, and the respective \emph{dual right} $r^*(i)$ and \emph{dual left}  $l^*(i)$ elements can be identified. These set of dual faces are slightly sophisticated and are used rarely in the text. Then, in order to simplify the notation, we define only the dual edges whose points lie within a given space element, i.e. the faces $\Gamma_i^{\dualk{k}} \in \partial T_{m}^{\dualk{k}}$ for $k=x$, $y$, $z$ such that $\Gamma_i^{\dualk{k}} \subset \partial T_{i}$
Once the main and staggered adaptive meshes are fully defined, in the following section the space of discrete solutions of our $\SIDG$-$\p_N$ method is outlined.


		\subsection{The space of piecewise polynomials  $\p_N$}
		\label{sec:pN}
		
In the DG framework, the discrete solution is defined in the space of piecewise polynomials of maximum degree $N$ with discontinuities along the faces of the elements of the computational mesh. 
Then, given a reference mesh $\tilde{\Omega}$, the corresponding vector space of piecewise polynomials of maximum degree $N \in \nat^+_0$ 
is called $\p_N(\tilde{\Omega})$. In particular, in the staggered-DG formulation with adaptive mesh refinement, we look for piecewise polynomials $P(\x)$, $\x\in\tilde{\Omega}$  where $\tilde{\Omega}$ is one of the considered meshes, i.e.  $\tilde{\Omega}\in \{ \Omega_h, \Omega_h^{\dualk{k}} \}_{k=x,y,z}$. 
Similar to \cite{FambriDumbser,DumbserCasulli2013}, the discrete solution in every single space element $\tilde{T}_i 
\in \tilde{\Omega}$ is written according to the same but shifted \emph{one-dimensional} polynomial basis $\basis_{N}=\left\{ \varphi_l(x) \right\}_{l=0,N}$ along its own control volume $x \in \Delta \tilde{s}_k (\tilde{\x}_i)$, for \emph{each} spatial dimension $k=x,y,z$.

\begin{figure}
	\centering 
\resizebox{0.8\textwidth}{!}{\input{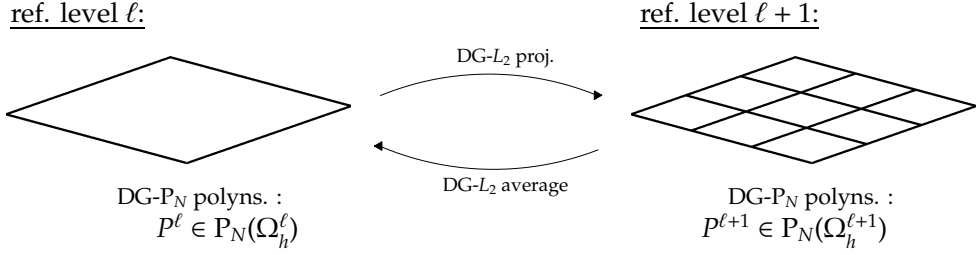}}
	\caption{Simple sketch of the polynomial mapping between two adjacent refinement levels $\ell$ and $\ell+1$, i.e. the $L_2$ projection and average.}
	\label{fig:mapping}
\end{figure}

In order to transform the numerical solution from a coarser level $\ell$  
to the finer $\ell+1$ 
and \emph{vice versa} (see Figure \ref{fig:mapping}) suitable $L_2$ projection $\mathcal{P}$ and average $\mathcal{A}$ operators are defined in the form
\begin{align} 
&\p_N(\Omega_h^{\ell+1}) \ni P_{c}^{\ell+1} =  \mathcal{P}\left(P_{i}^{\ell}\right),& \;&   \SP{\omega^{\ell+1}}{P_{c}^{\ell+1}}_{T_{c}} = \SP{\omega^{\ell+1}}{P_{i}^{\ell}}_{T_{c}} , & \forall T_{c} \subset T_{i}, \;\;  T_c \in \Omega_h^{\ell+1},\;T_i \in \Omega_h^{\ell} \\ 
& \p_N(\Omega_h^{\ell}) \ni  P_{i}^{\ell} =  \mathcal{A}\left(\{P_{c}^{\ell+1}\}\right),& \;& \SP{\omega^{\ell}}{ P_{i}^{\ell}}_{T_{i}} = \sum \limits_{T_{c}\subset T_{i} }  \SP{\omega^{\ell} }{ P_{c}^{\ell+1}}_{T_{c} }
, &   T_c \in \Omega_h^{\ell+1},\;T_i \in \Omega_h^{\ell}   
\end{align}
where  $\SP{\cdot}{\cdot}$ is the scalar product
\begin{align}
 & \SP{f}{g}_{T_i} = \int \limits_{T_i}  \mathit{f}\, \mathit{g} d \mathbf{x} & \mathit{f}, \mathit{g} \in L^2 
\end{align} 
and
\begin{align*}
 & \omega_{\ov{T_i}}(s) = \varphi(\xi_1) \otimes \varphi(\xi_2) \otimes \varphi(\xi_3) ,  & \text{with} & \qquad s = \x_i + \sum_{k}\left(\xi_k -\frac{1}{2}\right) \Delta x_k^{\ell(T_i)}(\x_i) \hat{x}^{\kapped}. 
& 0 \leq \xi_1,\xi_2,\xi_3 \leq 1,  \label{eq:omega}
\end{align*}
defines the \emph{shifted} basis polynomials, centered in the barycenter $\x_i$ of $T_i\in\Omega_h^{\ell}$ and
written in tensorial form.
If $\varphi$ is the vector of the basis elements $(\varphi_1,\varphi_2, \ldots, \varphi_{N+1}) $, then $\omega_{\ov{T_i}}$ is the  vector of the \emph{shifted} basis elements  $(\omega_1,\omega_2, \ldots, \omega_{(N+1)^d})_{\ov{T_i}} $ generating the three dimensional basis 
$\basis_{N,i}^{d} \equiv \basis_{N,i} \otimes \basis_{N,i} \otimes \basis_{N,i}$, for any $T_i\in\Omega_h^{\ell}$, for any piecewise polynomial $P^{\ell}\in \p_N(T_i)$. 
Then $P$ can be expanded along the basis in the form
\begin{align}
P(x)  = \sum_l \frac{\SP{\omega_l(x)}{P(x)}_{T_i}}{||\omega_l |\ov{T_i}^2} \omega_l \equiv  \sum_l P_l \omega_l  \equiv  \mathbf{P}\cdot \omega(x)    
\end{align}
where $\mathbf{P}$ is called to be the vector of the \emph{degrees of freedom} of $P$, $||\cdot||$ is the norm related to the chosen scalar product.
 Moreover, if one chooses $\basis_{N}^{d}$ to be 
 a \emph{nodal} basis 
then the following equivalence holds: 
$
P_l   \equiv   P(x^{\text{GL}}_l)
$
with a large gain in computational effort, specially at higher space-dimensions when three-dimensional coefficient matrices can be written in a tensorial combination of one dimensional matrices. 
In our scheme, the basis functions are given by the Lagrange interpolation polynomials passing through the Gauss-Legendre quadrature nodes. 

In this work, although each physical variable will be defined referring to \emph{a proper and given computational mesh} (main or dual), it becomes useful defining the transformation rules for projecting 
from one vector space $P \in P_N(\Omega_h)$ to a corresponding staggered one $P^*\in P_N(\Omega^{\dualk{k}}_h)$, i.e. 
 satisfying  the classical $L_2$ projection equations, i.e.
\begin{align}
&P^*= \pi^{\dualk{k}}(P): &   \SP{ \psi}{P^*}_{T^{\dualk{k}}_m}   \;\; =  \;\; \SP{ \psi}{ P}_{T^{\dualk{k}}_m}&  & \forall T_m^{\dualk{k}}\in\Omega_h^{\dualk{k}} \label{eq:pistar}\\
&P= \pi(P^*):& \SP{\omega}{P}_{T_i}  \;\; =  \;\; \SP{ \omega}{ P^*}_{T_i} ;&  &\forall T_i\in\Omega_h  \label{eq:pi}
\end{align}
where the vectors of the \emph{shifted} polynomial basis element $\psi$ are given by
\begin{align}
 & \psi_{\ov{T_{i,j}^{\dualk{k}}}}(s) = \varphi(\xi_1) \otimes \varphi(\xi_2) \otimes \varphi(\xi_3) ,  & \text{with} & \qquad s = \x^{\dualk{k}}_m + \sum_{k} \left(\xi_k -\frac{1}{2}\right) \Delta x_k^{\ell(T_{i,j}^{\dualk{k}})}(\x^{\dualk{k}}_m) \hat{x}^{\kapped}, & 0 \leq \xi_1,\xi_2,\xi_3 \leq 1,  \label{eq:psi}
\end{align}
i.e. polynomials $\psi_{\ov{T_{i,j}^{\dualk{k}}}}$  centered in the barycenter $\x^{\dualk{k}}_m$ of $T^{\dualk{k}}_m\in\Omega_h^{\dualk{k}}$ and
written in tensorial form; $k$ labels the direction of the staggering for $\Omega_h^{\dualk{k}}$ ($k=x$, $y$ or $z$).
 
Following \cite{FambriDumbser}, because of the tensorial definition of the polynomial basis, a generic coefficient matrix $\mathbf{Z}$ operating over a generic vector of degrees of freedom $X$, appearing in (\ref{eq:pi2}), can be written as 
\begin{align}
&\operator{Z}{}   \equiv \operator{Z}{xyz}\equiv \operator{Z}{x}\otimes \operator{Z}{y}\otimes\operator{Z}{z},& \dof{x}  \equiv \{\hat{x}_l \}_{l=1,..,(N+1)^d} \equiv \{\hat{x}_{ll'l''}\equiv \hat{x}_{l} \otimes  \hat{x}_{l'} \otimes  \hat{x}_{l''} \}_{l,l',l''=1,..,N+1},
\end{align}
where the components $\operator{Z}{k}$ operate along the respective $k$-th index of 
 the  generic vector of degrees of freedom 
$
\dof{X} 
$
.
Whenever a $k$-th direction is specified, an identity operator is assumed in the other directions, i.e.
\begin{align}
& 
\virg{\;\operator{Z}{x}\;} \equiv \operator{Z}{x} \otimes \unit \otimes \unit, & \virg{\;\operator{Z}{y}\;} \equiv \unit \otimes \operator{Z}{y} \otimes  \unit, & & \virg{\;\operator{Z}{z}\;}\equiv  \unit \otimes\unit \otimes\operator{Z}{z}, \end{align}
\begin{align}
& \virg{\;\operator{Z}{x y}\;} \equiv \operator{Z}{x} \cdot\operator{Z}{y}  \equiv \operator{Z}{x} \otimes \operator{Z}{y}\otimes \unit, & \virg{\;\operator{Z}{x z}\;} \equiv \operator{Z}{x} \cdot\operator{Z}{z}  \equiv \operator{Z}{x} \otimes \unit \otimes \operator{Z}{z}, & & \virg{\;\operator{Z}{y z}\;} \equiv \operator{Z}{y} \cdot\operator{Z}{z}  \equiv \unit  \otimes  \operator{Z}{y} \otimes \operator{Z}{z},\end{align} 
%
%
Then, after evaluating the definite integrals of the known basis functions, one gets 
\begin{align}
&  \left.
\begin{array}{l}
\SP{\omega_l}{P}_{T_i} = \mathbf{M}_{lm\ov{T_i}} \hat{P}_m \\
\SP{\omega_l}{P^*}_{T_i} = \sum\limits_{\Gamma_j^{\kapped}} \mathbf{M}^{\sigma(i,j,k)}_{lm} \hat{P}^*_{m\ov{T_{i,j}^{\kapped}}}  
\end{array} \right\} & \Rightarrow \hat{P}_{l} = \mathbf{M}_{lm\ov{T_i}}^{-1} \sum\limits_{\Gamma_{i,j}^{\kapped}} \mathbf{M}^{\sigma(i,j,k)}_{mn} \hat{P}^*_{n\ov{T_{i,j}^{\kapped}}}  \label{eq:pi2} 
\end{align}
where the summation is intended over the $j$-th faces of $T_i$ in the direction of the staggering, i.e. the $k$-th, backward and forward; $\sigma(i,j,j)$ labels  the mutual position and refinement of $T_i$ and the corresponding dual element $T_{i,j}^{\kapped}$. 
Some new  matrices have been defined to be
\begin{align}
&\mathbf{M}_{\ov{T_i}} =  \left\{ M_{pq\ov{T_i}} \right\}_{p,q=1,..,(N+1)^d} = \left\{\int \limits_{T_i} \,\omega_{p\ov{T_i}}\left(\xi\right) \omega_{q\ov{T_i}}\left(\xi\right)  d\xi \right\}_{p,q=1,..,(N+1)^d} \nonumber \\
& \mathbf{M}^{\sigma(i,j,k)} = \left\{  \mathbf{M}^{\sigma(i,j,k)}_{pq}\right\}_{p,q=1,..,(N+1)^d} = \left\{ \int \limits_{T_i \cap T_{i,j}^{\dualk{k}} } \omega_{p\ov{T_i}}\left(\xi\right) \psi_{q\ov{T_{i,j}^{\dualk{k}}}}\left(\xi\right)  d\xi \right\}_{p,q=1,..,(N+1)^d}\label{eq:L2proj}
\end{align}
where $\mathbf{M}$ 
 is called \emph{mass matrix} which is diagonal according to the chosen polynomial basis, it is positive definite and it gives the projection coefficients between the piecewise polynomials  of   $\p_N(T_i)$; $\mathbf{M}^{-1}$ 
  is its inverse; $\omega_{q}$ 
	is the p-th basis polynomial centered for $\p_N(T_i)$; $\psi_{p}$ 
is the p-th basis polynomial for $\p_N(T_{i,j}^{\dualk{k}})$; $p$ ranging from $p=1$, $2$, $\ldots$ $(N+1)^d$; $\mathbf{M}^{\sigma(i,j,k)}$ gives the projection coefficients between the piecewise polynomial basis of   $\p_N(T_i)$ and  $\p_N(T_{i,j}^{\dualk{k}})$. It is important to stress the fact that, with a little linear algebra, the aforementioned three-dimensional tensors can be written in a simpler and versatile composition of one-dimensional matrices. 


%
%
%
%

		\subsection{Numerical method}  
		\label{sec:AMRDGSI}
		Once the main and the staggered grids have been defined, the space of discrete solutions $\p_N$ is outlined, then the numerical method follows the procedure of \cite{FambriDumbser}. 
		
		A weak formulation of the governing equations (\ref{eq:NSmom}-\ref{eq:NSinc}) is considered in the form   
		\begin{align} 
		& \int \limits_{T_m^{\dualk{k}}} \psi_p \left( \frac{\partial v_k}{\partial t} + \nabla \cdot \mathbf{F}_{v_k} + \partial_{x_k} p \right) d\x = 0,  &(x_k,v_k) = (x,v_x),(y,v_y),(z,v_z),  \label{eq:DGNSmom}\\
& \int \limits_{T_i} \omega_p \, \left( \nabla \cdot \mathbf{v}\right) d\x = 0,  & T_m^{\dualk{k}}\in \Omega_h^{\dualk{k}} , \;\; T_i \in \Omega_h \label{eq:DGNSinc}
\end{align}
where it is important to highlight that the momentum equations are integrated along the \emph{dual} AMR meshes $\Omega_h^{\dualk{k}}$ according to the $k$-th space-direction, and the incompressibility condition along the main AMR mesh $\Omega_h$; $\psi_p \in \basis_{N}(T_{i,j}^{\dualk{k}})$, i.e. the $p$-th element of the basis for the \emph{dual} vector space $\p_N(T_{i,j}^{\dualk{k}})$; $ \omega_p \in \basis_{N}(T_i)$, i.e. the $p$-th element of the basis for the vector space $\p_N(T_i)$. In this formalism, it is clear that the chosen weak formulation of our staggered DG scheme consists 
in an $L_2$ projection of the incompressible Navier-Stokes equations onto the chosen vector space of staggered piecewise polynomials $\p_N$. 
 Then, according to \cite{DumbserCasulli2013,FambriDumbser}, the discrete solution of the staggered DG discretization is defined as follows:  
the velocity components $v_k$ are approximated in the respective vector spaces of the staggered dual meshes $\p_N(\Omega_h^{\dualk{k}})$, while the pressure is defined on the main grid, i.e. 
in $\p_N(\Omega_h)$. Therefore, the respective expansion over the polynomial basis reads 
\begin{align}
 &v_{h}^{\kapped}(\mathbf{x},t)\ov{T_m^{\dualk{k}}} = v_m^{\dualk{k}}(\mathbf{x},t) = \sum\limits_{l=1}^{(N+1)^d} \psi_{l}(\x)  \, \hat{v}_{l}(t) & &    \text{for}\; \mathbf{x}\in T_m^{\dualk{k}}, \;\;\; m=1,..., N^{\kapped}_{\text{faces}}\;\; \label{eq:vDGpoly} \\
& && \psi_k \in \basis_N(T_m^{\dualk{k}}), \;\; \hat{v}_l \in \real  \nonumber\\
& p_{{h}}(\mathbf{x},t)\ov{T_{i}} = p_{i}(\mathbf{x},t) = \sum\limits_{l=1}^{(N+1)^d} \omega_{l}(x)  \, \hat{p}_{l}(t) && \text{for}\; \mathbf{x}\in T_i, \;\;\; i=1,..., \Nel\;\; \label{eq:pDGpoly}\\
& && \omega_l \in \basis_N(T_i), \;\; \hat{p}_l \in \real  \nonumber
\end{align}
After substituting (\ref{eq:vDGpoly}) into the incompressibility condition (\ref{eq:DGNSinc}) one obtains 
\begin{equation}
  \int \limits_{\partial T_{i}} \omega_{l}(\x) \, \mathbf{v}_{h} \cdot \vec{n} \, dS 
	- \int \limits_{T_{i}} \nabla \omega_{l}(\x) \cdot \mathbf{v}_{h} \, d\x = \sum_k \left[\int \limits_{\Gamma_m^{\kapped} \in \partial T_{i}} \omega_{l}(\x) \, \mathbf{v}_{h} \cdot \vec{n} \, dS 
	- \int \limits_{T_m^{\dualk{k}} \cap T_{i}} \nabla \omega_{l}(\x) \cdot \mathbf{v}_{h} \, d\x \right] = 0. \label{eq:DGinc2}
\end{equation} 
Indeed, exactly because of the chosen $C$-staggering, $\mathbf{v}_{h} \cdot \vec{n}$ is well defined along $\partial T_i$, i.e.  all the velocity components are \emph{continuous} along the respective faces $\Gamma_{i,j}^{\kapped}$; moreover $\omega_k$ is continuous in $T_i$.
On the other hand, due to the staggering, the pressure $p_h$ is discontinuous within $T_{m}^{\dualk{k}}$. Then, after substituting the ansatz (\ref{eq:pDGpoly}) into the momentum equation (\ref{eq:DGNSmom}),  the pressure gradient terms can be written as
\begin{align}
\int_{T_m^{\dualk{k}}} \psi_l \partial_{x_k} p_h d\x \equiv \int_{T_m^{\dualk{k}}\cap T_{l(m)}} \psi_l \partial_{x_k} p_{l(m)} d\x + \int_{T_m^{\dualk{k}} \setminus T_{l(m)}} \psi_l \partial_{x_k} p_{r(m)} d\x + \int_{\Gamma_{m}^{\dualk{k}}} \psi_l \left( p_{r(m)} - p_{l(m)} \right) dS 
\end{align}
where $l(m)$ is the index of the left element with respect to $\Gamma_m^{\kapped}$. 

Notice that, whenever  $T_m^{\dualk{k}}$ is an $\US$ element, then $r(m)=l(m)\equiv i$ and then: (i) the surface integral of equation (\ref{eq:DGinc2}) vanishes because $\Gamma_m^{\kapped}=\varnothing$; (ii) the second volume integral vanishes because $(T_m^{\dualk{k}} \setminus T_{l(m)})\equiv \varnothing$; (iii) the jump term vanishes because $\Gamma_m^{\dualk{k}}\equiv \varnothing$\;
\footnote{In principle, for $\US$ elements, $\Gamma_m^{\dualk{k}}$ is something indefinite, and it can be set to be $\Gamma_m^{\dualk{k}}\equiv \varnothing$. This definition is a way of stating the pressure of being continuous inside $T_m^{\dualk{k}}$. Indeed, for $\US$ elements it holds $T_m^{\dualk{k}} \subset T_{r(m)\equiv l(m)}$.}.  
Otherwise, if $T_m^{\dualk{k}}$  is a \emph{standard} dual element, the second volume integral is evaluated exactly along the right element of $\Gamma_m^{\kapped}$, i.e.  $T_{r(m)} \equiv  (T_m^{*} \setminus T_{l(m)}) $
independently of the refinement level.

Then, after evaluating the definite integrals of system (\ref{eq:DGNSmom}-\ref{eq:DGNSinc}), a semi-implicit higher order DG finite element formulation of the incompressible Navier-Stokes equations on  staggered AMR grids reads\footnote{ Symbol $\cyclic{a}{b}{c}$ denotes the standard cyclic (or circular) permutations of $(a,b,c)$: i.e. $\left\{(a,b,c), (b,c,a), (c,a,b)\right\}$.}
%
\begin{align} 
  \Mass{x$_1$ x$_2$ x$_3$} \cdot \left( \dof{V}^{n+1}_{m}- \dof{Fv}^{n}_{m}  \right) + \frac{\Delta t}{\Delta x_1}\Mass{ x$_2$ x$_3$} \cdot\left( \tild{R}{x$_1\sigma(m,r(m))$} \cdot\dof{P}^{n+\theta}_{r(m)} - \tild{L}{x$_1\sigma(m,l(m))$}\cdot \dof{P}^{n+\theta}_{l(m)} \right)& = 0, 
 \label{eq:vSDG1}  \\
&   \forall T_m^{\dualk{x_1}} \in \Omega_h^{\dualk{x_1}},\ \forall (x_1,x_2,x_3) \in \cyclic{x_1}{x_2}{x_3}
	\nonumber \\
%
 \sum \limits_{ x_1=x,y,z} \ \sum \limits_{\Gamma_{i}^{\dualk{x_1}} \subset T_i}
\Delta x_2 \Delta x_3 \;
\Mass{x$_2$x$_3$} \left(\barr{R}{x$_1\sigma^*(i,r^*(i))$} \cdot\dof{V}^{n+1}_{r^*(i)} \!-\! \barr{L}{x$_1\sigma^*(i,l^*(i))$}\cdot \dof{V}^{n+1}_{l^*(i)}\right)
  & = 0. 
  \label{eq:pSDG1}  \\ 
	&
 \text{with}\;\; (x_1,x_2,x_3) \in \cyclic{x_1}{x_2}{x_3} 
 \nonumber
\end{align}
The momentum equation (\ref{eq:vSDG1}) has been normalized with respect to the volume of $T_m^{\dualk{k}}$ divided by the time step $\Delta t$ which is intended in the classical sense $\Delta t=t^{n+1}-t^n$; $\dof{V}$ is the vector of the degrees of freedom of the $x_1$-th velocity component; $l^*(i)$ and $l^*(i)$ are the integer indexes of the  right and left \emph{dual} elements $T_{m}^{\dualk{x_1}}$ and $T_{n}^{\dualk{x_1}}$ that shares the dual face $\Gamma_i^{\dualk{x_1}}$%
\footnote{In this notation the integrals of eq. (\ref{eq:DGinc2}) break into the summations 
\begin{align}
\int\limits_{\partial T_i} \;\; \longrightarrow \sum_{x_1} \;\; \int \limits_{\Gamma_m^{\ixed{x_1}} \in \partial T_{i}} \longrightarrow \sum_{x_1}^{(0)}\;\;\sum^{(1)} \limits_{\Gamma_i^{\dualk{x_1}} \subset  T_{i}} \;\; \sum_{n=l^*(i),r^*(i)}^{(2)} \;\;\int\limits_{ \Gamma_n^{\ixed{x_1}}}\;\; ,
\end{align}
where the first summation (1) is intended over all the dual faces $\Gamma_i^{\dualk{x_1}}$ within $T_{i}$; the second summation $(2)$ is intended over the right and left (or forward and backward) dual elements $T_{r^*(i)}^{\dualk{x_1}}$ and  $T_{l^*(i)}^{\dualk{x_1}}$ that shares the corresponding dual face $\Gamma_i^{\dualk{x_1}}$. 
}
  $\sigma(m,j)$ and $\sigma^*(j,m)$ labels the mutual position and refinement of face $\Gamma_m^{\ixed{x_1}}$   respect to element $T_j$, the dual face $\Gamma_i^{\dualk{x_1}}$ respect to the dual element $T_m^{\dualk{x_1}}$, respectively; the implicitness factor $\theta \in [0.5,1]$ has been introduced for the implicit time-discretization of the pressure gradients in the form $X^{\theta}=\theta X^{n+1} + (1-\theta)X^n$; note that if $\theta=1/2$ the well known second order Crank-Nicolson time discretization is obtained. $\dof{Fv}^{n}$ collects the advective-diffusive terms that can  be approximated by means of a fully explicit or a proper semi-implicit discretization within am operator-splitting scheme. Further details will appear in the next sections.
The following one-dimensional matrices have been introduced  
\begin{align} 
\tild{R}{$\sigma(m,r(m))$}  & = \frac{1}{ \min(|T_m^{\dualk{x_1}}|,|T_{r(m)}|)}\left\{  \int\limits_{\Gamma_{m}^{\dualk{x_1}}} \psi_{p\ov{T_{m}^{\dualk{x_1}}}}  \omega_{q\ov{T_{r(m)}}} dS  
+\int \limits_{T_{m}^{\dualk{x_1}} \cap T_{r(m)}}  \,\psi_{p\ov{T_{m}^{\dualk{x_1}}}}  \partial_{x_1} \omega_{q} {\ov{T_{r(m)}}} d\x \right\}_{p,q=1,..,(N+1)^d} \nonumber\\
& \equiv \left\{ R_{pq}^{\sigma(m,r(m))} \right\}_{p,q=1,..,(N+1)^d} \equiv \Mass{x$_2$x$_3$}\tild{R}{x$_1\sigma(m,r(m))$}\nonumber\\
\tild{L}{$\sigma(m,l(m))$}  & = \frac{1}{ \min(|T_m^{\dualk{x_1}}|,|T_{l(m)}|)} \left\{  \int\limits_{\Gamma_{m}^{\dualk{x_1}}} \psi_{p\ov{T_{m}^{\dualk{x_1}}}}  \omega_{q\ov{T_{l(m)}}} dS 
- \int \limits_{T_{m}^{\dualk{x_1}} \cap T_{l(m)}}  \,\psi_{p\ov{T_{m}^{\dualk{x_1}}}}  \partial_{x_1} \omega_{q} {\ov{T_{l(m)}}} d\x \right\}_{p,q=1,..,(N+1)^d} \nonumber\\
& \equiv \left\{ L_{pq}^{\sigma(m,l(m))} \right\}_{p,q=1,..,(N+1)^d} \equiv \Mass{x$_2$x$_3$}\tild{L}{x$_1\sigma(m,l(m))$}\nonumber
\\
\barr{R}{$\sigma^*(i,r^*(i))$}  &\equiv \left\{ \left[ L_{pq}^{\sigma(r^*(i),l(r^*(i)))}\right]^T \right\}_{p,q=1,..,(N+1)^d}  \equiv \left[\tild{L}{$\sigma(r^*(i),l(r^*(i)))$} \right]^T \equiv \Mass{x$_2$x$_3$}\left(\tild{L}{x$_1\sigma(r^*(i),l(r^*(i)))$}\right)^T\equiv \Mass{x$_2$x$_3$}\barr{R}{x$_1\sigma^*(i,r^*(i))$}\nonumber\\
\barr{L}{$\sigma^*(i,l^*(i))$}  &\equiv \left\{ \left[ R_{pq}^{\sigma(l^*(i),r(l^*(i)))}\right]^T \right\}_{p,q=1,..,(N+1)^d}  \equiv \left[\tild{L}{$\sigma(l^*(i),r(l^*(i)))$} \right]^T \equiv \Mass{x$_2$x$_3$}\left(\tild{L}{x$_1\sigma(l^*(i),r(l^*(i)))$}\right)^T\equiv \Mass{x$_2$x$_3$}\barr{R}{x$_1\sigma^*(i,r^*(i))$}\label{eq:matrices}
\end{align}
where $|\cdot|$ denotes the volume of the corresponding space-element. 
Because of the symmetry between the aforementioned matrices, that is verified after a little algebra (see \cite{FambriDumbser} for the purely Cartesian case), it becomes simpler to introduce the following notation
\begin{align}
& \mathcal{R}\equiv \tild{R}{}  \equiv \barr{L}{T} ,  & \mathcal{L} \equiv \tild{L}{} \equiv \barr{R}{T}.
\end{align}
It should be noticed that equations (\ref{eq:vSDG1}-\ref{eq:pSDG1}) constitute a particular coupled equation system with a typical \emph{saddle point structure} which is typical for
any discrete formulation of the incompressible Navier-Stokes equations. The total number of unknowns of the system is large, since we have $\Nel$ space-elements, 
each element with $(N+1)^d$ degrees of freedom. In principle, a direct solution can become cumbersome, since it involves four unknown quantities per degree of freedom: 
three velocity components and the scalar pressure. The complexity of the problem can be considerably reduced with a very simple manipulation. 

A well known numerical strategy, widely used for solving such a complex linear systems, is the application of the \emph{Schur complement} to the saddle point system (\ref{eq:vSDG1}-\ref{eq:pSDG1}). In the context of staggered grids, this procedure has been successfully adopted in the field of ocean modeling and free-surface dynamics \cite{CasulliWalters,Casulli2009,CasulliVOF}, physiological fluid flows in the arterial system  \cite{CasulliDumbserToro,Blood3D2014}, compressible fluids \cite{DumbserIbenIoriatti,DumbserCasulli2016} and a novel family of higher order DG methods \cite{DumbserCasulli2013,TavelliDumbser2014,TavelliDumbser2014b,TavelliDumbser2015,TavelliDumbser2016,FambriDumbser} which inspired this work. 
After multiplying equation (\ref{eq:DGNSmom}) by the inverse of the element mass matrix $\Mass{}$ and a direct substitution of equation (\ref{eq:DGNSmom}) into (\ref{eq:DGNSinc}), 
the following linear algebraic system for the pressure degrees of freedom as the only unknowns is obtained: 
\begin{align}
\mathbb{H}\cdot \mathpzc{P}^{n+\theta} \equiv \mathpzc{b}^n \label{eq:dPoisson}
\end{align}
in which operator $\mathbb{H}$ is the block coefficient matrix representing the discrete Laplace operator for the pressure Poisson equation in the space of solutions $\p_N$ within the chosen 
staggered-mesh framework; $\mathpzc{P}^{n+\theta}=\theta \mathpzc{P}^{n+1} + (1-\theta)\mathpzc{P}^{n}$ is the complete vector of the degrees of freedom for the pressure; $\mathpzc{b}^{n}$ collects all the known terms, i.e. the non-linear advection and diffusion, see Section \ref{sec:AD}. In particular,  $\mathpzc{P}^{n+\theta}$ multiplied by the $i$-th row of $\mathbb{H}$ reads 
\begin{align}
\mathbb{H}^i\cdot \mathpzc{P}^{n+\theta} \equiv &  \sum \limits_{\cyclic{x_1}{x_2}{x_3}} \ \sum \limits_{\Gamma_{i}^{\dualk{x_1}} \subset T_i}
\frac{ \Delta x_2 \Delta x_3\Mass{x$_2$x$_3$} \left(  \Hoperator{R}{x$_1$}  \cdot\dof{P}^{n+\theta}_{r(r^*(i))} +  \Hoperator{C}{x$_1$}  \cdot\dof{P}^{n+\theta}_{r(l^*(i))\equiv l(r^*(i))\equiv i} +  \Hoperator{L}{x$_1$}  \cdot\dof{P}^{n+\theta}_{l(l^*(i))}    \right)}{\Delta x_1}  . 
  \label{eq:psyst} 
\end{align} 
where 
\begin{align}
\Hoperator{R}{} &= \mathcal{L}^T\Mass{x$_1$}   \mathcal{R} &&=\barr{R}{x$_1\sigma^*(i,r^*(i))$} \Mass{x$_1$}   \tild{R}{x$_1\sigma(r^*(i),r(r^*(i)))$}& \equiv \Hoperator{L}{T} \label{eq:R} \\
\Hoperator{L}{} &=  \mathcal{R}^T\Mass{x$_1$}   \mathcal{L} &&=\barr{L}{x$_1\sigma^*(i,l^*(i))$} \Mass{x$_1$}   \tild{L}{x$_1\sigma(l^*(i),l(l^*(i)))$}& \equiv \Hoperator{R}{T} \label{eq:L}\\
\Hoperator{C}{} &=  \mathcal{L}^T\Mass{x$_1$}   \mathcal{L}  + \mathcal{R}^T\Mass{x$_1$}   \mathcal{R}  &&=\barr{R}{x$_1\sigma^*(i,r^*(i))$} \Mass{x$_1$}   \tild{L}{x$_1\sigma(r^*(i),l(r^*(i)))$} + \barr{L}{x$_1\sigma^*(i,l^*(i))$} \Mass{x$_1$} \tild{R}{x$_1\sigma(l^*(i),r(l^*(i)))$} & \equiv \Hoperator{C}{T} \label{eq:C}
\end{align}
It can be explicitly shown, see \cite{FambriDumbser}, that the Laplace operator can be decomposed into a more familiar composite matrix product, i.e.
\begin{align}
\mathbb{H} \equiv \mathcal{D} ^T \mathbb{M}^{-1}\mathcal{D},
\end{align}
where $\mathcal{D}$ is the weak gradient operator obtained after an $L_2$ projection of the gradient $\nabla$ in the vector space $\p_N(\Omega_h^*)$, i.e.
\begin{align}
\mathcal{D}_{x_1}^m \mathpzc{P} \equiv \frac{\SP{\psi}{\partial_{x_1} p_h}_{T_m^{\dualk{x_1}}}}{\SP{\psi}{\psi}_{T_m^{\dualk{x_1}}} }= \frac{\mathcal{R} \cdot \dof{P}_{r(m)} - \mathcal{L} \cdot \dof{P}_{l(m)}}{\Delta x_1}; \label{eq:D}
\end{align}

The real advantages of the chosen staggered framework arise in the properties of $\mathbb{H}$, which defined the discrete pressure Poisson system, which is at the same time the most complex linear 
system to be solved in our algorithm. 
Equations (\ref{eq:R}-\ref{eq:C}) together with the very simple algebraic analysis given by \cite{FambriDumbser} for the purely Cartesian case state that the discrete Laplace operator $\mathbb{H}$ 
satisfies the following peculiar properties: 

(i) $\mathbb{H}$ is \textbf{symmetric}, i.e.  $\mathbb{H}\equiv \mathbb{H}^T$; 

(ii) $\mathbb{H}$  is \textbf{always} \emph{at least}  \textbf{positive semi-definite}, independent of the chosen boundary conditions; in particular, $\mathbb{H}$ is \emph{positive definite} \emph{up to the solutions of}
\begin{align}
&\mathcal{D} \, \mathpzc{P}'  = 0  &\Longleftrightarrow &&\mathpzc{P}' \in \text{Ker}(\mathcal{D}),
\end{align}
where $\text{Ker}(\cdot)$ is the kernel set;

(iii) for \emph{uniform Cartesian grids}, $\mathbb{H}$ is \textbf{ block hepta-diagonal} for the three-dimensional case and  only \textbf{block penta-diagonal} for the two dimensional case.

(iv) whenever any \emph{pressure boundary condition} is imposed, then $\mathbb{H}$ is shown to be \textbf{strictly positive definite}.

\noindent It should be mentioned the resulting algebraic system (\ref{eq:dPoisson}) would occur also after choosing a different mesh staggering, but the corresponding non-zero blocks (\emph{stencil}) of the local equation (\ref{eq:psyst}) would be larger and the computational effort needed for solving the discrete Poisson equation (\ref{eq:dPoisson}) would presumably increase. When AMR is active, then  the mesh is non-conforming and each element can have more than $2d$ neighbors. As a consequence, the stencil size of (\ref{eq:psyst}) becomes slightly larger. In particular, whenever a neighbor element of $T_i$ is refined, the stencil of $\mathbb{H}^i$  adds $(\err^{d-1}-1)$ non-zero blocks, $\err$ being the refinement factor, $\err^{d-1}$ being the maximum number of nearest neighbors in a given oriented space-direction. Then, given a refinement factor $\err$, the largest stencil size for $\mathbb{H}$ will be $(2d\err^{d-1} +1)$, $2d+1$ being the stencil size for the uniform Cartesian case, $2d$ being the number of faces of a Cartesian element. Remember that the refinement levels of two neighbor elements are allowed to differ at most by one. 
 
It should be mentioned that in the uniform Cartesian case, a rigorous theoretical analysis of $\mathbb{H}$  for the design of specific preconditioners, using the theory of matrix-valued symbols and Generalized Locally Toeplitz (GLT) algebras (see \cite{serra1998,GSz,glt}) has been very recently provided with promising results in terms of numerical efficiency \cite{SIDG_analysis2016} and showing  beneficial properties of the respective condition number. 

Finally, once the equation system for the degrees of freedom of the pressure (\ref{eq:psyst}) has been solved by means of a classical matrix-free conjugate gradient method, then the 
velocity components can be updated directly via (\ref{eq:vSDG1}). In all numerical examples shown in this paper no preconditioner has been used. 

		\subsection{Advection and diffusion}   
		
		\label{sec:AD}

Any stable explicit DG scheme can be adopted for solving the advective-diffusive terms $\dof{Fv}$ in equation (\ref{eq:vSDG1}).
In order to simplify the computation, in this work we considered a DG formulation on the collocated grid, with $v_h \in \p_N(\Omega_h)$, having
\begin{align}
\dof{Fv}^{n}_{i} &= \dof{V}^n_{i} - \frac{\Delta t}{|T_i|} \left(\Mass{xyz}\right)^{-1} \cdot\left(\, \int \limits_{ \partial T_{i}} \boldsymbol{\omega} \mathbf{F}_v\cdot \vec{n} dS - \int \limits_{T_{i}} \nabla \boldsymbol{\omega} \cdot \mathbf{F}_v \, d\x \right).\label{eq:Fv}
\end{align} 
Because of the discontinuities of the piecewise polynomials along the element faces, a very simple and classical Rusanov flux (or local Lax-Friedrichs flux - LLF) \cite{Rusanov1961a,Toro09} has been used for evaluating the surface integral of (\ref{eq:Fv}) in the form
\begin{align}
\mathbf{F}_q\cdot\vec{n} &= \frac{1}{2}\left(\mathbf{F}^+_q + \mathbf{F}^-_q \right)\cdot\vec{n} - \frac{1}{2}s_q \left(q^+-q^-\right)\;\;\; \text{with}\;q=u,v,w.
\end{align}
in which the penalty term $s_q$ is the maximum value of the Jacobian of the  flux tensor $\mathbf{F}_q$
\begin{align}
s_{q} &= 2 \, \text{max}\left(|q^+|,|q^-|\right) + 2\nu\frac{2N+1}{\Delta x_q\sqrt{\left.{\pi}\middle/{2}\right. }}\;\;\;\;\;\; \text{with} \; (q,\Delta x_q) = (u,\Delta x), (v,\Delta y),(w,\Delta z). 
\end{align}
The chosen Rusanov flux has been modified in order to account for both hyperbolic and parabolic terms (see \cite{MunzDiffusionFlux,DumbserNSE,HidalgoDumbser}). 
In this context, it should be emphasized that although the advection-diffusion equation contains parabolic terms, an appropriate numerical flux can be defined by the solution of the 
corresponding generalized Riemann problem, see \cite{MunzDiffusionFlux} for details. 

Since equation (\ref{eq:Fv}) is an explicit DG scheme, then a rather severe $\CFL$ time restriction with $\CFL<1$ becomes necessary for ensuring stability: 
\begin{align}
\Delta t = \text{CFL}\left[\left(2N+1\right)\left(\frac{|u_\text{max}|}{\Delta x_{\text{min}}} + \frac{|v_\text{max}|}{\Delta y_{\text{min}}} + \frac{|w_\text{max}|}{\Delta z_{\text{min}}}\right) + \left(2N+1\right)^2\left(\frac{2\nu}{\Delta x^2_{\text{min}}} + \frac{2\nu}{\Delta y^2_{\text{min}}}+\frac{2\nu}{\Delta z^2_{\text{min}}}\right)\right]^{-1}. \label{eq:CFL}
\end{align}
Once $\dof{Fv}^{n}_{i}$ has been computed on the main grid $\Omega_h$, the numerical solution is projected back to the dual space of solutions $p_N(\Omega_h^{\dualk{k}})$ 
according to the $L_2$ projection operators (\ref{eq:pistar}-\ref{eq:pi}), with the projection matrices (\ref{eq:L2proj}).

Notice that, the parabolic nature of the viscous terms introduces a quadratic dependence of the minimum step-size on $h/(2N+1)$, which can become extremely severe when 
more than one refinement level is considered for a viscous fluid flow at low Reynolds number. In such circumstances, the stability condition (\ref{eq:CFL}) can become too 
restrictive. Then, following \cite{FambriDumbser}, an implicit discretization of the diffusive fluxes is taken for the advection diffusion system  
\begin{align}
\frac{\partial \mathbf{v}}{\partial t} + \nabla \cdot \mathbf{F_c} - \nabla \cdot \boldsymbol{\sigma} &= 0.\label{eq:AD1}\\
 \boldsymbol{\sigma} & = - \nu \nabla \mathbf{v}.\label{eq:AD2}
\end{align}
In particular, the velocity components and the stress tensor components are the unknowns of the system. The velocity components are discretized within the piecewise polynomials of maximum degree $N$ along the main grid $\Omega_h$, the component of the stress tensor along the \textit{dual grid} $\Omega_h^{\dualk{k}}$. Then the final algorithm is obtained following exactly the same strategy adopted for the incompressible Navier-Stokes equations: 
\begin{enumerate}
	\item the governing equations (\ref{eq:AD1}) and (\ref{eq:AD2}) are projected along the chosen space of solutions, i.e. $\p_N(\Omega_h)$ and $\p_N(\Omega_h^{\dualk{k}})$ respectively; 
\item the unknown variables are assumed to live in the corresponding space of solutions, i.e. $v^{\kapped}_h \in \p_N(\Omega_h)$ and $F_h^{\kapped} \in \p_N(\Omega_h^{\dualk{k}})$, and then substituted into the discrete DG 
  equations;
\item the Schur complement of the resulting algebraic coupled system is solved, after substituting the equation for the tensor components $F_h^{\kapped}$ into the momentum equations; then the decoupled system for the only velocity components $v^{\kapped}_h$ is obtained, in a very familiar structure, i.e.
\begin{align} 
&\left(\mathbb{M} + \nu  \mathbb{H}\right) \cdot \mathpzc{V}_{(k)}^{n+\halb} = \mathpzc{b_2}, & k=x,y,z
, \label{eq:implicitviscosity} 
\end{align}
where $\mathbb{M}^i\equiv \Mass{}$ is the mass matrix, $\mathbb{H}$ is exactly the afore-defined discrete Laplace operator, $\mathpzc{V}_{(k)}^{n+\halb}$ is the vector of the complete set of degrees of freedom for the $k$-th velocity component; $\mathpzc{b_2}$ collects all the known terms, i.e. the explicit advective terms of equation (\ref{eq:Fv}) multiplied by a mass matrix; $n+\halb$ labels a fictitious time step $t^{n+\halb}$, labeling the intermediate stage within the global algorithm; 
\item system (\ref{eq:implicitviscosity}) is very efficiently solved by means of a classical conjugate gradient method. Throughout this paper, we do not employ any preconditioner;
\item the velocity components are projected back to the dual space $\p_N(\Omega_h^{\dualk{k}})$ into $\dof{Fv}^{n}_{m}$ of equation (\ref{eq:vSDG1}). 
\end{enumerate}

Notice that, although equation system (\ref{eq:implicitviscosity}) must be solved for each velocity component, the coefficient matrix $\left(\mathbb{M} + \nu  \mathbb{H}\right)$  has even better 
properties than $\mathbb{H}$: (i) it is \textbf{symmetric}; (ii) it is \textbf{strictly positive definite}; (iii) for the \emph{purely Cartesian case}, $\left(\mathbb{M} + \nu  \mathbb{H}\right)$ is only \textbf{ block hepta-diagonal} for the three-dimensional case and only  \textbf{block penta-diagonal} for the two dimensional case.
(iv) the viscosity acts as a perturbation coefficient of the mass matrix, which is \textbf{purely diagonal} (not only block diagonal); then the system is surely much better behaved compared to the 
discrete pressure system. 
\noindent
After the aforementioned operator-splitting procedure for the semi-implicit time discretization, since the parabolic terms are treated fully implicitly, the  final $\CFL$ time-restriction relaxes to 
\begin{align}
\Delta t = \frac{\text{CFL}}{\left(2N+1\right)}\left(\frac{|u_\text{max}|}{\Delta x_{\text{min}}} + \frac{|v_\text{max}|}{\Delta y_{\text{min}}} + \frac{|w_\text{max}|}{\Delta z_{\text{min}}}\right)^{-1}.\label{eq:CFL1}
\end{align}

\subsection{Final algorithm}
		\label{sec:FA}
Finally, even in the staggered AMR framework, the resulting numerical scheme can be summarized, similar to \cite{FambriDumbser}, as follows: 
\begin{enumerate}
\item the velocity is projected to the main grid $\Omega_h$ with $\pi$ defined in equation (\ref{eq:pi}); advection terms are discretized explicitly while the viscous terms are treated implicitly for each velocity component within $\p_N(\Omega_h)$
\begin{align} 
\left(\mathbb{M} + \nu  \mathbb{H}\right) \cdot \mathpzc{V}^{n+\frac{1}{2}}_{\#} =  \mathbb{M}  \cdot\mathpzc{Fv}^n_{\#},   \label{eq:LOdiffusionV}    
\end{align}
symbol '$\#$' labels the degrees of freedom respect to the main space of solutions $\p_N(\Omega_h)$;
\item the velocity is projected back to the dual space $\p_N(\Omega_h^{\dualk{k}})$  according to the projection operator $\pi^{\dualk{k}}$ defined in equation (\ref{eq:pistar}) and the system for the discrete pressure Poisson equation is solved 
\begin{gather}
\mathbb{H}\cdot \mathpzc{P}^{n+\theta}_{\#}   = \basis^{n+\frac{1}{2}}_{\#}, \label{eq:LOpressure}\\
\text{with} \;\;\; \basis^{n+\frac{1}{2}}_{\#}\equiv\mathbb{M}^{yz} \left(\mathcal{D}^{x}\right)^T  \cdot \mathpzc{U}^{n+\frac{1}{2}}_{*}    +
\mathbb{M}^{zx} \left(\mathcal{D}^{y}\right)^T \cdot \mathpzc{V}^{n+\frac{1}{2}}_{*}    + \mathbb{M}^{xy} \left(\mathcal{D}^{z}\right)^T \cdot \mathpzc{W}^{n+\frac{1}{2}}_{*},  \nonumber
\end{gather}
symbol '$*$' denotes the degrees of freedom with respect to the dual space of solutions $\p_N(\Omega_h^{\dualk{k}}))$;
\item the velocity components are directly updated according to
\begin{align} 
\mathpzc{V}^{n+1}_{*} &= \mathpzc{V}^{n+\frac{1}{2}}_{*} - \frac{\Delta t}{\Delta x_1}  \left(\mathbb{M}^{-1}  \mathcal{D}\right)^y \cdot\mathpzc{P}^{n+\theta}_{{\#}},   \label{eq:LOvelocityV}  
\end{align}
\end{enumerate}


\section{Numerical tests}
\label{sec:tests}
In order to verify the accuracy and robustness of the presented method together with the AMR framework, a series of non-trivial numerical tests is chosen. The  
test problems have been selected according to the following criteria: i) an analytical, numerical or experimental reference solution exists;  ii) the test problem involves all terms of the incompressible 
Navier-Stokes equations, i.e. the nonlinear convective terms as well as the viscous and pressure forces are equally important; iii) physical instabilities and energy dissipation from the 
largest to the smallest spatial scales are present (see the three-dimensional Taylor-Green vortex problem); iv) pressure, velocity and wall boundary conditions are considered; 
v) the AMR framework should provide a remarkable benefit in terms of resolution without introducing spurious mesh-effects; vi) whenever possible, a convergence table should be provided according 
to a smooth \emph{analytical} reference solution.

		\subsection{2D Taylor-Green vortex problem} 
\begin{figure} 
\centering 
			\includegraphics[width=0.6\textwidth]{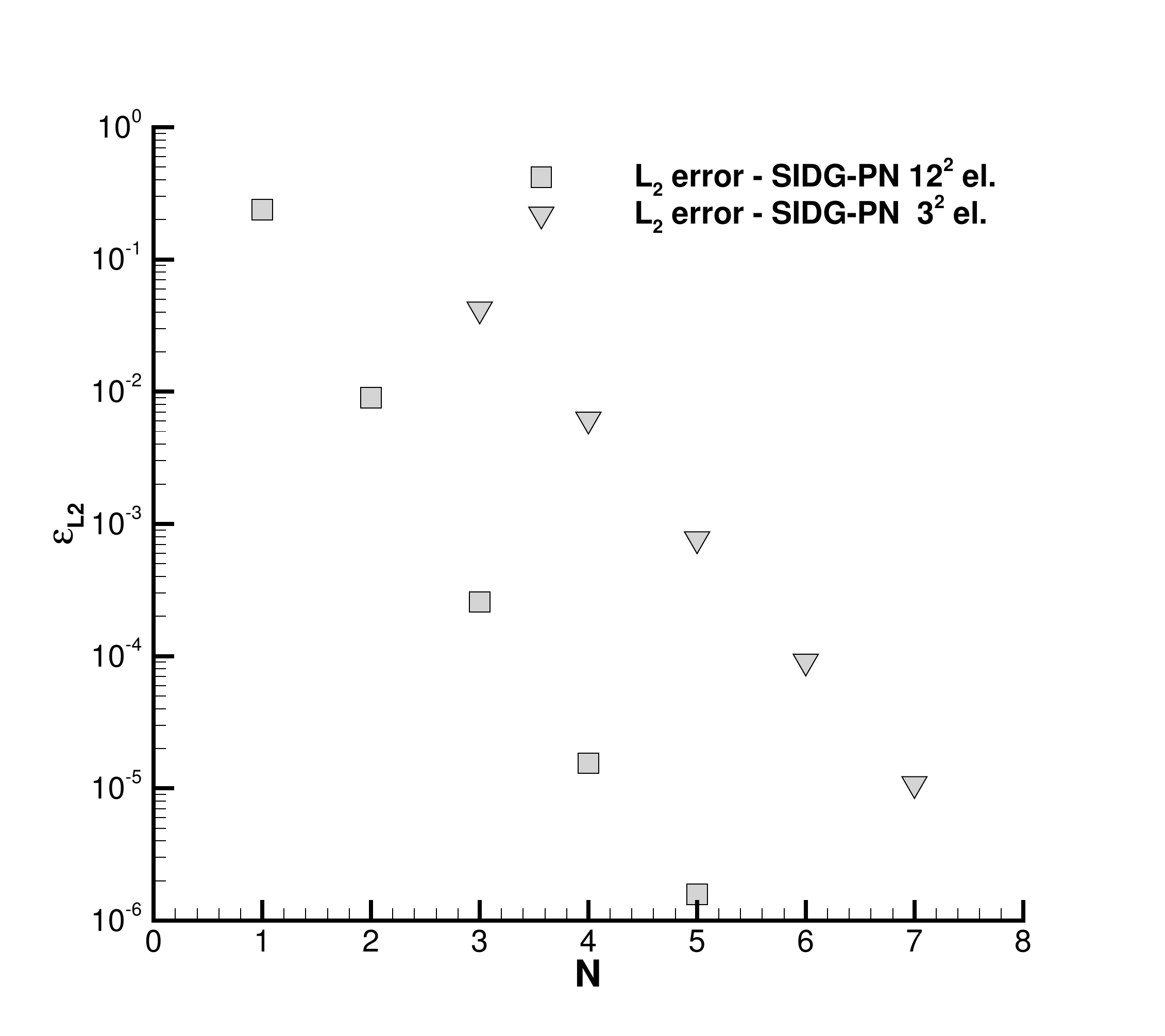} 
			\caption{Numerical $L_2$ error $\epsilon_{L_2}$ of the $u$ velocity component for the two dimensional Taylor-Green vortex problem computed with our $\SIDG$-$\p_{N}$ schemes as a function of the polynomial degree $N$ on a fixed grid of $12^2$ elements (squares) and $3^2$ elements (triangles).}\label{fig:spectral}
\end{figure}

A numerical convergence study should be carried out based on a smooth analytical reference solution of the governing PDE system, but there are only very few known 
analytical solutions for the incompressible Navier-Stokes equations available. Here, the two dimensional Taylor-Green vortex problem has been selected for testing the order of accuracy 
of our schemes. An exact solution of the two-dimensional incompressible Navier-Stokes equations is given by 
\begin{align}
&u(x,y,t) = \sin(x) \cos(y) e^{-2\nu t}, \qquad  v(x,y,t) = - \cos(x) \sin(y) e^{-2\nu t}, \nonumber \\
&p(x,y,t) = \frac{1}{4} \left( \cos(2x) + \cos (2y)\right) e^{-4\nu t}, \nonumber 
\end{align}
within the spatial domain $\Omega = [0,2\pi]^2$ and periodic boundary conditions everywhere. For the present test we use $\nu = 10^{-1}$. The flow dynamics corresponds to an exponential 
decay in time due to dissipative processes.
Without dissipation, the fluid flow is an unstable stationary equilibrium of inertial and pressure forces.
In order to compute properly a convergence table within the AMR framework, the ratio between finer and coarser elements should remain approximately constant in time and also for 
increasing mesh resolution. 
In this particular test the $L_1$, $L_2$ and $L_{\infty}$ errors have been evaluated at a final time $t_{\text{end}}=0.1$.
The resulting convergence study is summarized in Table \ref{tab:TGV2D} for polynomial degrees $N=1,\ldots,6$, fixing a refinement factor $\err=3$ and using up to $\ell_{\text{max}}=1$ refinement levels. 
Since the presented method is only second-order accurate in time, very small time-steps have been used. 
Furthermore, Figure \ref{fig:spectral} shows the dependence of the $L_2$ error on the polynomial degree $N$ for a given mesh. One can notice 
that an \emph{exponential decay} of the numerical error, called \emph{spectral decay}, with the polynomial degree is verified. 
These results indicate that for smooth problems the present algorithm works properly even when using very high order approximation polynomials and very coarse meshes. 
It should be emphasized that with very few higher order elements, e.g. $N=6$, the resulting numerical error is much smaller than the error obtained with an extremely 
refined grid at lower polynomial degrees (see Table \ref{tab:TGV2D} and Figure \ref{fig:spectral}). 

		\begin{table} 
\caption{Numerical convergence table computed for the two dimensional Taylor-Green vortex problem using staggered $\SIDG$-$\p_N$ schemes with AMR for polynomial degrees $N=1,\ldots,6$;  $\ell_{\text{max}}=1$; $\err=3$.} 
 \centering
\scriptsize
 \begin{tabular}{c r ccc ccc c c    r ccc ccc c   c }
   \hline
	&$N_{\text{el}}$ &  $\epsilon_{L^1}$&  $\epsilon_{L^2}$ & $\epsilon_{L^\infty}$ &  $\mathcal{O}_{L^1}$   & $\mathcal{O}_{L^2}$   & $\mathcal{O}_{L^\infty}$  &      &$N_{\text{el}}$ &  $\epsilon_{L^2}$&  $\epsilon_{L^1}$ & $\epsilon_{L^\infty}$ &  $\mathcal{O}_{L^1}$ &  $\mathcal{O}_{L^2}$   & $\mathcal{O}_{L^\infty}$  &      &  \\
		\hline
		&\multicolumn{7}{l}{$N=1$}&   & \multicolumn{8}{l}{$N=2$} & \\
   \cline{2-8}    \cline{10-17} \\
	& 12$^2$	&  1.133E+00 &	2.374E-01 &	8.664E-02&	---	 &  ---	 & ---	 &  &   6$^2$	&  3.938E-01	& 8.213E-02	 & 3.511E-02 &  ---  & ---	  & ---	 &	 	\\
	& 24$^2$	&  3.785E-01 &	7.630E-02 &	3.250E-02&	1.58 &	1.64 &	1.41 &  &   12$^2$	&  4.015E-02	& 8.997E-03	 & 4.312E-03 &  3.29 &	3.19  &	3.03 &   	\\  
	& 36$^2$	&  2.935E-01 &  6.074E-02 &	2.489E-02&	0.63 &	0.56 &	0.66 &  &   18$^2$	&  1.744E-02	& 4.012E-03	 & 2.420E-03 &  2.06 &	1.99  &	1.42 & 		\\
	& 48$^2$	&  2.085E-01 &   4.29E-02 &	1.74E-02 &	1.19 &	1.21 &	1.25 &  &   24$^2$	&  5.898E-03	& 1.495E-03	 & 1.057E-03 &  3.77 &	3.43  &	2.88 &  	\\ 
	\hline
	&\multicolumn{7}{l}{$N=3$}&   & \multicolumn{8}{l}{$N=4$} & \\
   \cline{2-8}    \cline{10-17} \\
	& 3$^2$	&  1.91E-01		&  4.19E-02		&  1.72E-02&	---	 &  ---	 & ---	 &  &   3$^2$	&  4.33E-02	&  	1.04E-02	&  5.84E-03 &  ---  & ---	  & ---	 &	 	\\
	& 6$^2$	&  1.27E-02		&  2.81E-03		&  1.48E-03&	3.92 &  3.90 & 3.54 &  &   6$^2$	&  1.88E-03	&  	4.20E-04	&  3.10E-04 & 4.53  & 4.63  & 4.23 &   	\\  
	& 9$^2$	&  3.05E-03		&  6.32E-04		&  5.10E-04&	3.51 &  3.68 & 2.63 &  &   9$^2$	&  3.08E-04	&  	7.05E-05	&  5.52E-05 & 4.46  & 4.40  & 4.25 & 		\\
	& 15$^2$&  5.18E-04	    &  	1.13E-04	&  	8.34E-05&	3.47 &  3.37 & 3.55 &  &   12$^2$	&  9.18E-05	&  	2.04E-05	&  1.68E-05 & 4.20  & 4.32  & 4.14  & \\  
	\hline
	&\multicolumn{7}{l}{$N=5$}&   & \multicolumn{8}{l}{$N=6$} & \\
   \cline{2-8}    \cline{10-17} \\
	& 3$^2$	&  6.40E-03	 & 1.45E-03	 & 5.97E-04&	---	 &  ---	 & ---	 &  &   3$^2$	&  5.37E-04	& 	1.27E-04	& 6.91E-05 &  ---  & ---  & ---	 &	 	\\
	& 6$^2$	&  1.01E-04	 & 2.24E-05	 & 1.05E-05&	5.99 & 6.02	 & 5.82 &  &   4$^2$	&  1.61E-04	& 	3.47E-05	& 1.34E-05 &  4.20 & 4.50	& 5.71  &   	\\  
	& 9$^2$	&  9.98E-06	 & 2.42E-06	 & 1.68E-06&	5.71 & 5.48	 & 4.52 &  &   6$^2$	&  7.96E-06	& 	1.81E-06	& 9.86E-07 &  7.41 & 7.28 & 6.43  & 		\\
	& 12$^2$&  2.99E-06	& 7.14E-07	 & 4.04E-07&	4.19 & 4.24	 & 4.96 &  &   8$^2$	&  1.59E-06	&   3.23E-07	& 1.73E-07 &  5.59 &	5.99 &	6.05 &  	\\

   \hline
 \end{tabular}
\label{tab:TGV2D}
 \end{table}

		\subsection{2D laminar Blasius boundary layer}
\begin{figure} 
\centering 
\begin{tabular}{cc} 
			\includegraphics[height=0.25\textwidth]{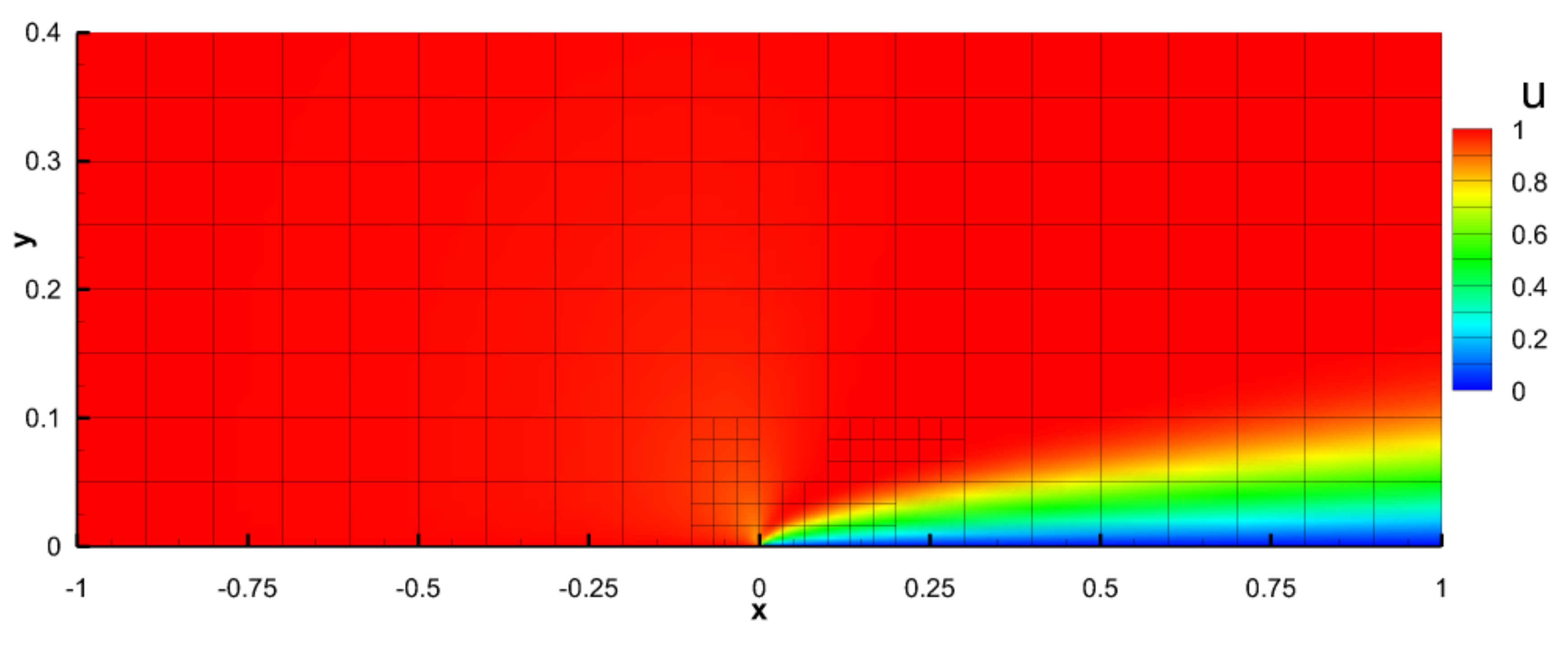} & 
			\includegraphics[height=0.35\textwidth]{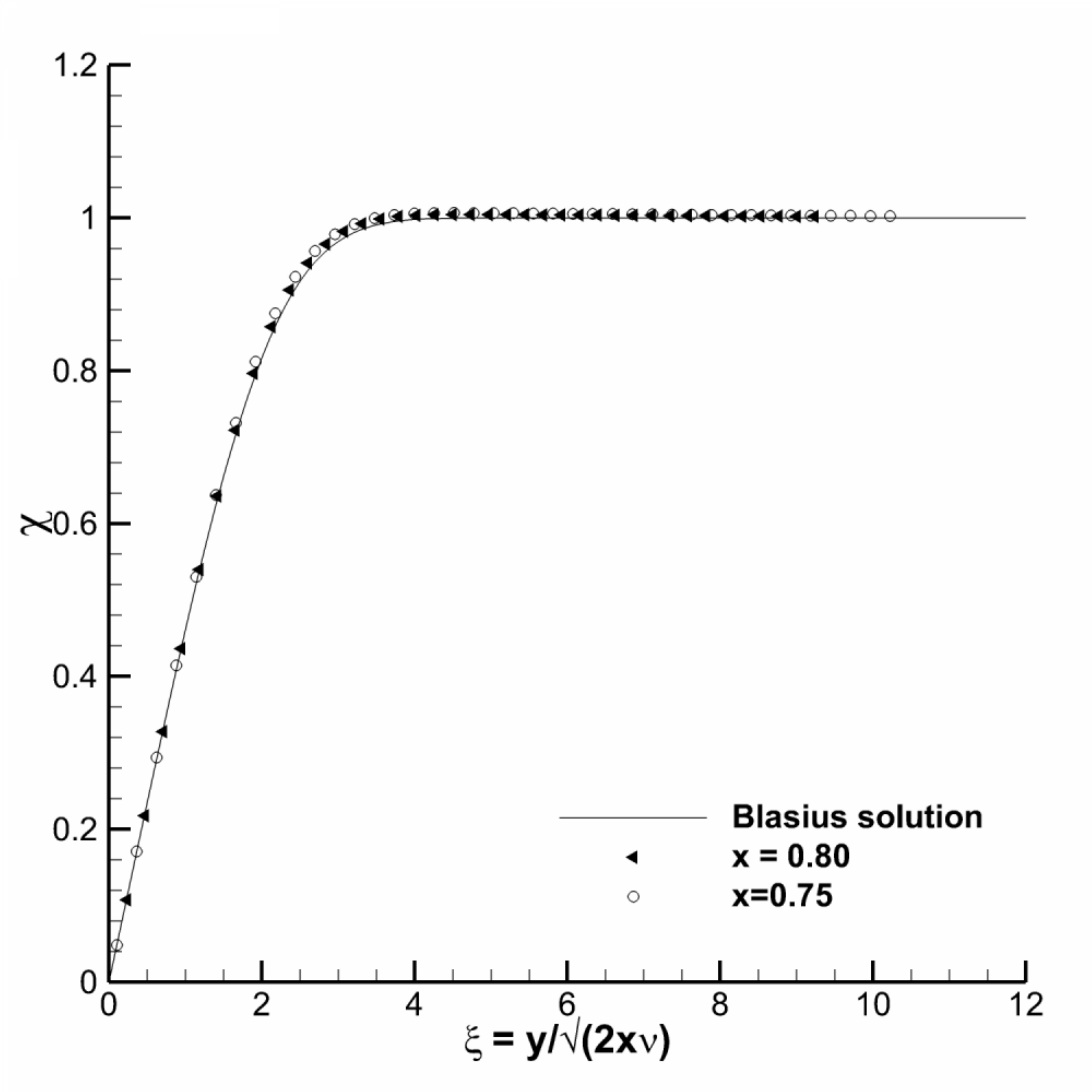} 
\end{tabular} 
\caption{The numerical solution obtained for the two dimensional laminar boundary layer test. The two dimensional view of the horizontal velocity field is shown on the left. The velocity profile interpolated along the vertical direction at two different axial position, i.e. $x=0.75$, $0.80$, is plotted next to the reference Blasius solution on the right. The results have been obtained 
with our staggered semi-implicit AMR-DG-$\p_{5}$ together with a refinement factor of $\mathfrak{r}=3$ and  $\ell_\text{max}=1$ levels. }\label{fig:BLayer}
\end{figure}

In this section a laminar two dimensional boundary layer over a flat plate is simulated. A rectangular domain $\Omega = [-1,1]\times[0.0,0.5]$ is initialized with a $20\times10$ grid at the 
coarsest AMR level $\ell=0$. The flow is driven by a constant velocity profile $\mathbf{v}= ( 1, 0)$ at the left entrance and a constant pressure condition ($p=0$) at the right and at the top outflow. 
At the bottom, i.e. at $y=0$, slip wall boundary conditions are applied for $x<0$ and the usual no-slip wall-boundary condition for $x \geq 0$. The resulting fluid flow is the classical 
laminar flow over a flat plate, for which a reference solution is given by the well-known Blasius solution \cite{Blasius1908}. 
For this test the kinematic viscosity is set to $\nu=10^{-3}$. 
Figure \ref{fig:BLayer} shows the numerical results obtained with the $\p_4$ version of our $\SIDG$ method using an adaptive mesh with $\ell_{\text{max}}=1$ and a refinement factor of 
$\err=3$. The grid is refined according to the gradient of the velocity magnitude $|\mathbf{v}|$. 

A very good agreement between the reference solution and the numerical solution obtained with the semi-implicit staggered DG scheme is observed. As it is shown in Figure \ref{fig:BLayer}, 
the refinement only takes place in the region of large velocity gradients, in particular at the leading edge of the boundary layer, as expected.


\subsection{Lid-driven cavity problem: 2D and 3D}
The lid-driven cavity problem is a well-known standard test problem for checking the correct implementation of boundary conditions, together with the accurate discretization of all 
the terms in the governing equations, i.e. the non-linear convection, the viscous forces, the pressure gradient and the incompressibility constraint. In this test the 
computational domain is the unit square $\Omega=[-0.5,0.5]^d$, representing a square cavity filled with an incompressible viscous fluid. A the upper boundary (the lid), a constant 
velocity vector $\mathbf{v}=(1,0,0)$ is imposed which drives the flow; no-slip wall boundary conditions are applied to the remaining boundaries. This is a non-trivial test mainly 
because of three facts: all the physical forces (pressure gradient, convection and diffusion) are of the same order of magnitude; there are two singularities in the velocity field at the 
upper corners, where the horizontal component of the velocity is a double valued function, i.e. $u=1$ at the top and $u=0$ at the lateral walls; the pressure field is only determined 
up to a constant, because boundary conditions are only imposed for the velocity but not for the pressure. At the initial time, the pressure is set to be $p=0$ everywhere. 
		
\begin{figure} 
\centering 
			\includegraphics[width=0.45\textwidth]{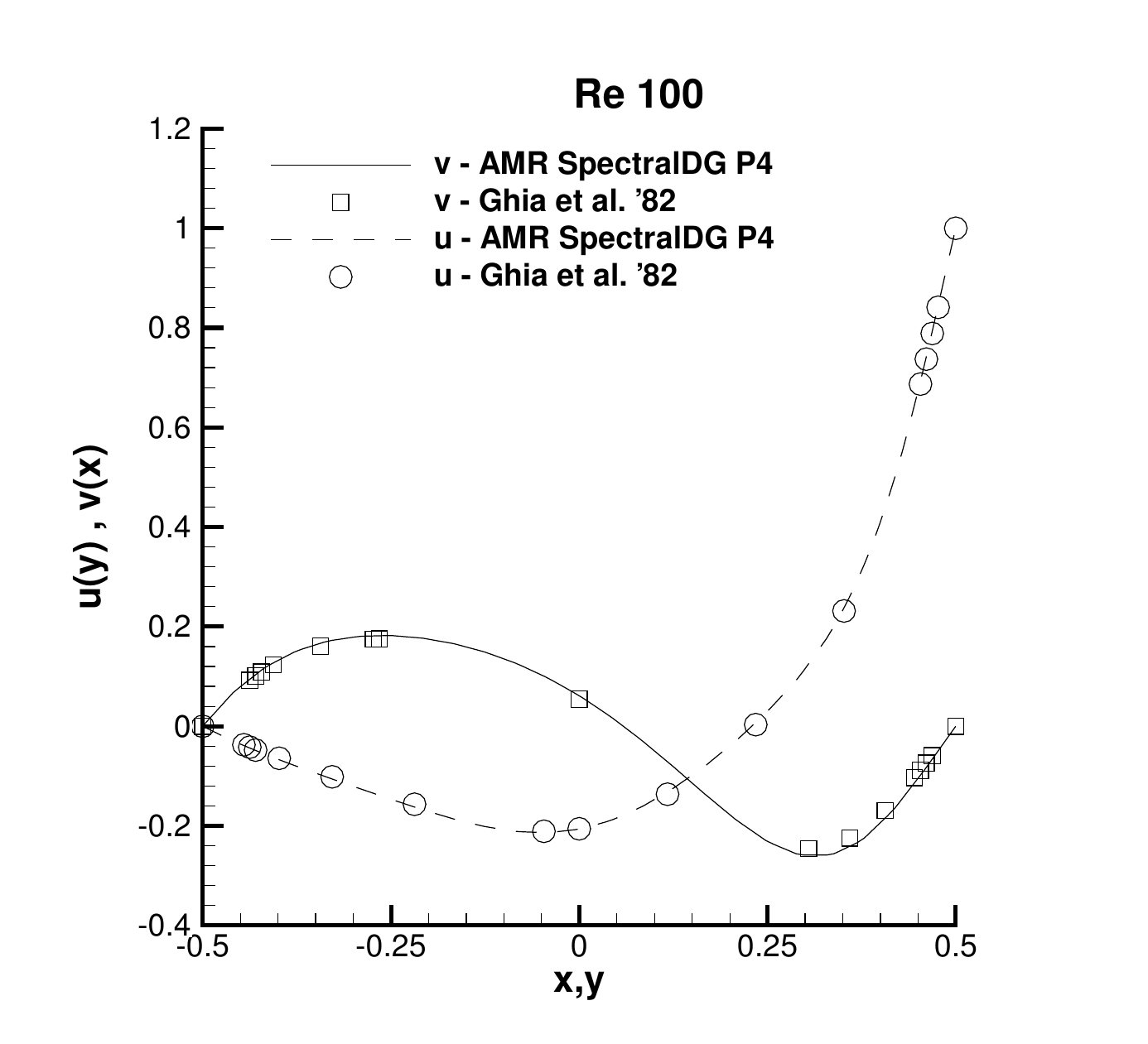}\;\includegraphics[width=0.45\textwidth]{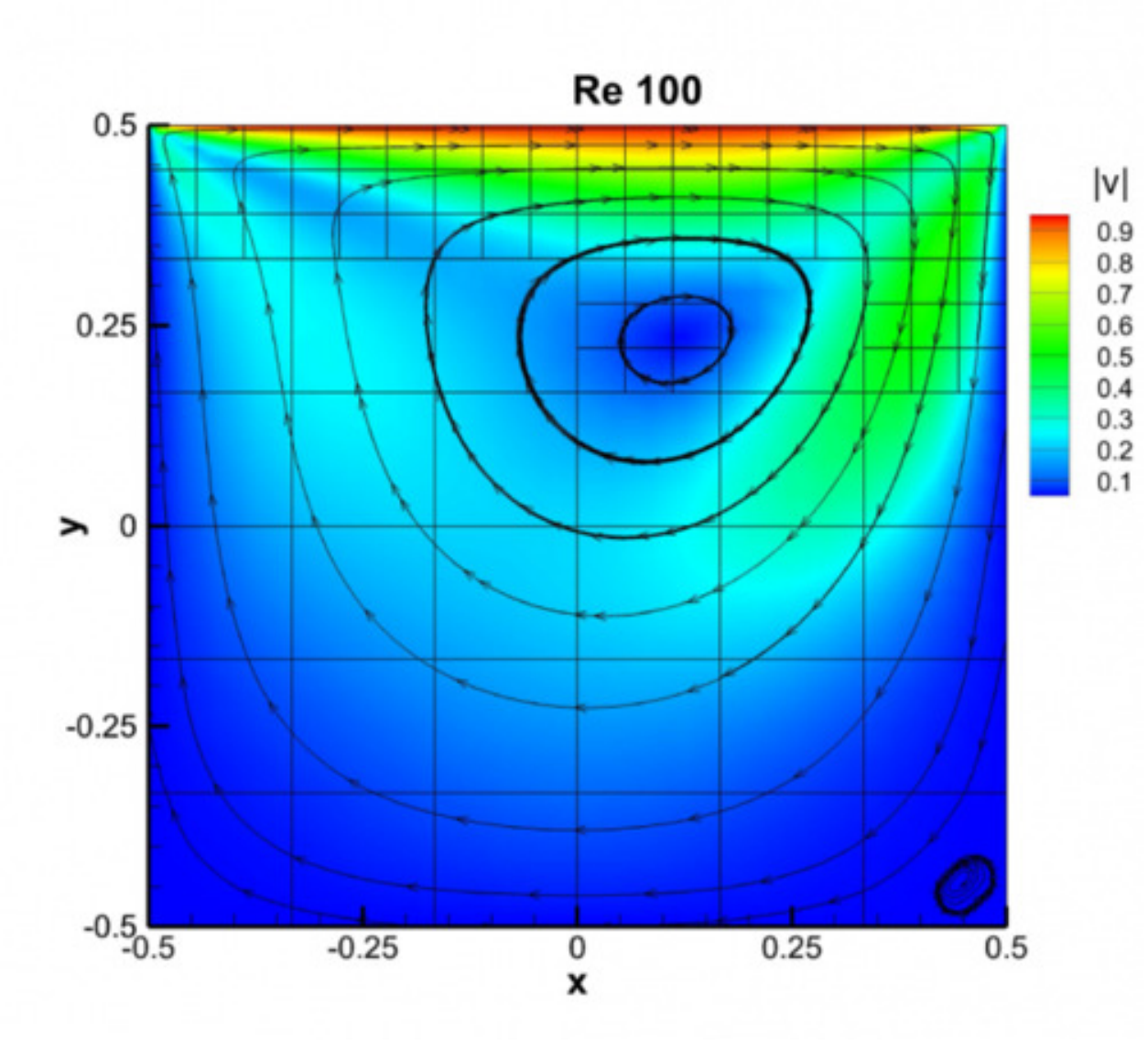}\\
			\includegraphics[width=0.45\textwidth]{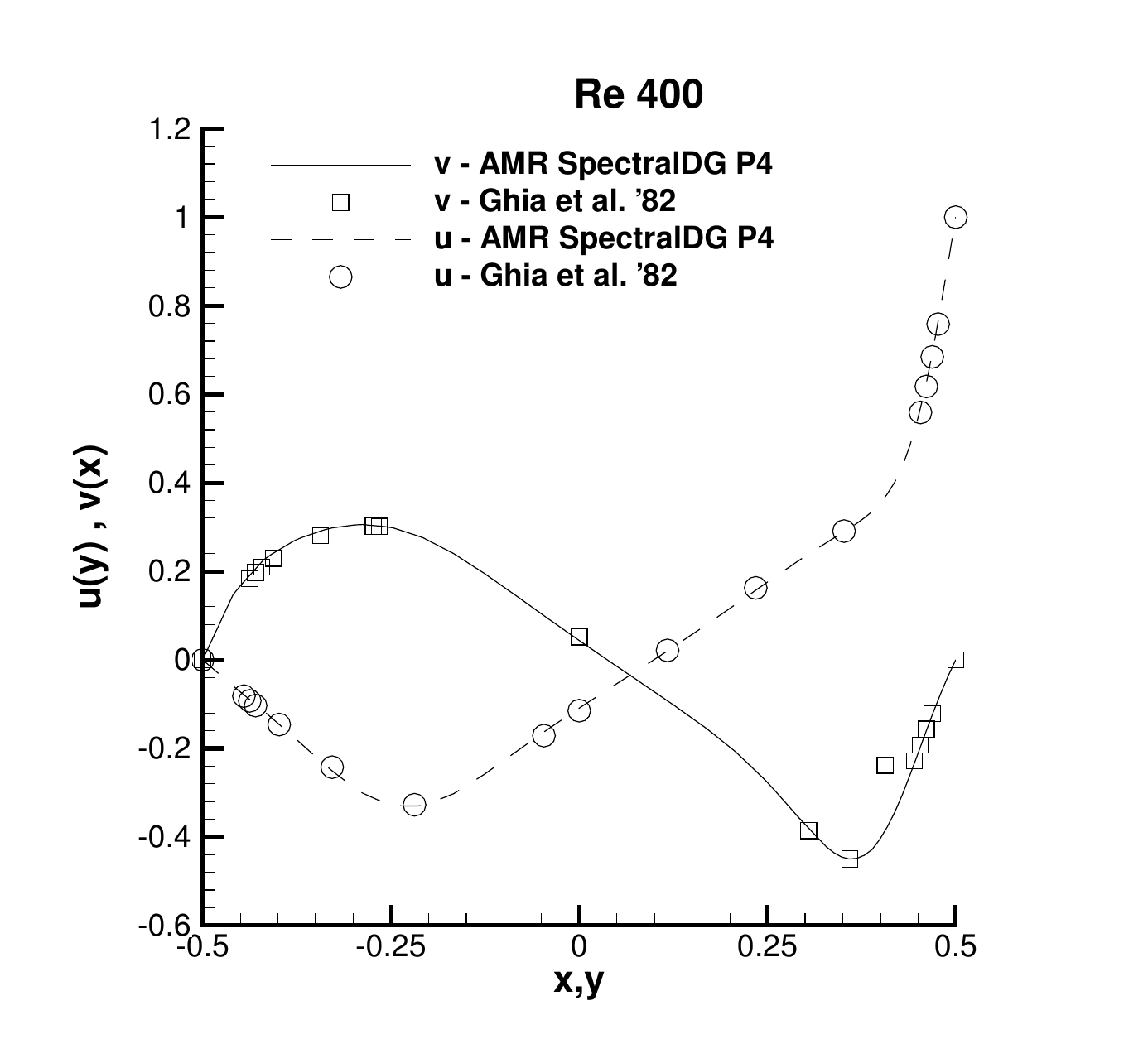}\;\includegraphics[width=0.45\textwidth]{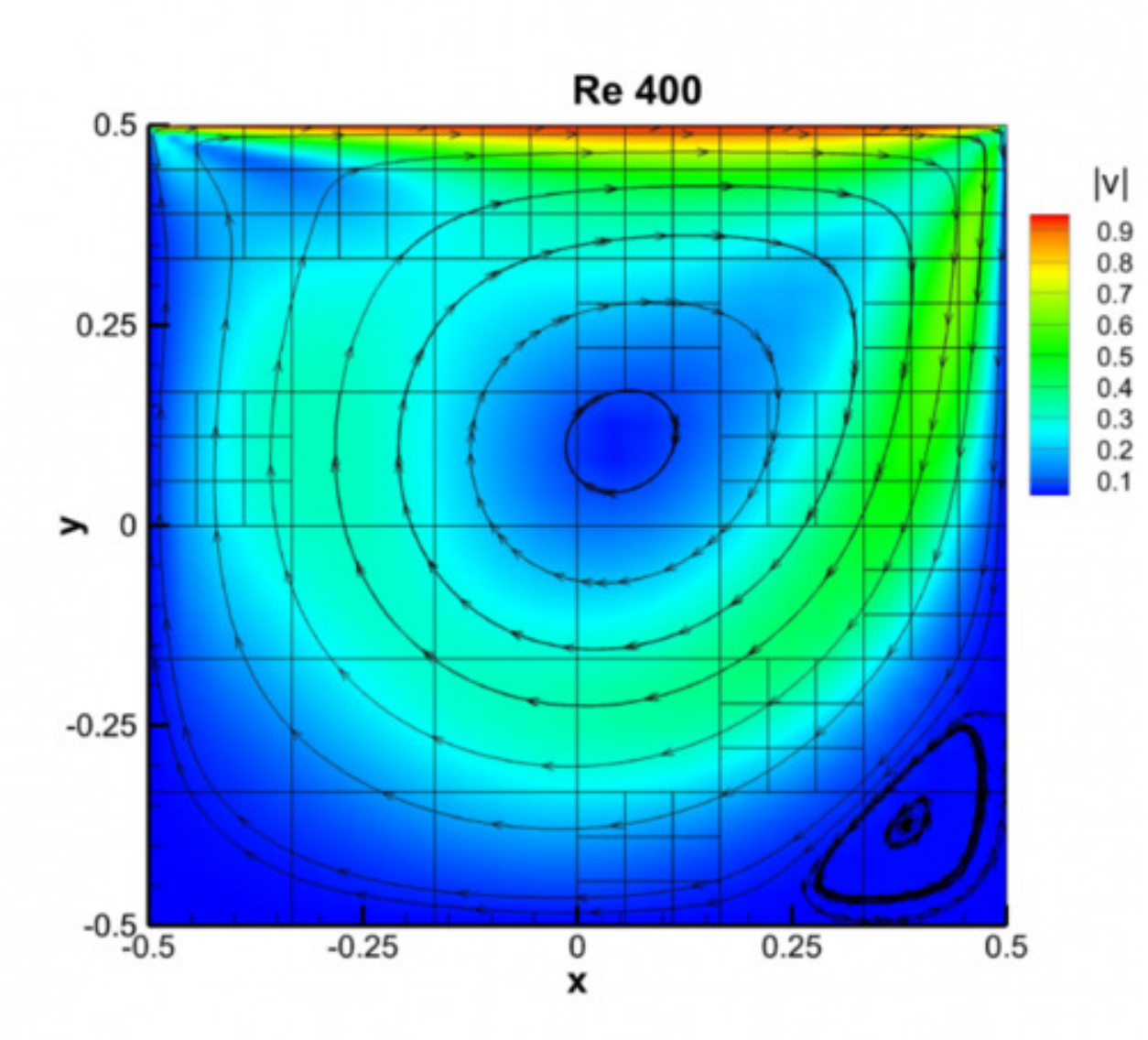}\\
			\includegraphics[width=0.45\textwidth]{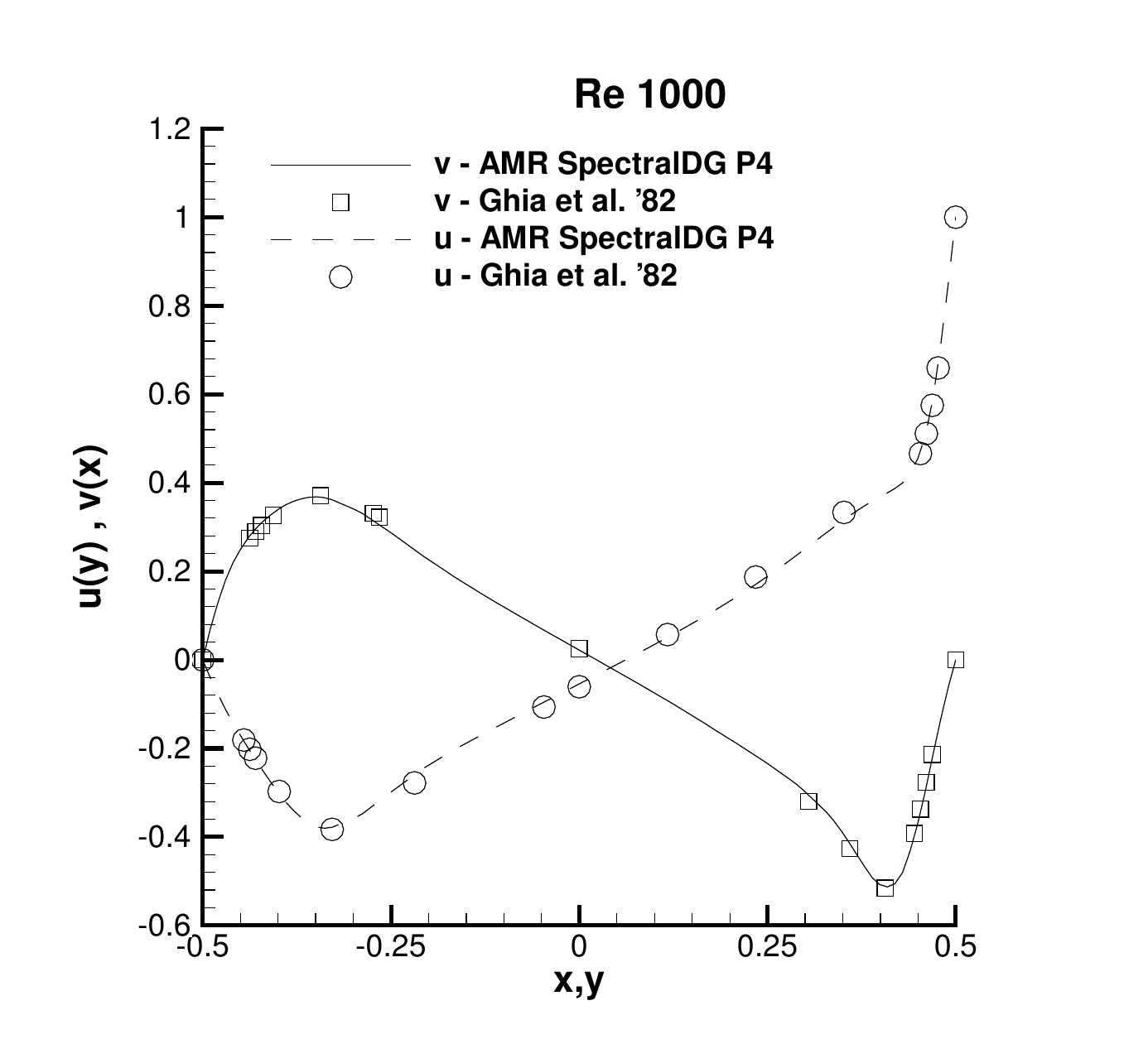}\;\includegraphics[width=0.45\textwidth]{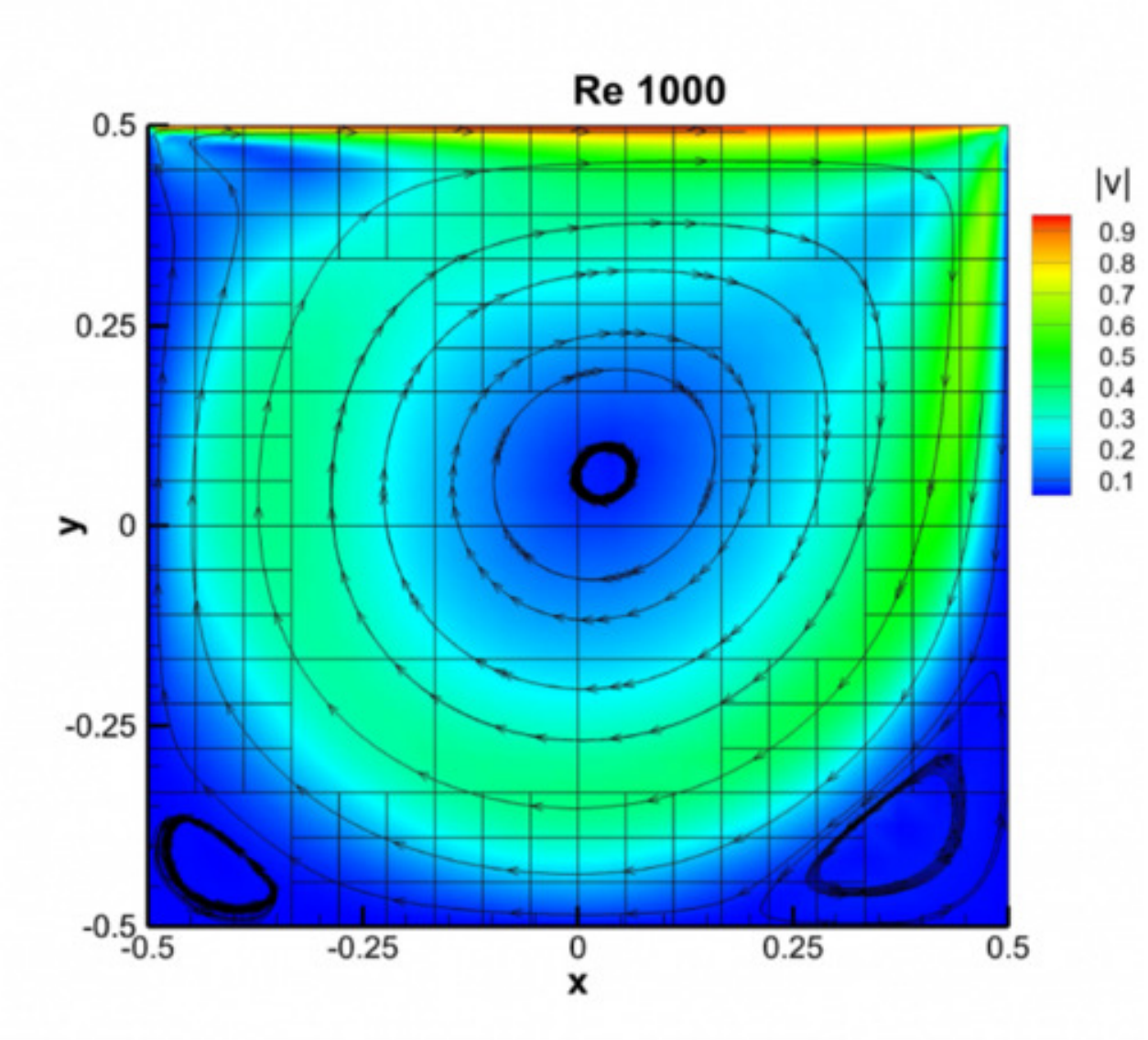}\\
\caption{The numerical solution obtained for the two dimensional lid-driven cavity problem compared with the numerical results of \cite{Ghia1982} at different Reynolds numbers, respectively, from the top to the bottom:  Re=$100$, Re=$400$ and Re=$1000$ using $6x6$ elements on the coarsest grid level. These results have been  obtained with the $\p_{4}$-version of our staggered semi-implicit spectral DG method.}\label{fig:LDCavity2Da}
\end{figure}

\begin{figure} 
\centering 
			\includegraphics[width=0.45\textwidth]{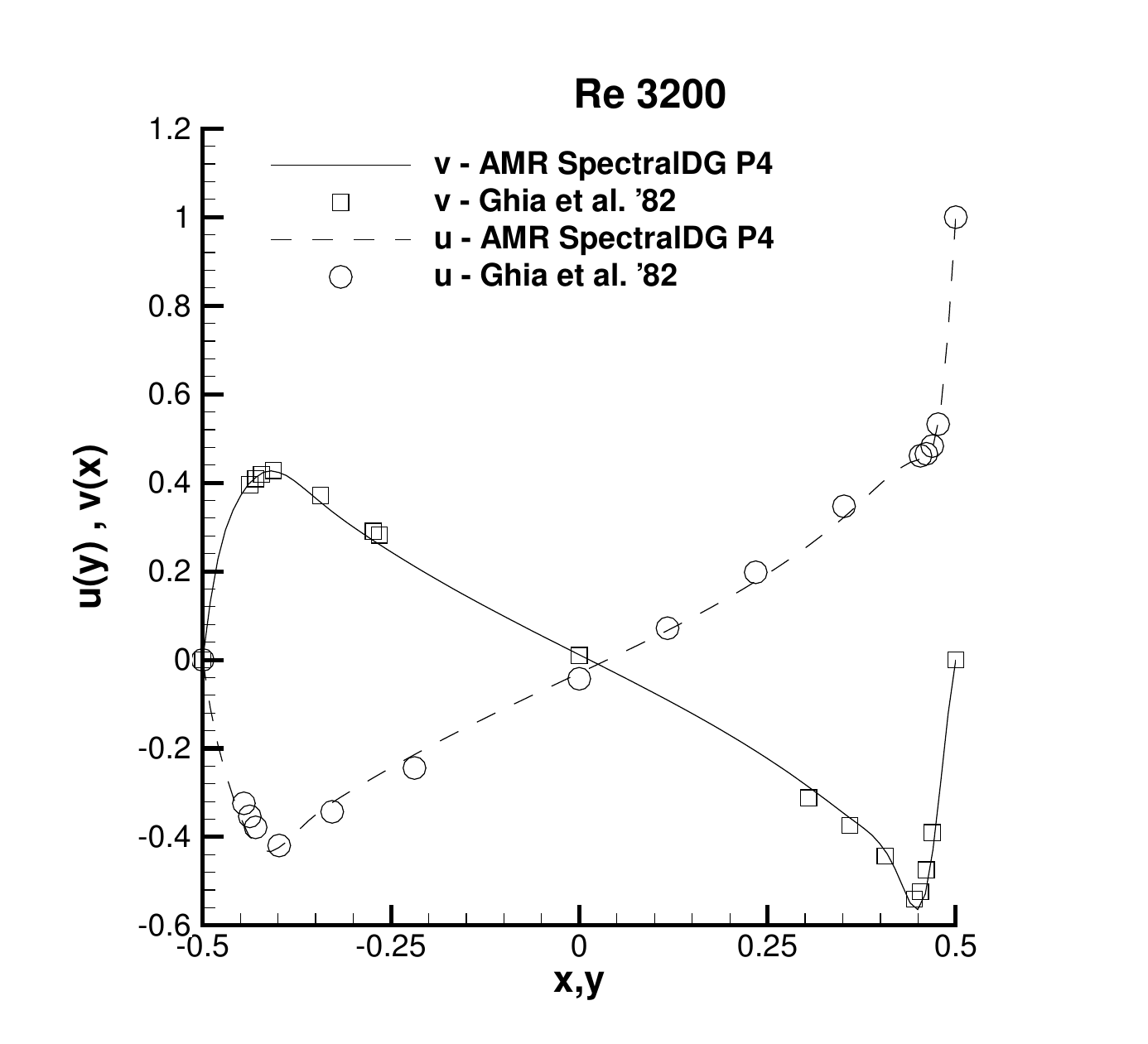}\;\includegraphics[width=0.45\textwidth]{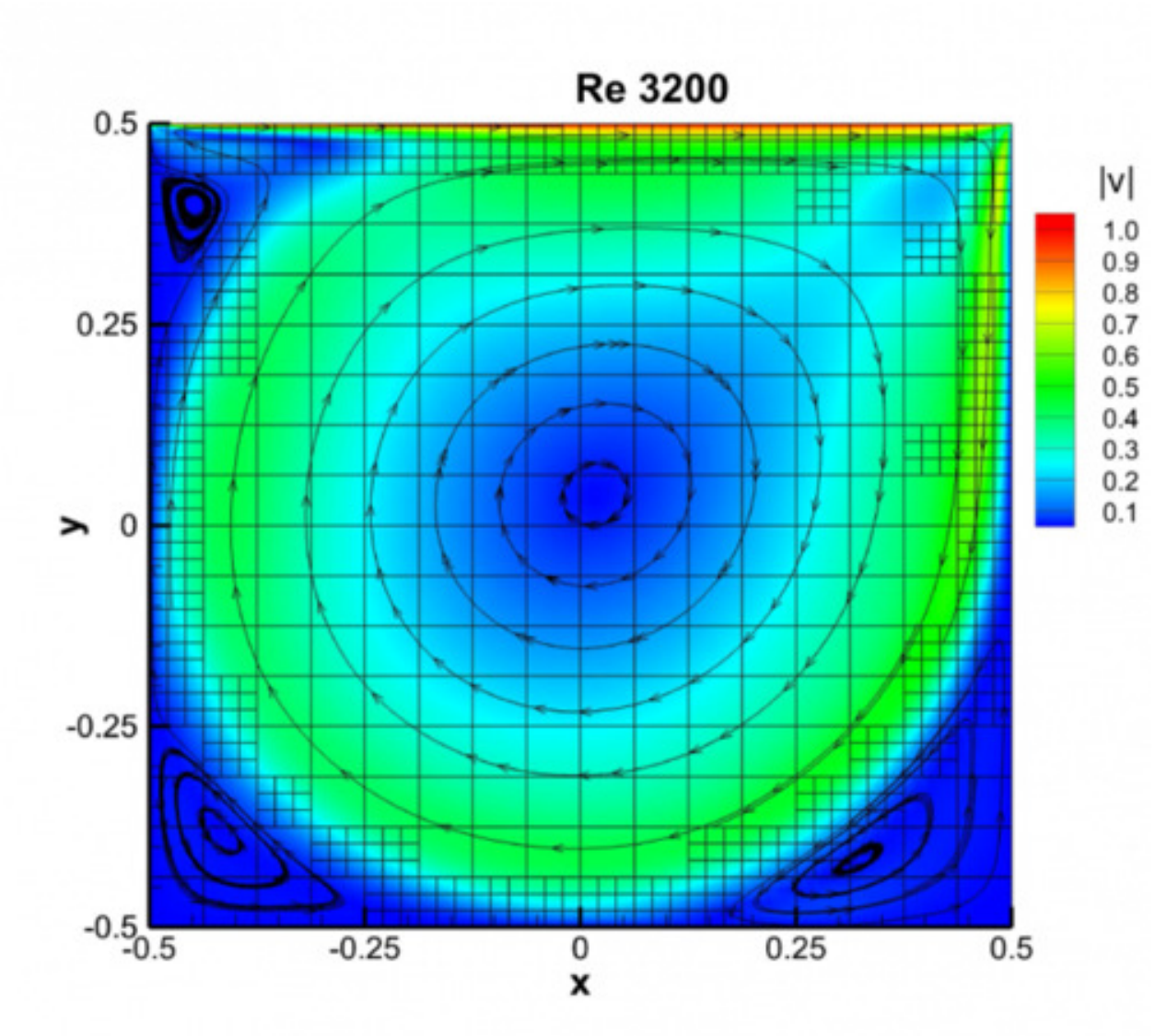}
\caption{The numerical solution obtained for the two dimensional lid-driven cavity problem compared with the numerical results of \cite{Ghia1982} at Reynolds number Re=$3200$ 
using $16x16$ elements on the coarsest grid level. These results have been obtained with the $\p_{4}$-version of our staggered semi-implicit spectral DG method.}\label{fig:LDCavity2Db}
\end{figure}

		Figures \ref{fig:LDCavity2Da} and \ref{fig:LDCavity2Db} show the numerical results in two space dimensions obtained at different Reynolds numbers in the range $Re\in[100,3200]$. The horizontal and vertical velocity profiles, interpolated along the vertical and the horizontal axes respectively, are plotted and compared with the reference solution given by Ghia et al. \cite{Ghia1982}. The AMR grid is adapted according to a refinement factor $\err=3$ and up to one single maximum refinement level ($\ell_{\max}=1$). As shown in Figures \ref{fig:LDCavity2Da} and \ref{fig:LDCavity2Db}, the computed results match the reference solution very well and the automatic adaptation is well-driven along the higher velocity gradients close to the walls.  In this case the \SIDG-$\p_{4}$ scheme has been used, corresponding to a total number of degrees of 
 freedom of $\ndof=5^2=25$ per space-element. On the level zero grid ($\ell_0$), the two-dimensional domain has been discretized within $6^3=36$ elements for Reynolds numbers 
$100$, $400$ and $1000$, while $16^2=256$ elements have been used for Reynolds number $3200$. 
		
\begin{figure} 
\centering 
			\includegraphics[width=0.48\textwidth]{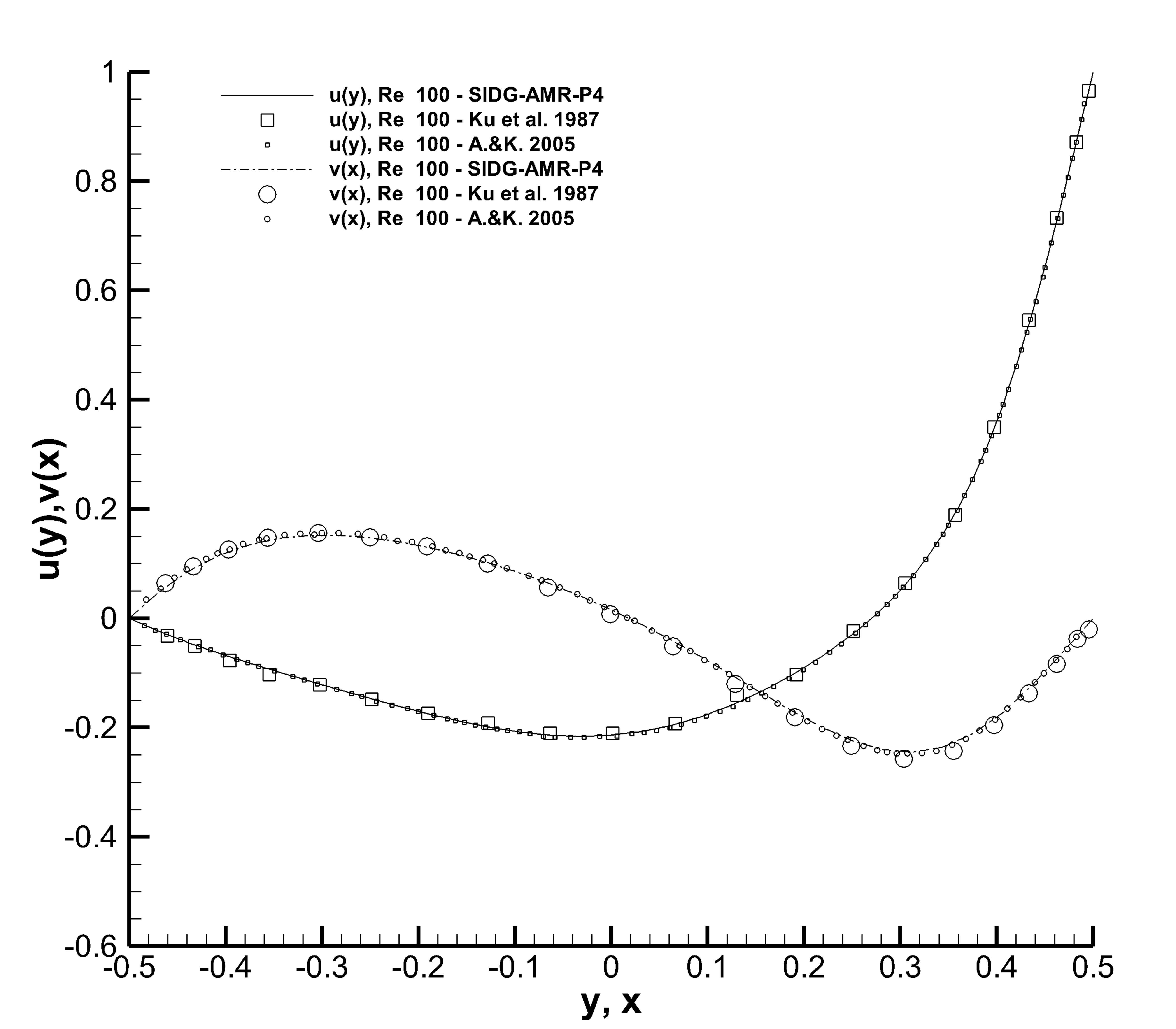}\;
			\includegraphics[width=0.48\textwidth]{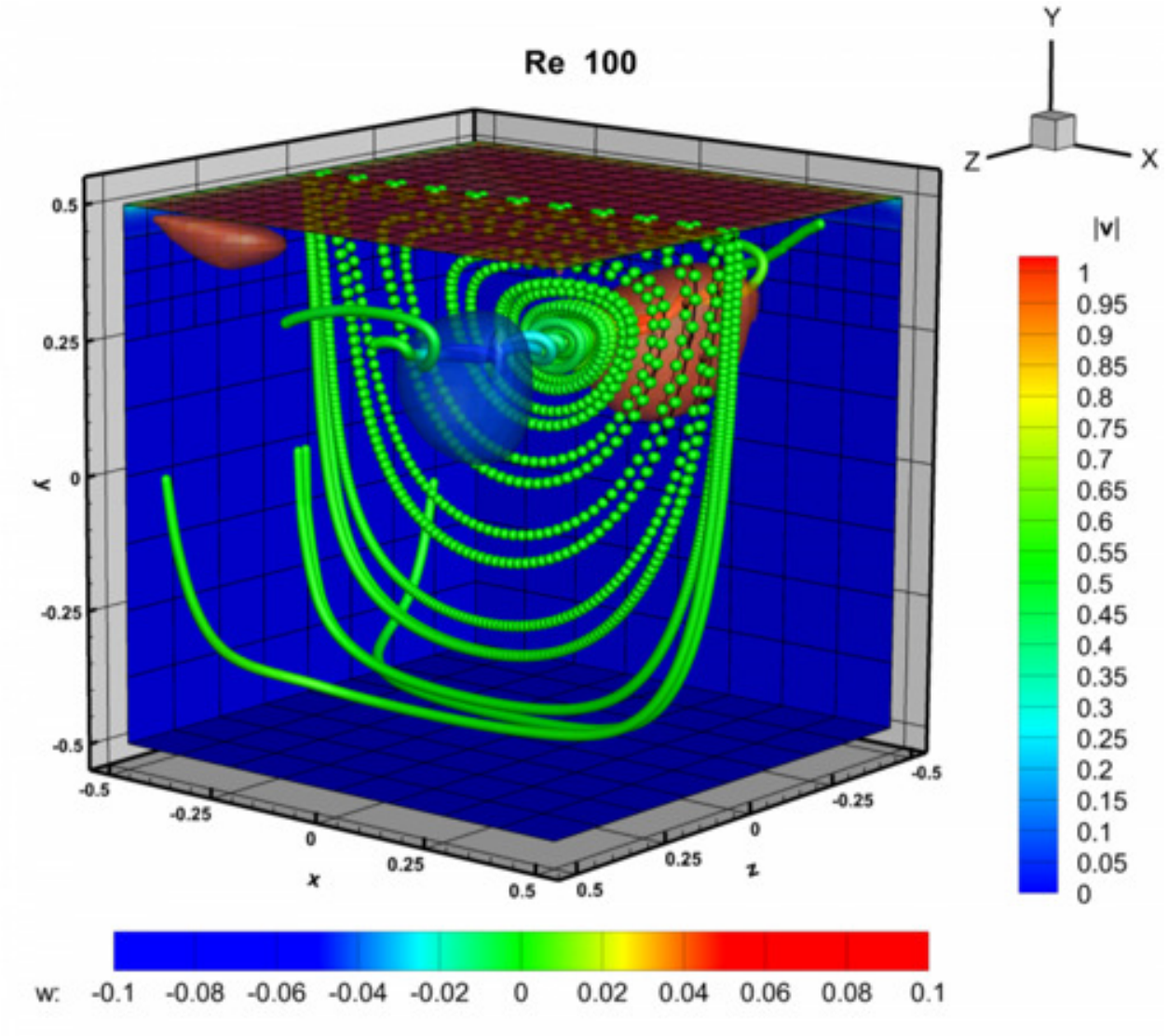}\\
			\includegraphics[width=0.48\textwidth]{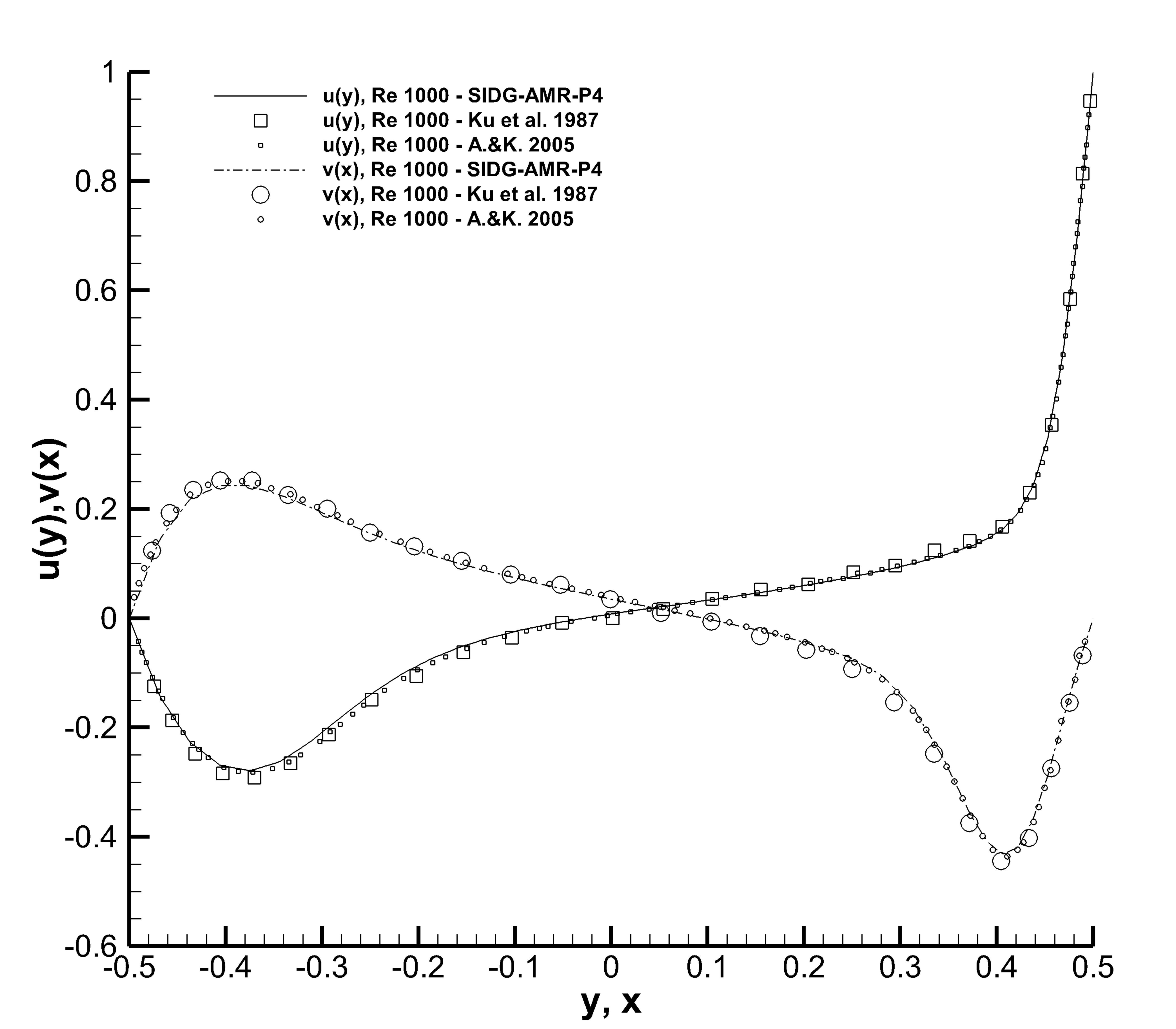}\;
			\includegraphics[width=0.48\textwidth]{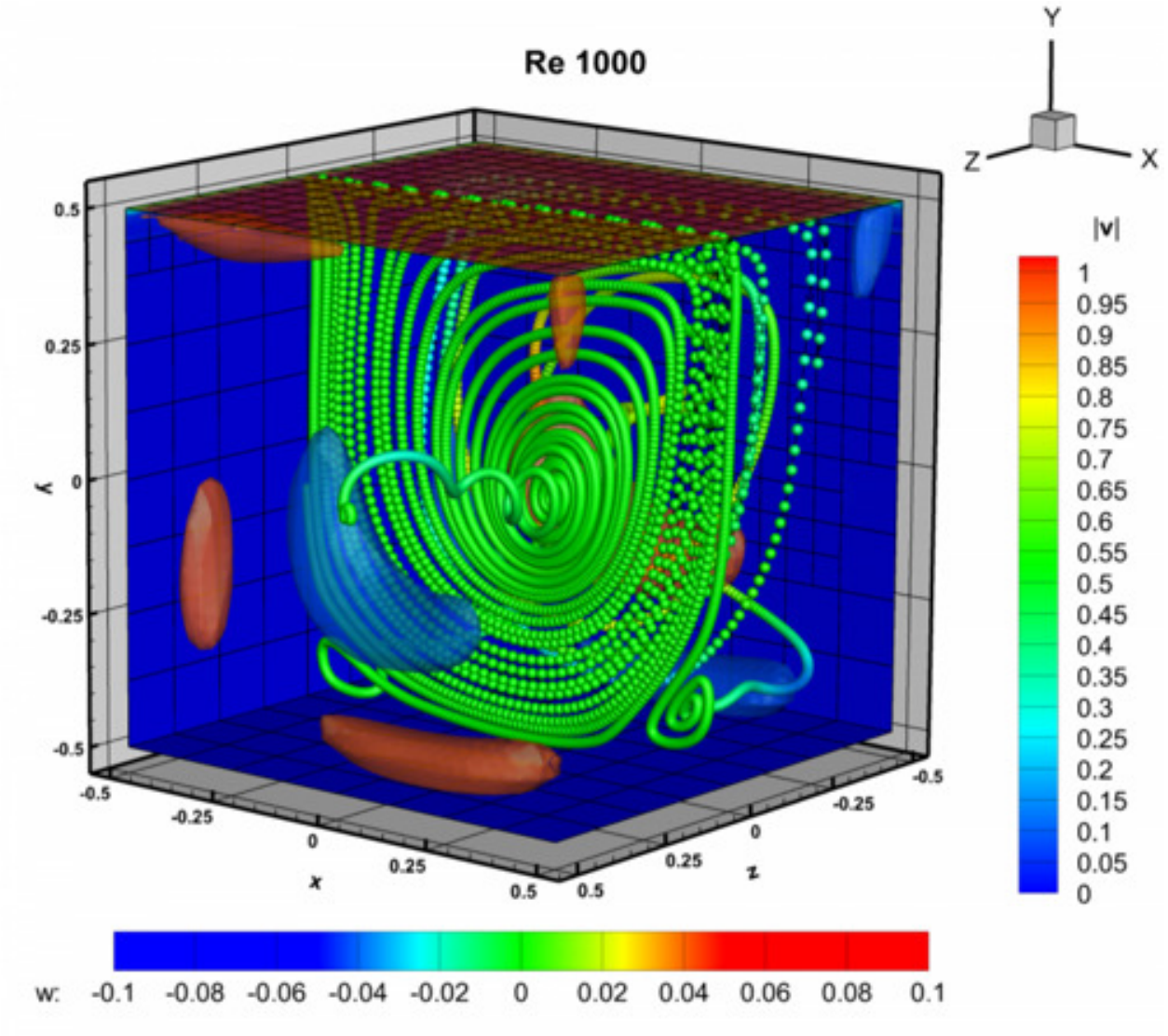} 
\caption{The numerical solution obtained for the three dimensional lid-driven cavity problem compared with, at the left, the numerical results of \cite{Albensoeder2005,Ku1987} at  Reynolds numbers Re=$100$ (top) and Re=$1000$ (bottom) using $8\times 8\times 8$ elements. The results have been obtained with our staggered semi-implicit DG-$\p_{4}$ 
method with AMR. In the 3d-view on the right, streamlines are colored with the $w$ velocity magnitude, $w$-peak iso-surfaces are shown, the boundary-slices are colored by the velocity magnitude $|\textbf{v}|$ together with the \emph{main} AMR-grid. }\label{fig:LDCavity3Da}
\end{figure}

\begin{figure} 
\centering 
			\includegraphics[width=0.4\textwidth]{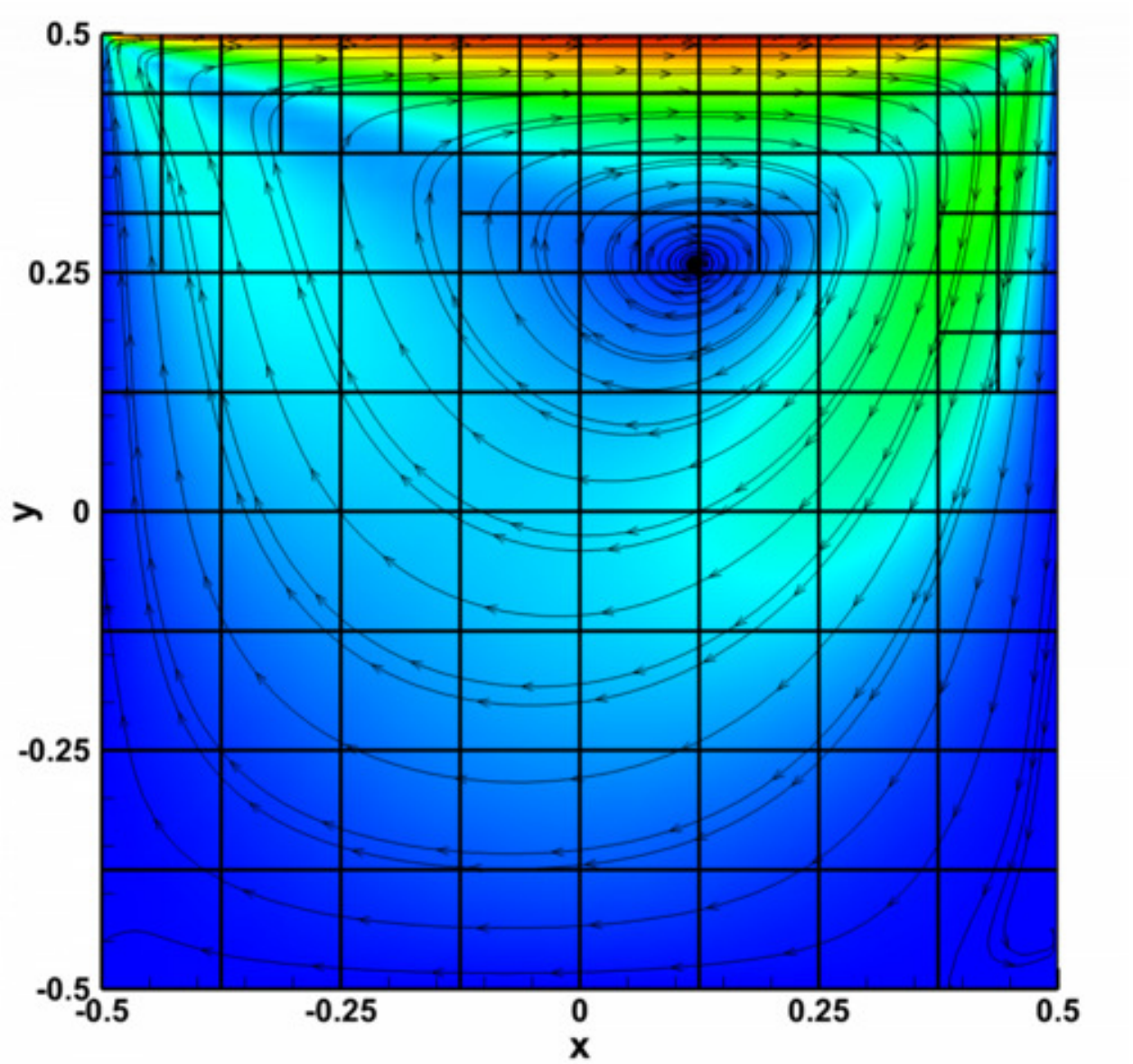}
			\includegraphics[width=0.4\textwidth]{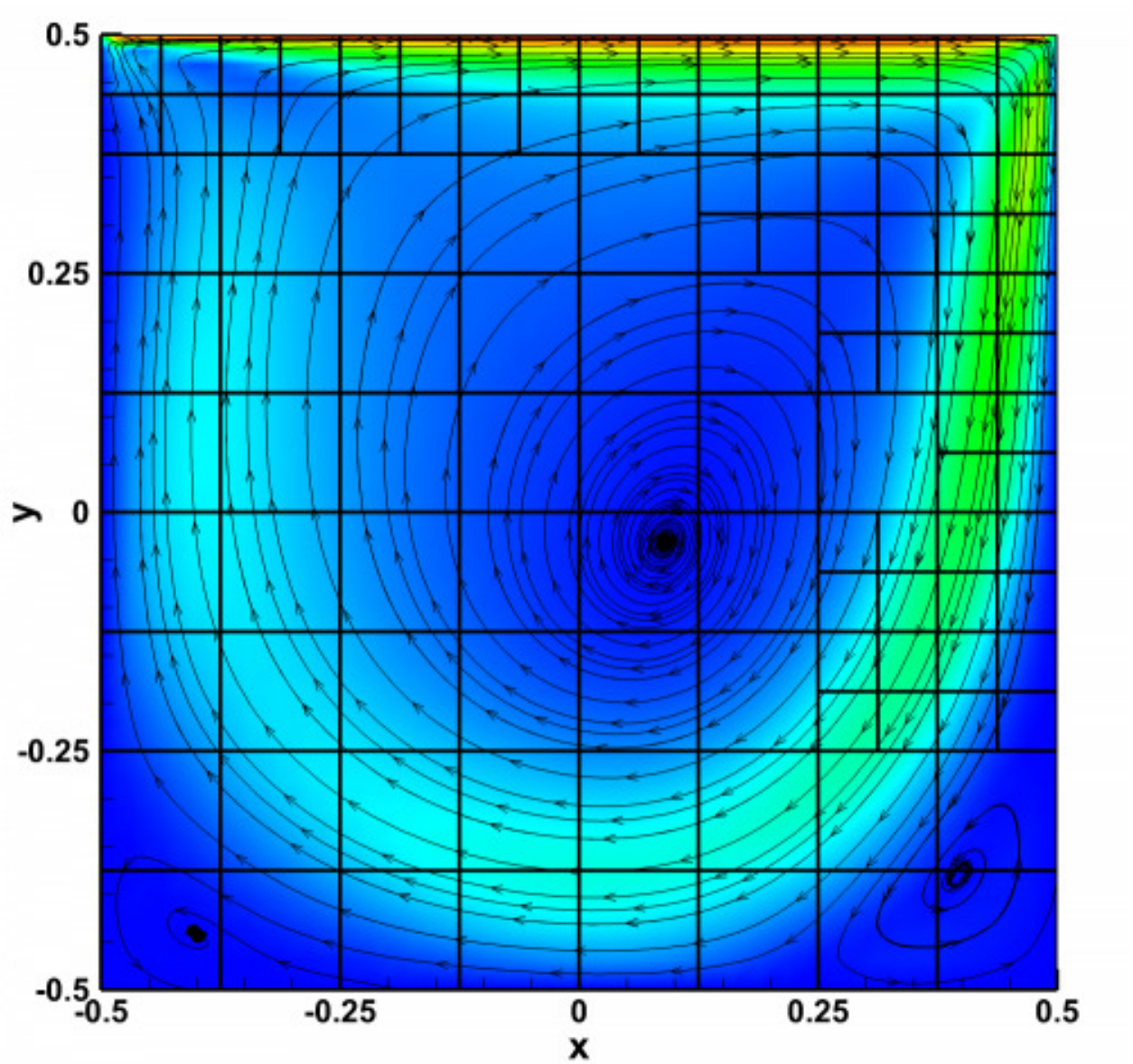}\\
			\includegraphics[width=0.4\textwidth]{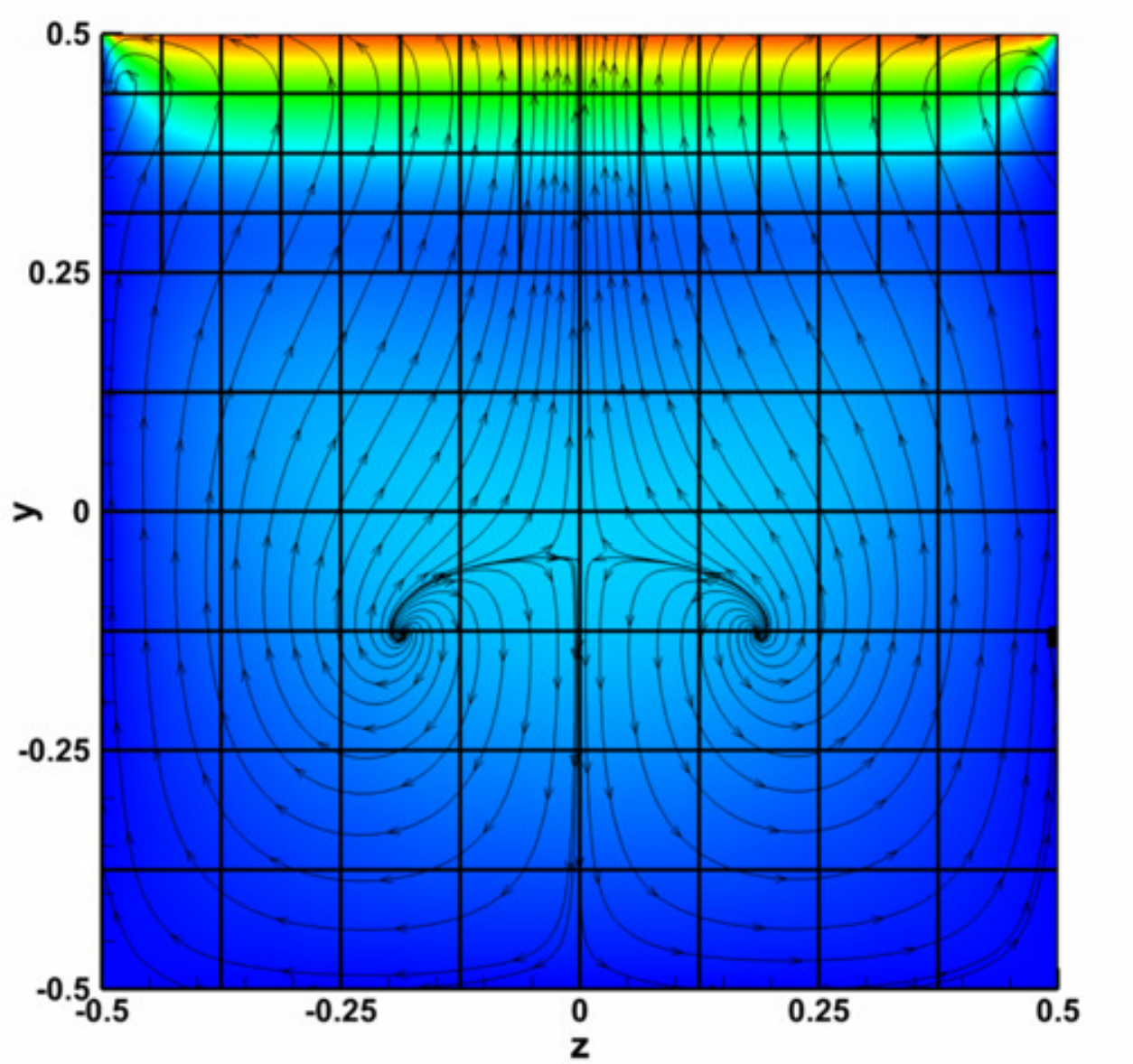}
			\includegraphics[width=0.4\textwidth]{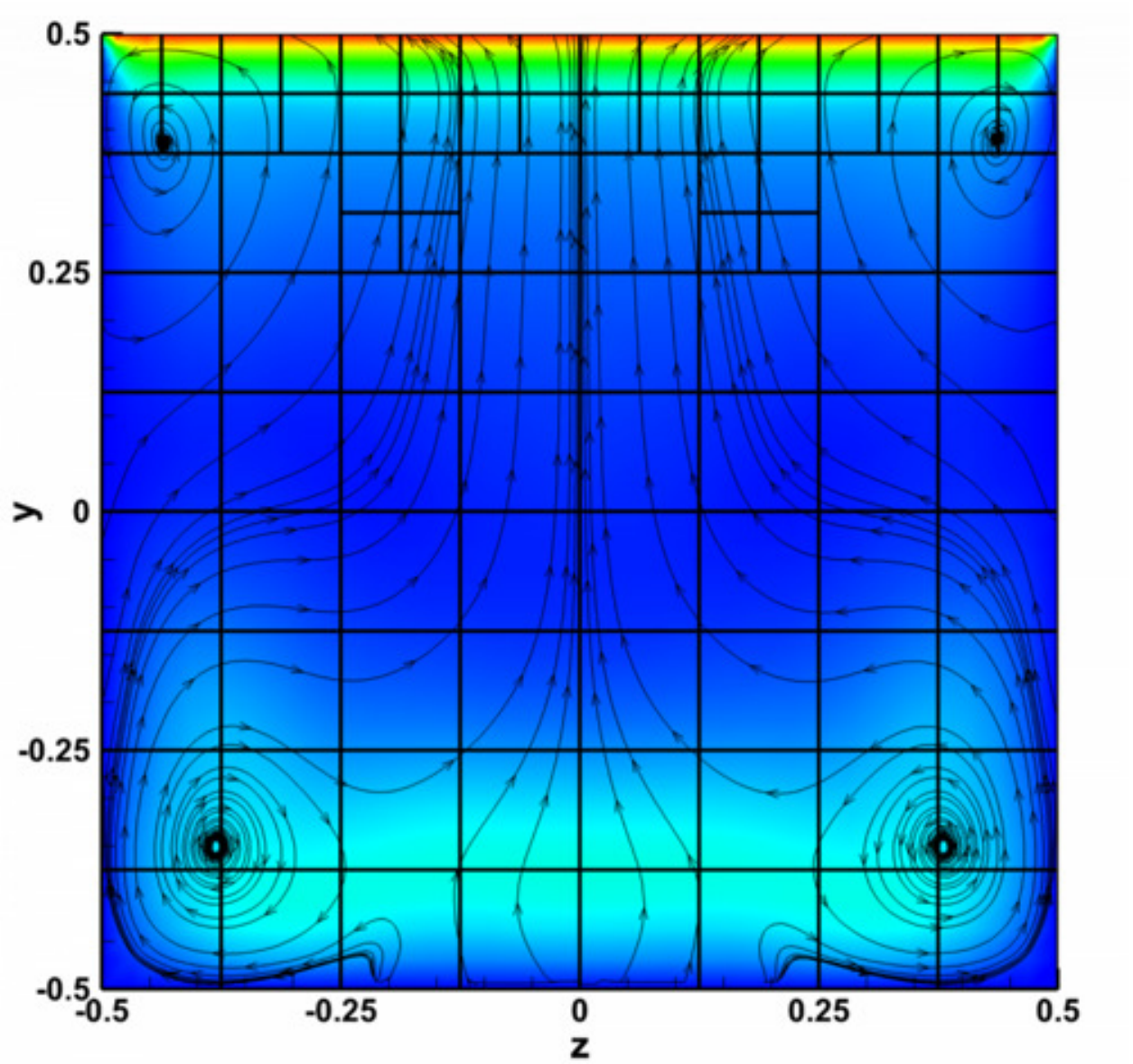}\\
			\includegraphics[width=0.4\textwidth]{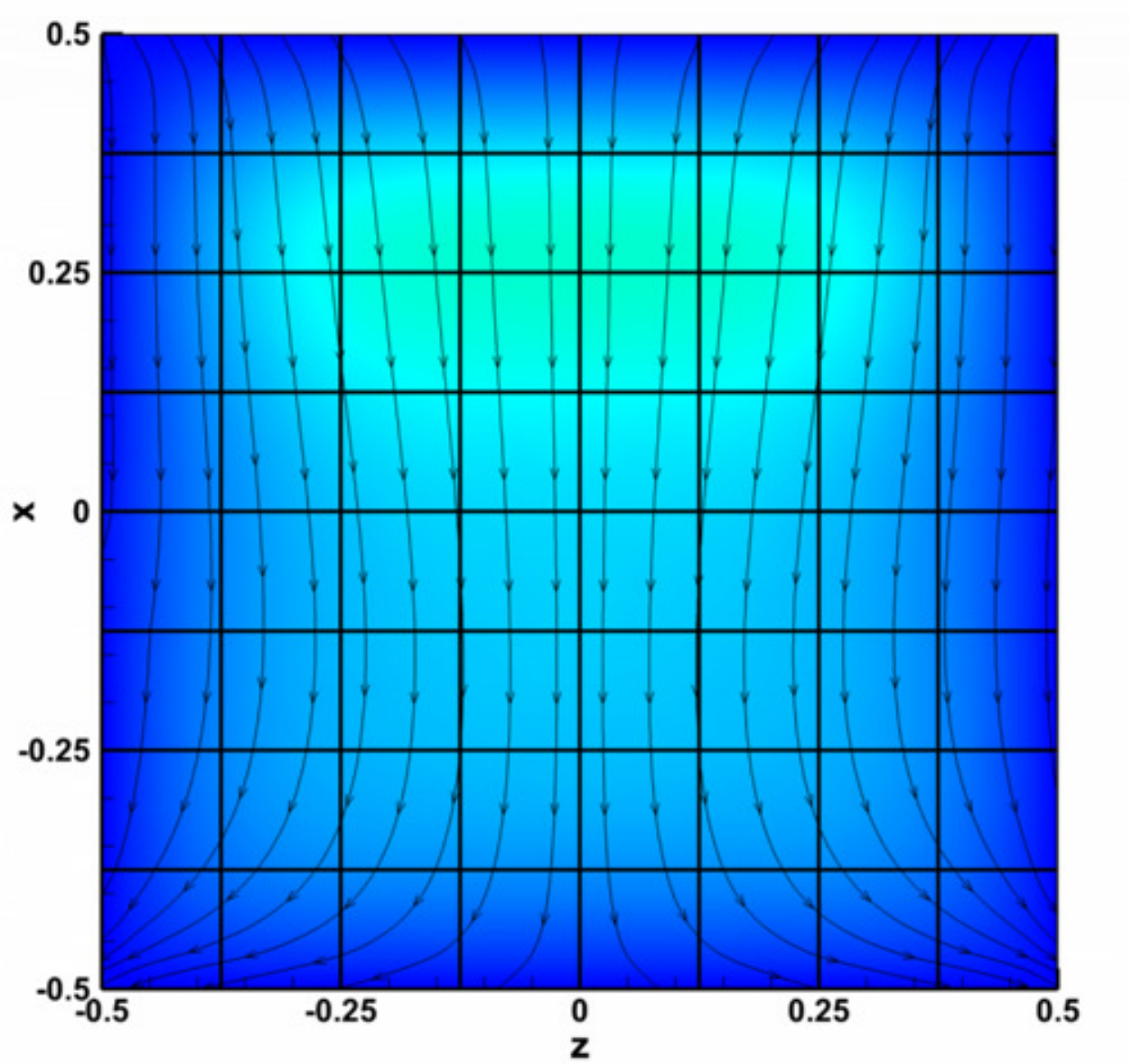}
			\includegraphics[width=0.4\textwidth]{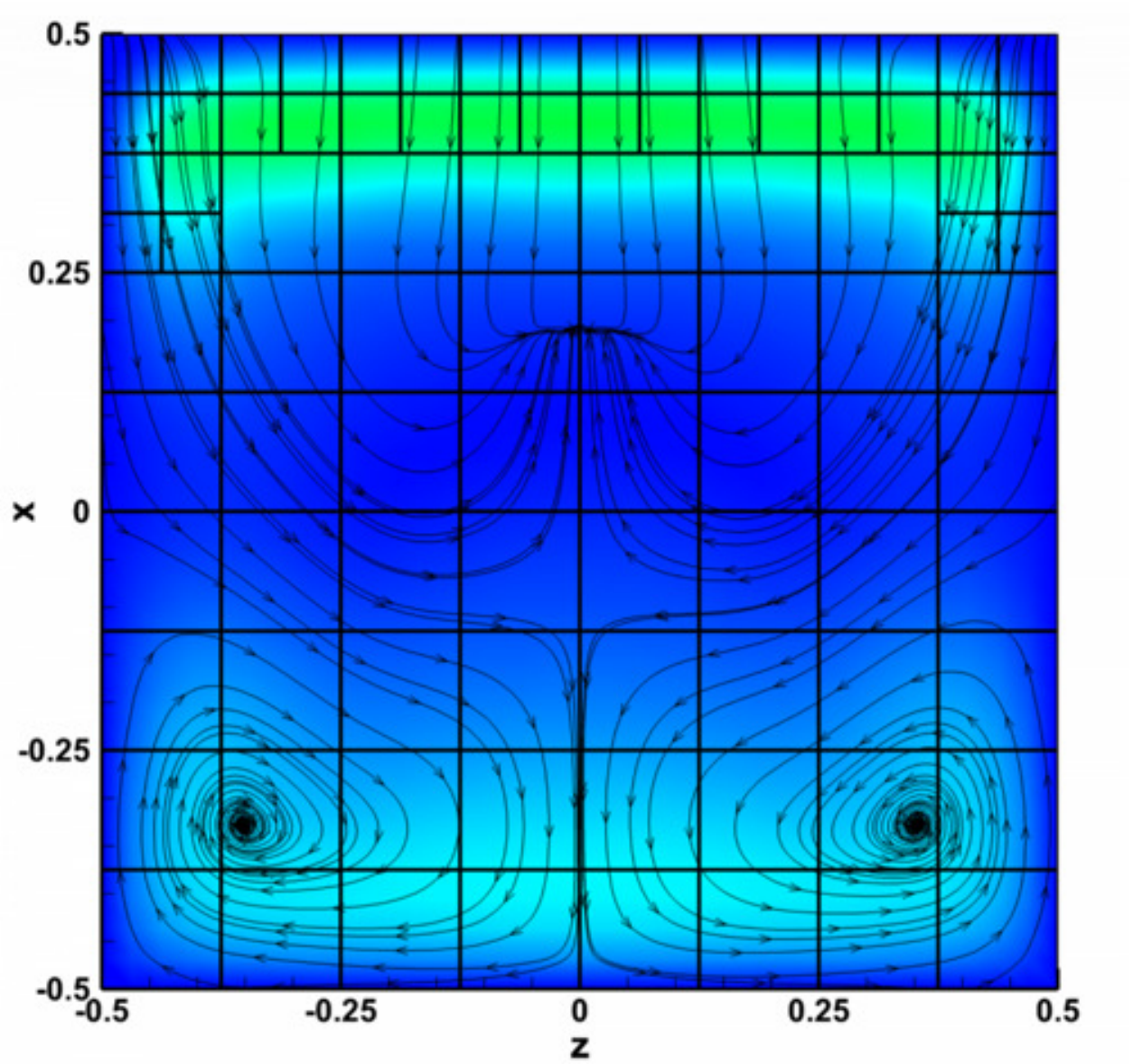}
\caption{The numerical solution obtained for the three dimensional lid-driven cavity problem  interpolated along the three orthogonal planes $x-y$, $z-x$ and $z-y$ from the top the bottom, at  Reynolds numbers Re=$100$ (left) and Re=$1000$ (right) using $8\times 8\times 8$ elements. Tangential streamlines and the velocity magnitude $|\textbf{v}|$ are depicted. The results have been obtained with our staggered semi-implicit AMR-DG-$\p_{4}$ 
method.}\label{fig:LDCavity3Db}
\end{figure}

In three space dimensions the fluid flow becomes much more complicated due to the appearance of secondary recirculations which make the flow field fully three-dimensional. Figures \ref{fig:LDCavity3Da}-\ref{fig:LDCavity3Db} show the respective numerical results obtained at Reynolds numbers $Re=100$ and $Re=1000$ with our \SIDG-$\p_4$ scheme, corresponding to a total number of degrees of freedom 
of $\ndof=5^3=125$ per space-element. The velocity profile interpolated along the vertical and horizontal axes, i.e. $(x,y,z) \in[-0.5,0.5]\times\{0\} \times\{0\}$ and $(x,y,z) \in\{0\} \times[-0.5,0.5]\times\{0\}$, is shown  and compared with the reference solutions of \cite{Albensoeder2005,Ku1987} next to the corresponding three-dimensional view of the fluid flow. A very good match 
between our numerical results and the reference data is obtained. The three dimensional recirculation together with the active AMR grid are highlighted in Figure \ref{fig:LDCavity3Da}, where  
the numerical solution is interpolated along the three orthogonal panes $x-y$,  $y-z$ and $x-z$. In this case, the physical domain has been discretized on the zeroth level using only $8^3=512$  
elements. For the AMR framework we set $\ell_{\text{max}}=1$ and $\err=2$. 
		
		\subsection{2D backward-facing step}
		
\begin{figure} 
\centering 
			\includegraphics[width=0.75\textwidth]{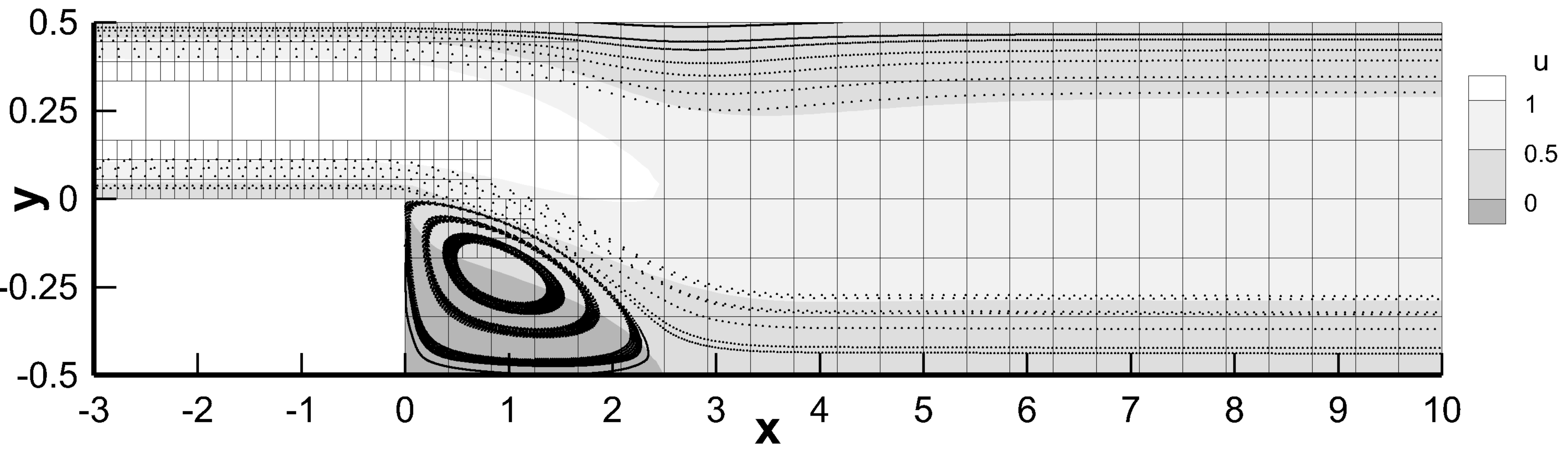}\\
			\includegraphics[width=0.75\textwidth]{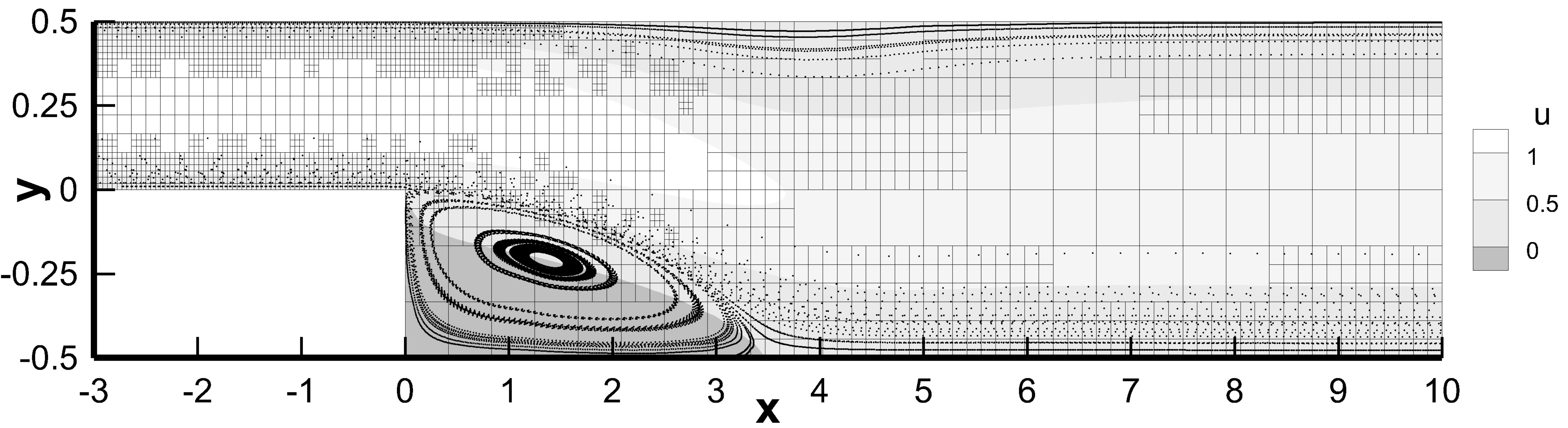}\\
			\includegraphics[width=0.75\textwidth]{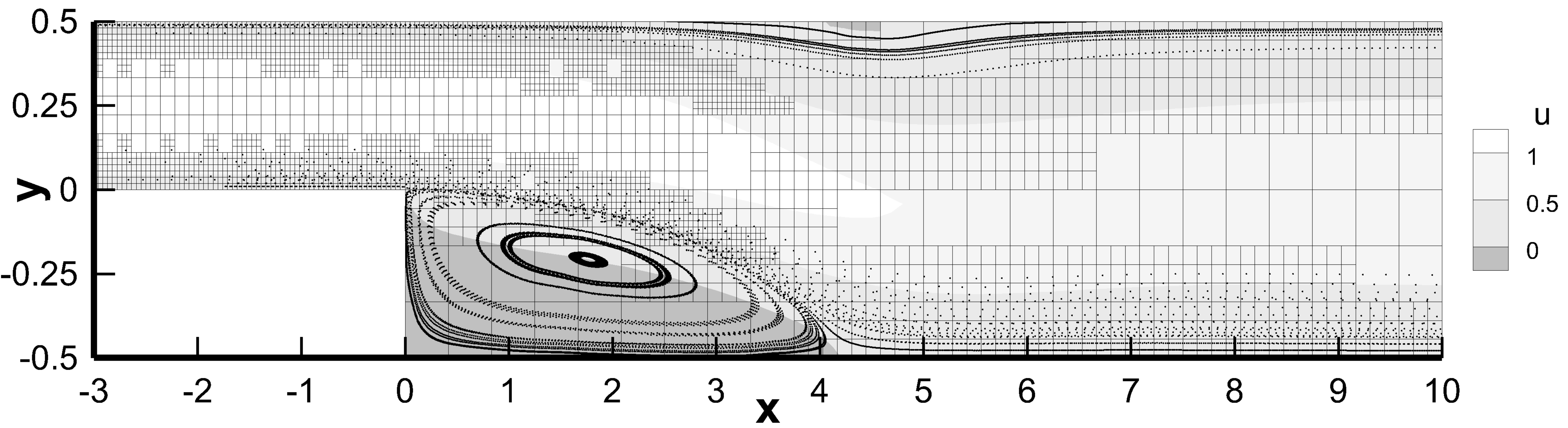}\\
			\includegraphics[width=0.75\textwidth]{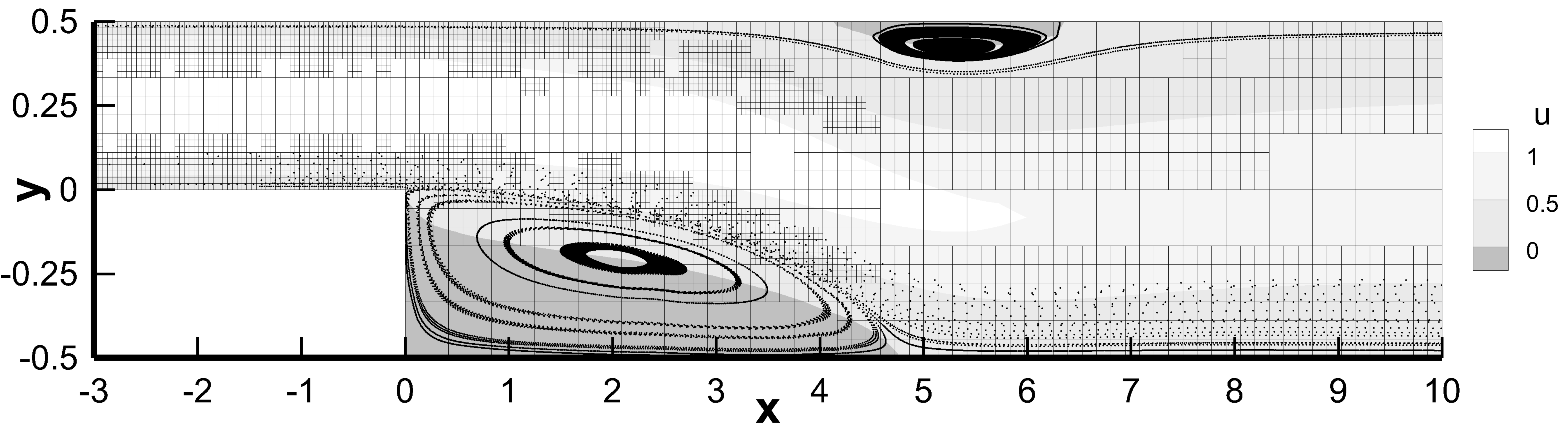}\\
			\includegraphics[width=0.75\textwidth]{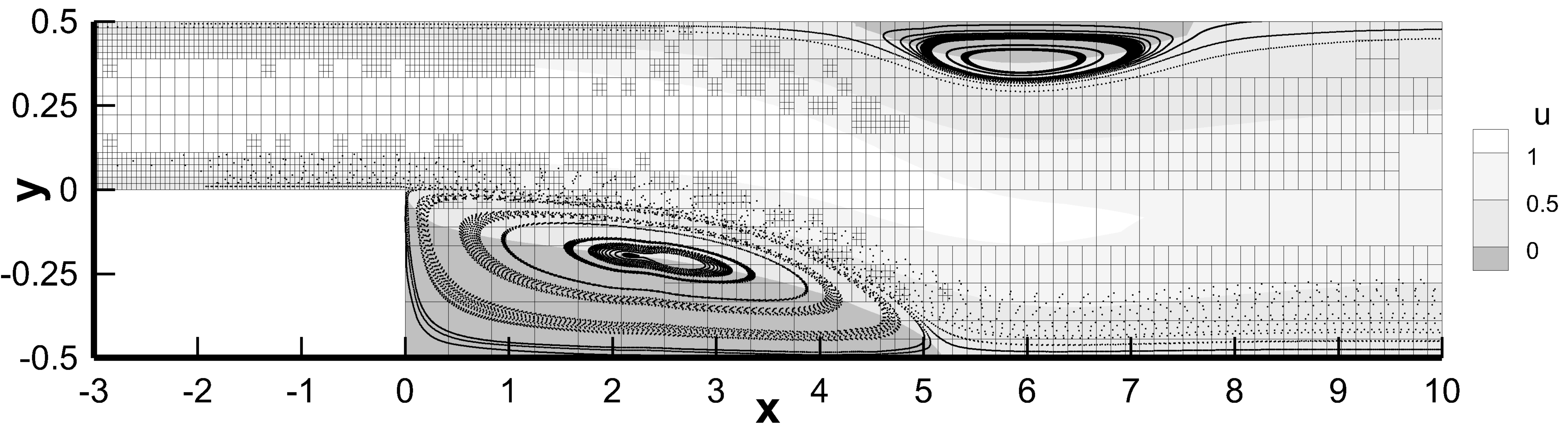}\\
			\includegraphics[width=0.75\textwidth]{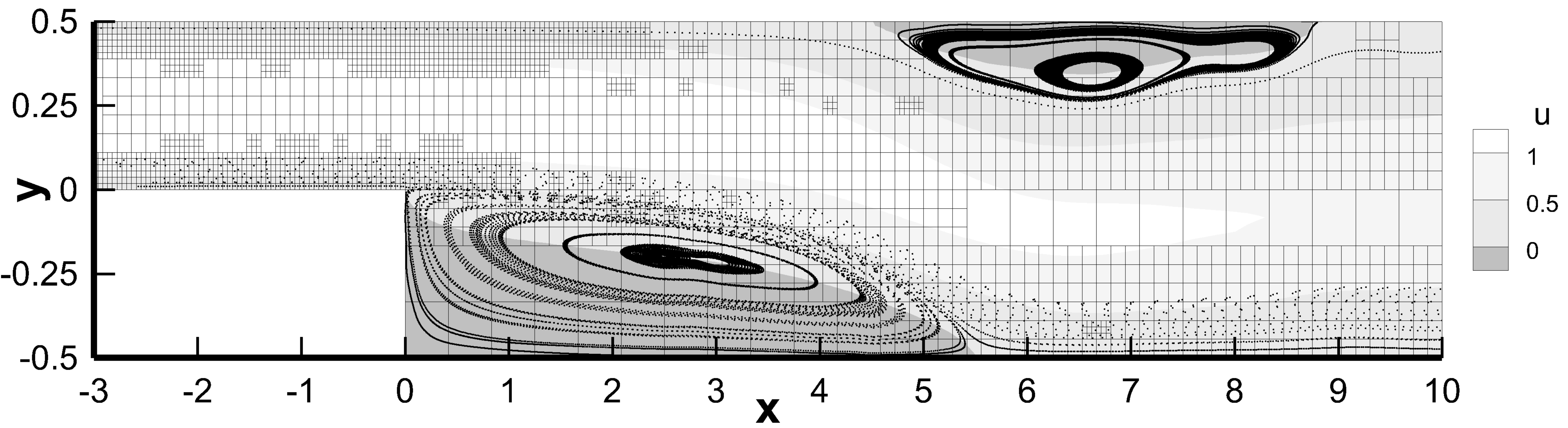}
\caption{Numerical solution obtained for the 2D backward facing step problem with the staggered semi-implicit AMR-DG-$\p_{4}$ method at different Reynolds numbers, from top to bottom, respectively:  
Re=$200$, Re=$300$, Re=$400$, Re=$500$, Re=$600$, and Re=$700$. 
}
\label{fig:BFStep2D}
\end{figure}

\begin{figure} 
\centering 
			\includegraphics[width=0.48\textwidth]{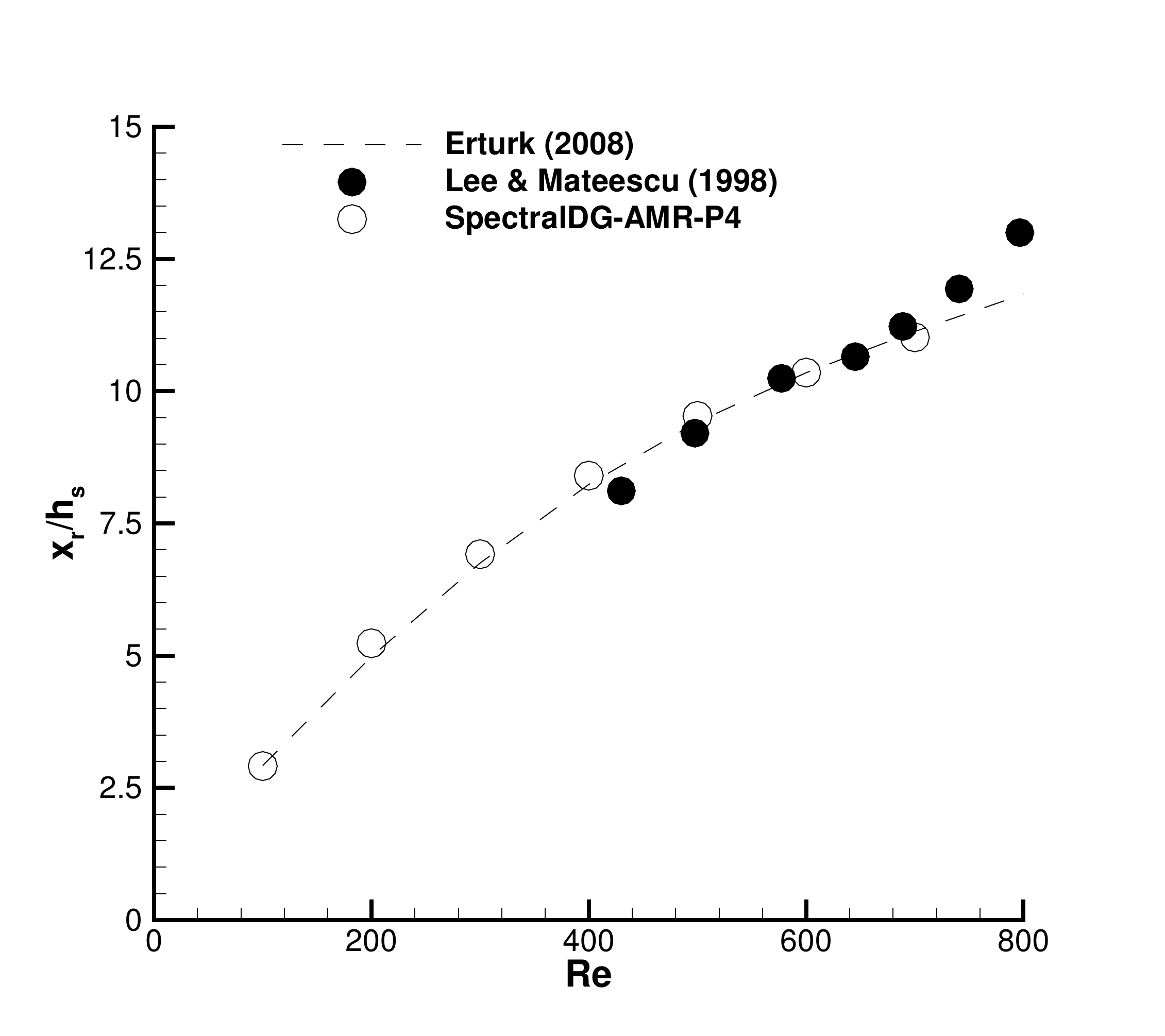}
\caption{Total lengths of the main recirculations obtained with the staggered semi-implicit spectral AMR-DG-$\p_{4}$ method for the 2D backward facing step problem next to the 2D 
numerical results of \cite{Erturk2008} and the experimental measurements of \cite{LeeMateescu1998} at different Reynolds numbers, up to $Re=800$.}\label{fig:BFStep_data} 
\end{figure}

		In this test a constant fluid flow is applied at the entrance of a two dimensional rectangular duct with a backward-facing step, characterized by an expansion ratio of $ER=h_{\text{out}}/h_{\text{in}}=2$ at $x=0$ and a backward facing step of height $h_\text{s}=0.5$. A sufficiently long axial domain is chosen to be $x\in[-10,20]$, in order to avoid spurious effects from the boundaries. 
		Indeed, in our set-up a constant velocity profile is applied at the entrance and a constant pressure $p=0$ is imposed at the outflow. Similar to the previous test, the backward-facing step became a  typical benchmark problem because it allows to cover a large domain of Reynolds numbers and flow regimes by simply varying the kinematic viscosity. The physical domain is discretized with elements of 
		size $\Delta x_{\ell=0}=5/12$ and $\Delta y_{\ell=0} =1/6$ on the coarsest main grid. For the AMR framework we use a maximum number of refinement levels $\ell_{\text{max}}=2$ and a refine factor of 
		$\err=3$, corresponding to a mesh size on the finest grid of $\Delta x_{\ell=2} =5/108$ and $\Delta y_{\ell=2} = 1/54$. A main recirculation zone is generated next to the step at the bottom, even 
		at low Reynolds numbers. Then, at higher Reynolds numbers, the fluid flow becomes more complicated and new recirculation zones appear in the duct. Figure \ref{fig:BFStep2D} shows the numerical results obtained ad different Reynolds numbers within $Re \in [100,800]$,  where the main recirculations together with the active AMR grid are highlighted. Corresponding reference solutions exist for the axial  extent of the main circulation zone, available in \cite{LeeMateescu1998} (experimental study) and \cite{Erturk2008} (numerical study), and are compared with the computed results obtained through our adaptive \SIDG-$\p_4$ scheme in Figure \ref{fig:BFStep_data}. As it is shown, our numerical results match the two-dimensional numerical reference data of \cite{Erturk2008} very well, as well as the 
		experimental data of \cite{LeeMateescu1998} in the low Reynolds number regime. Some differences with the experimental results arise for higher Reynolds numbers. The correct interpretation of  these discrepancies lies in the three-dimensional behaviour of the fluid flow at higher Reynolds numbers due to the side-wall boundary effects, see \cite{Tylli2002,Armaly1983,Mouza2005,Rani2007}. 
		
 
		\subsection{2D double shear layer} 
	
A classical two dimensional test that is well-suited for demonstrating the resolution of high order methods and the usefulness of adaptive grids is the double shear layer. In this work, 
the so called '\emph{thin}' double shear layer is chosen to be the initial and unstable state at $t=0$ given by 
\begin{align}
u  = \left\{ \begin{array}{rc} u_0 \tanh\left( (y - 0.5)/\delta \right)  & y > 0\\ - u_0 \tanh\left( (y + 0.5)/\delta \right)  & y \le 0 \end{array}\right.
\end{align} 
with a perturbation in the vertical velocity component in the vicinity of the shear layer that reads 
\begin{align}
v  = \left\{ \begin{array}{rc} v_0 \sin\left( 2\pi x \right) e^{-(y-0.5)^2/2 \sigma^2}  & y > 0\\ - v_0 \sin\left( 2\pi x \right) e^{-(y+0.5)^2/2 \sigma^2 }  & y \le 0 \end{array}\right.
\end{align}
The rectangular spatial domain $\Omega=[-0.5,0.5]\times[-1,1]$ has been discretized by a grid composed of $20 \times 40$ elements at the coarsest level $\ell=0$. We use 
up to $\ell_{\text{max}}=2$ refinement levels and a refine factor of $\err=2$ in this test problem. 
In the present test the chosen parameters are $\delta = 10^{-2}$ for the shear layer thickness, $u_0=10$, $v_0 = 0.5$ and $\sigma^2 = 0.05$ for the variance of the Gaussian perturbation. The example
has been run with a kinematic viscosity of $\nu=10^{-4}$.

\begin{figure} 
\centering 
			\includegraphics[width=0.44\textwidth]{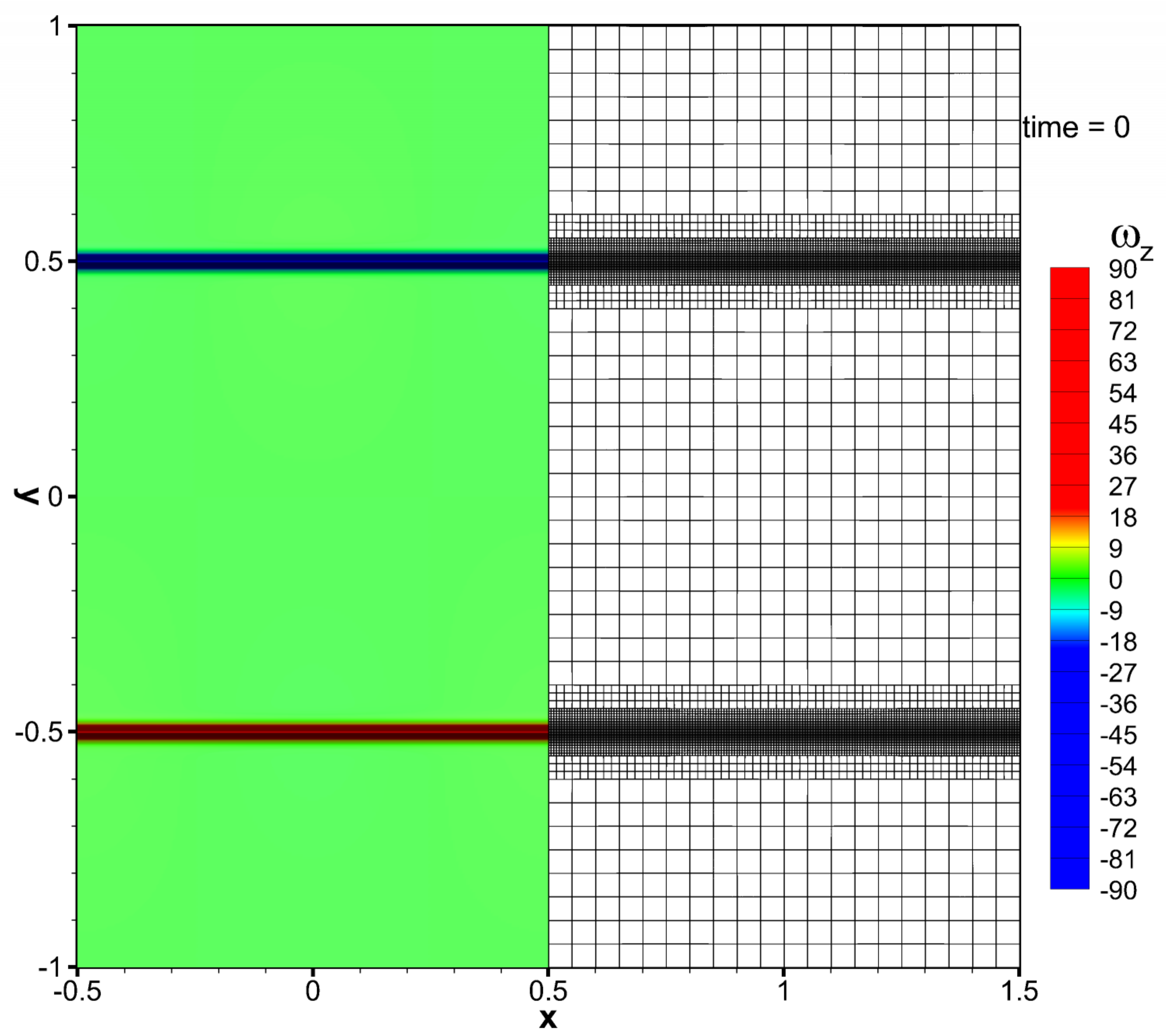}\;\includegraphics[width=0.44\textwidth]{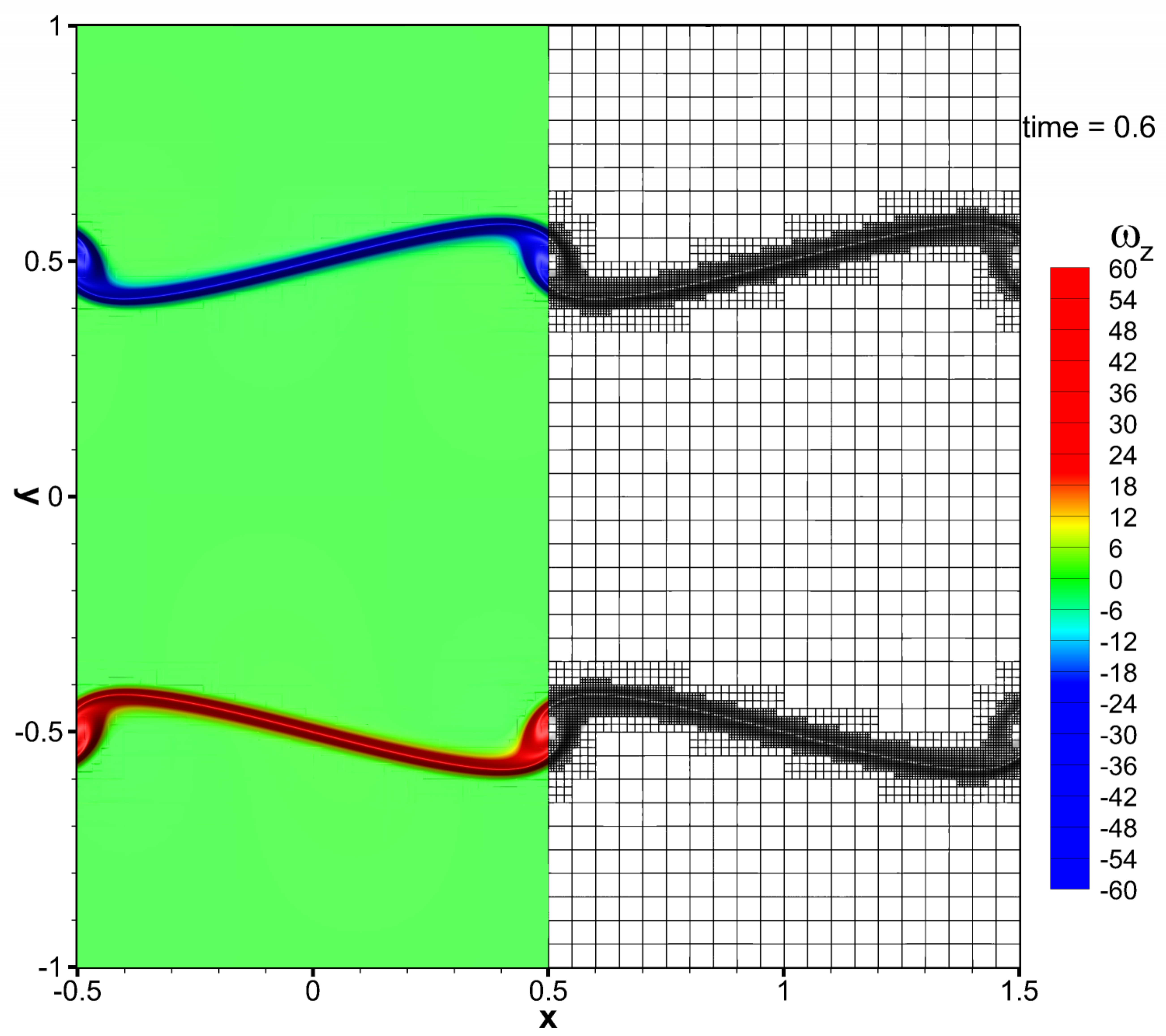}\\
			\includegraphics[width=0.44\textwidth]{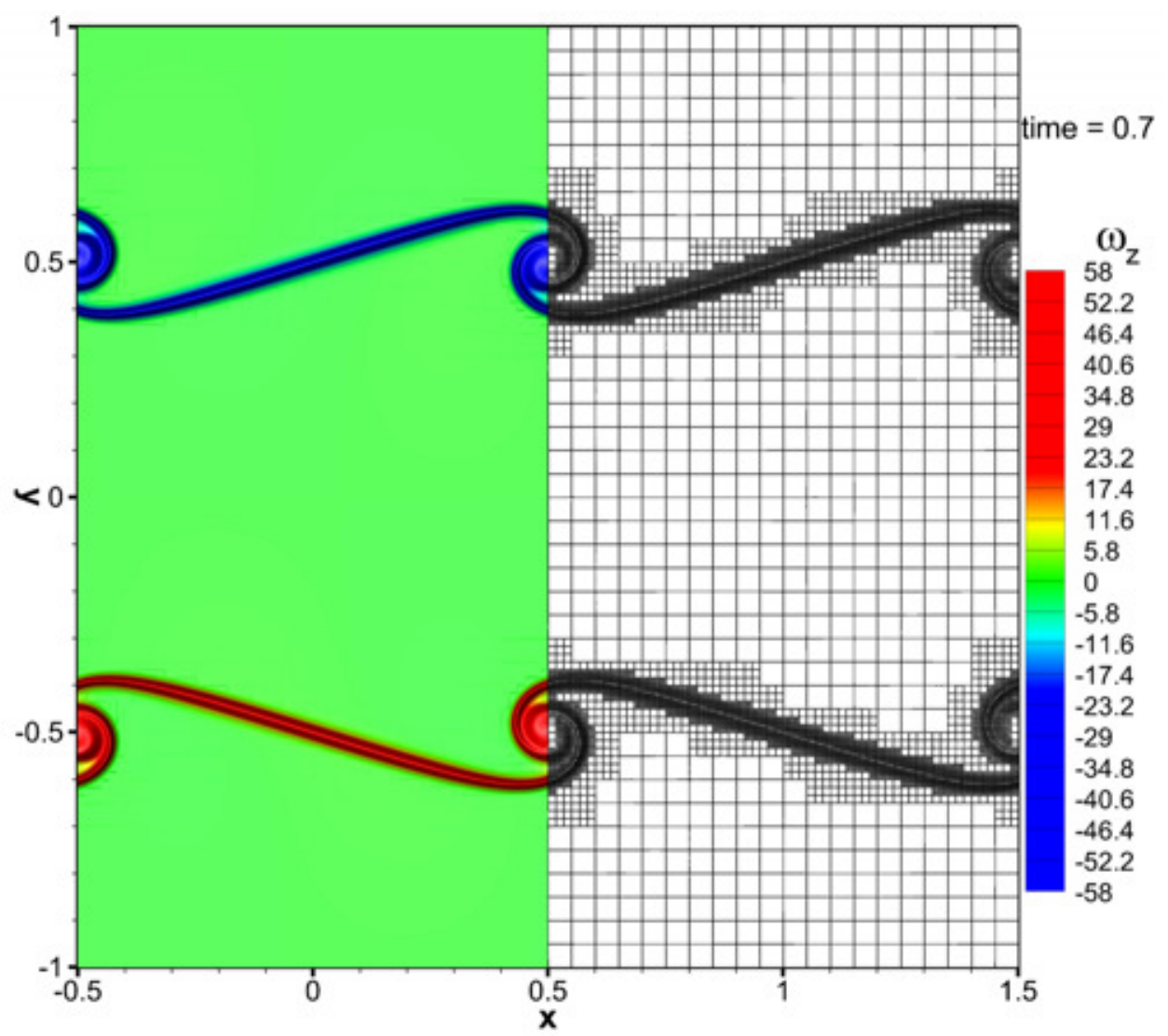}\;\includegraphics[width=0.44\textwidth]{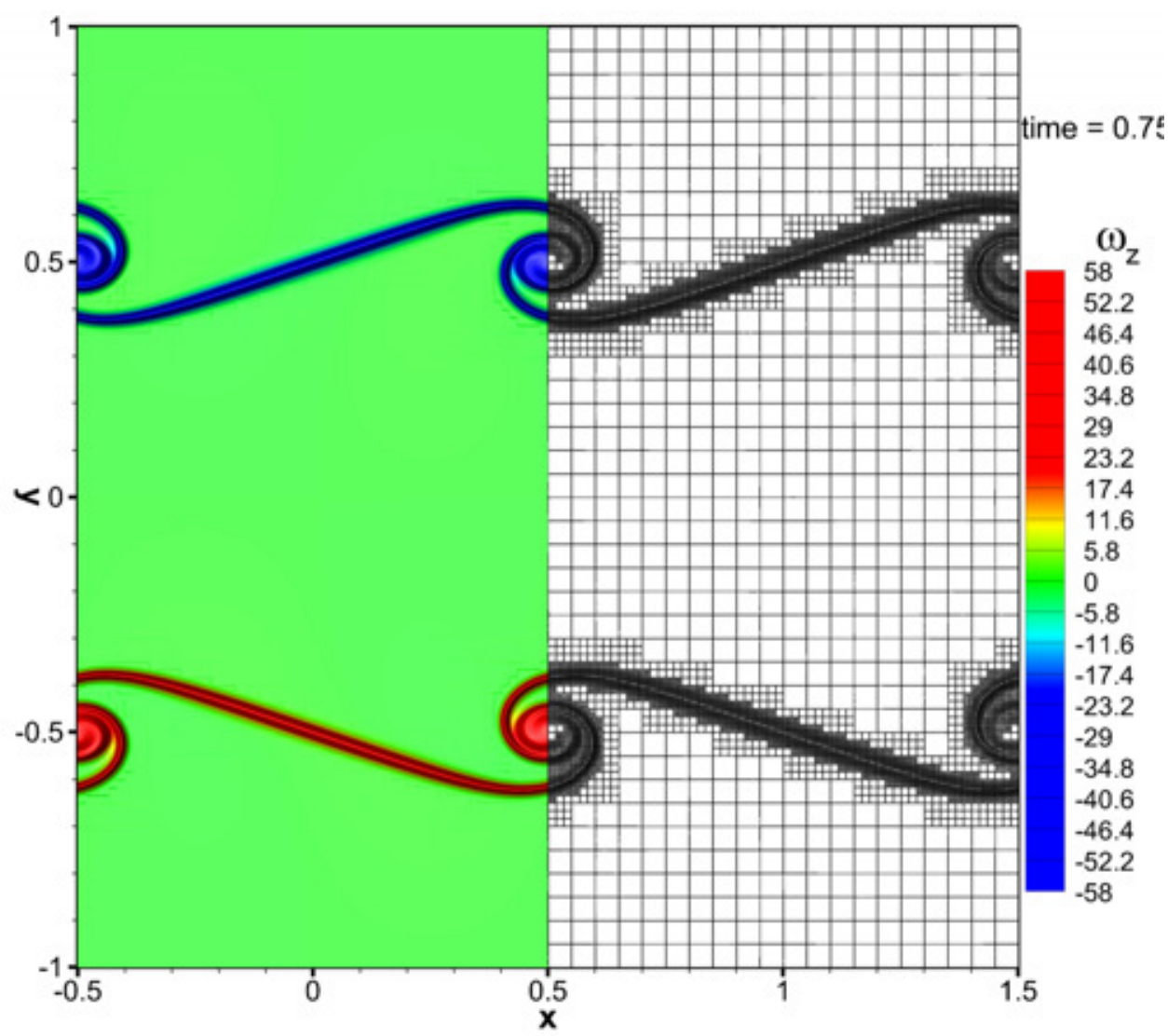}\\
			\includegraphics[width=0.44\textwidth]{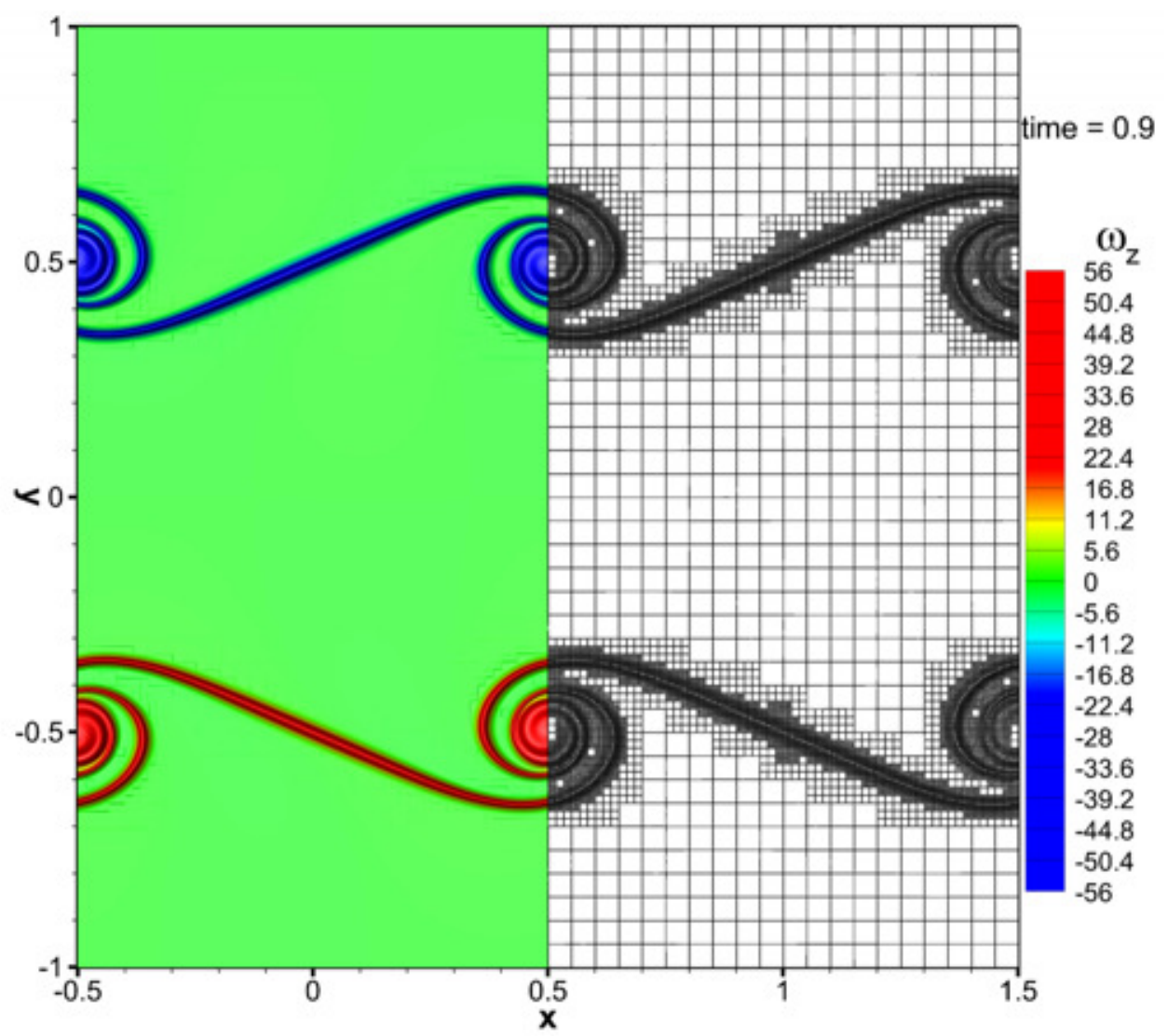}\;\includegraphics[width=0.44\textwidth]{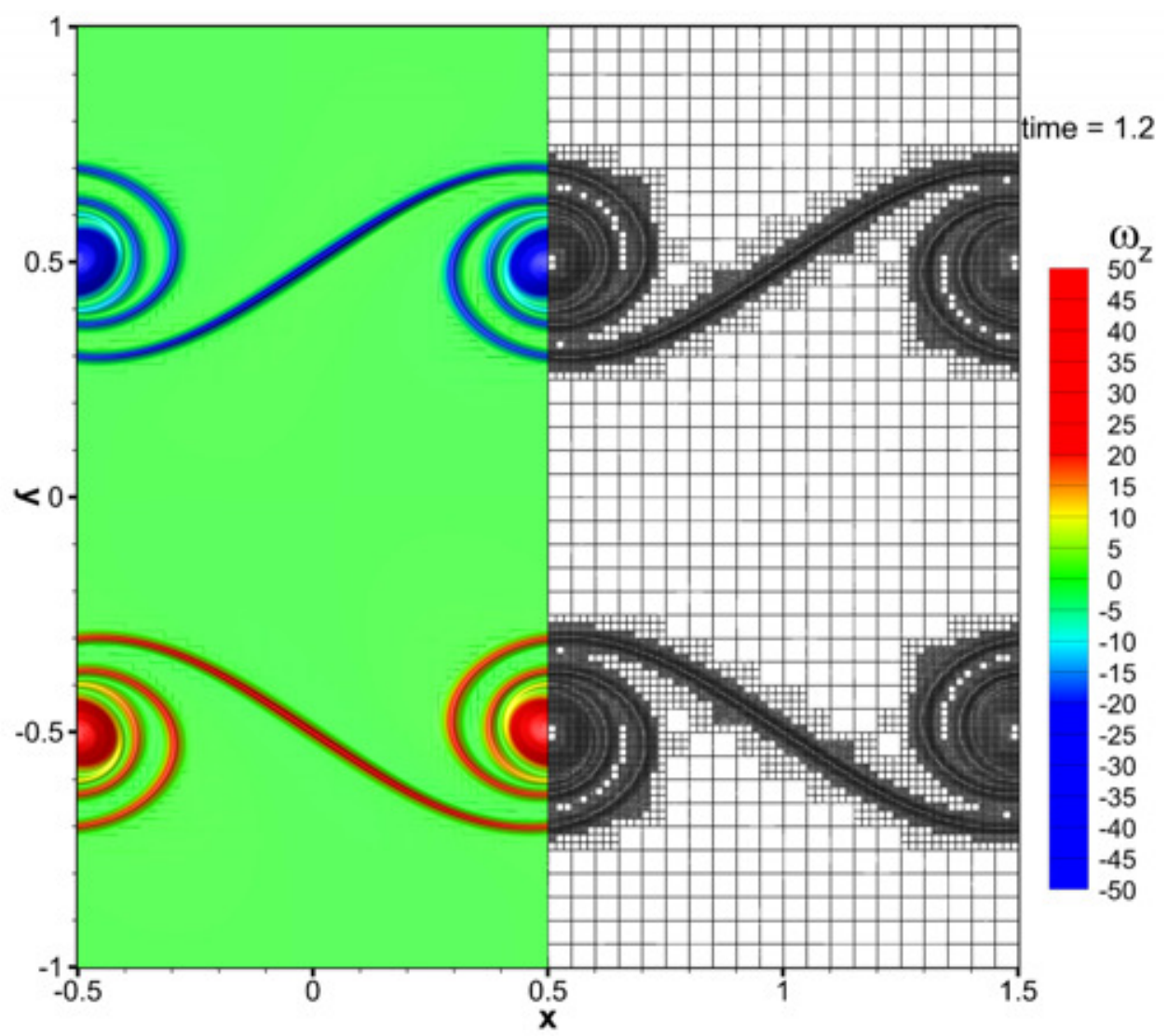} 
\caption{The numerical solution obtained for the two dimensional \emph{'thin'} double shear layer problem at different time steps. The staggered \SIDG-$\p_{9}$ method has been used together with a refinement factor $\mathfrak{r}=3$ and  $\ell_\text{max}=2$ refinement levels. Every figure depicts the vorticity field $\omega_z$ next to the respective AMR grid.}\label{fig:DSL2D}
\end{figure}

As mentioned in \cite{Brown1995,Minion1997} a not sufficiently accurate solution of the flow field may cause spurious oscillations that arise in different locations along the shear layers. 
Figure \ref{fig:DSL2D} shows the time evolution of the z-component of vorticity $\omega_x$ obtained with our $\SIDG$-$\p_9$ scheme, next to the active main AMR grid. The computed results show 
to be in agreement with previously published results in literature. Moreover, in this test the space-time AMR is shown to give major benefits, resulting in a very high resolution obtained with 
a still rather coarse mesh at $\ell=0$. Notice that the coarsest level corresponds to a total number of $N_0=800$ elements, a total number of degrees of freedom of $N^{\text{dof}}_0 = 80'000$ 
and a characteristic mesh size of $h_0=1/20=0.05$, while a uniform grid on the finest level corresponds to $N_2=64'800$, $N^{\text{dof}}_2 = 6'480'000$ and $h_2=h_0/9$. No spurious oscillations are generated 
in our simulation, which means that the non-conforming elements that appear in the AMR framework are treated properly by our numerical method.

		\subsection{3D Taylor-Green vortex problem}
\begin{figure} 
\centering 
			\includegraphics[width=0.6\textwidth]{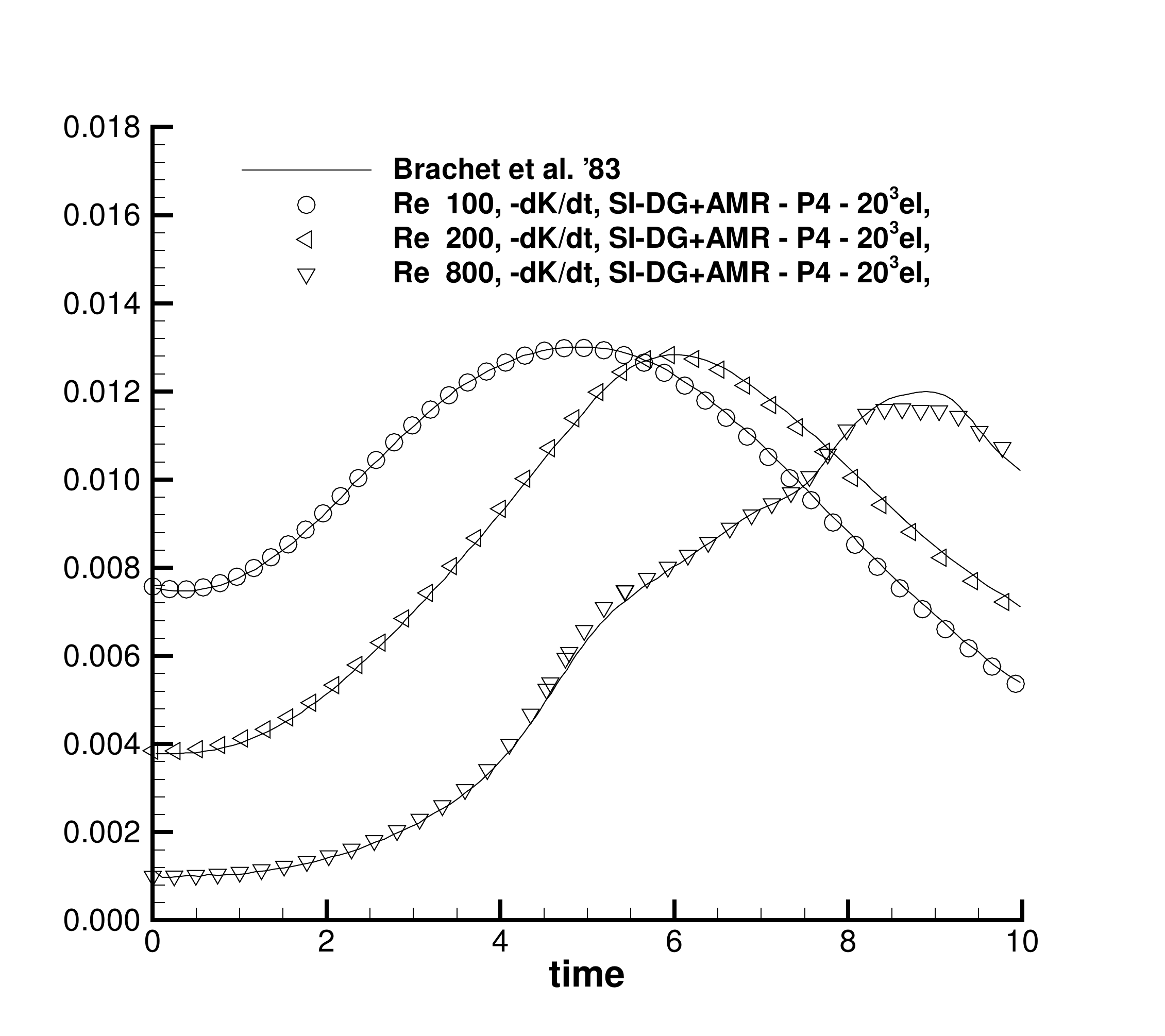}
			\caption{Time evolution of the kinetic energy dissipation rate $\epsilon(t)$ obtained with our $\SIDG$-$\p_4$ scheme  
			at   Reynolds numbers $Re=100$, $200$ and $800$. The DNS reference solutions of Brachet et al. \cite{Brachet1983} are plotted as continuous lines.} \label{fig:TG3D}
\end{figure}

\begin{figure} 
\centering 
			\includegraphics[width=0.32\textwidth]{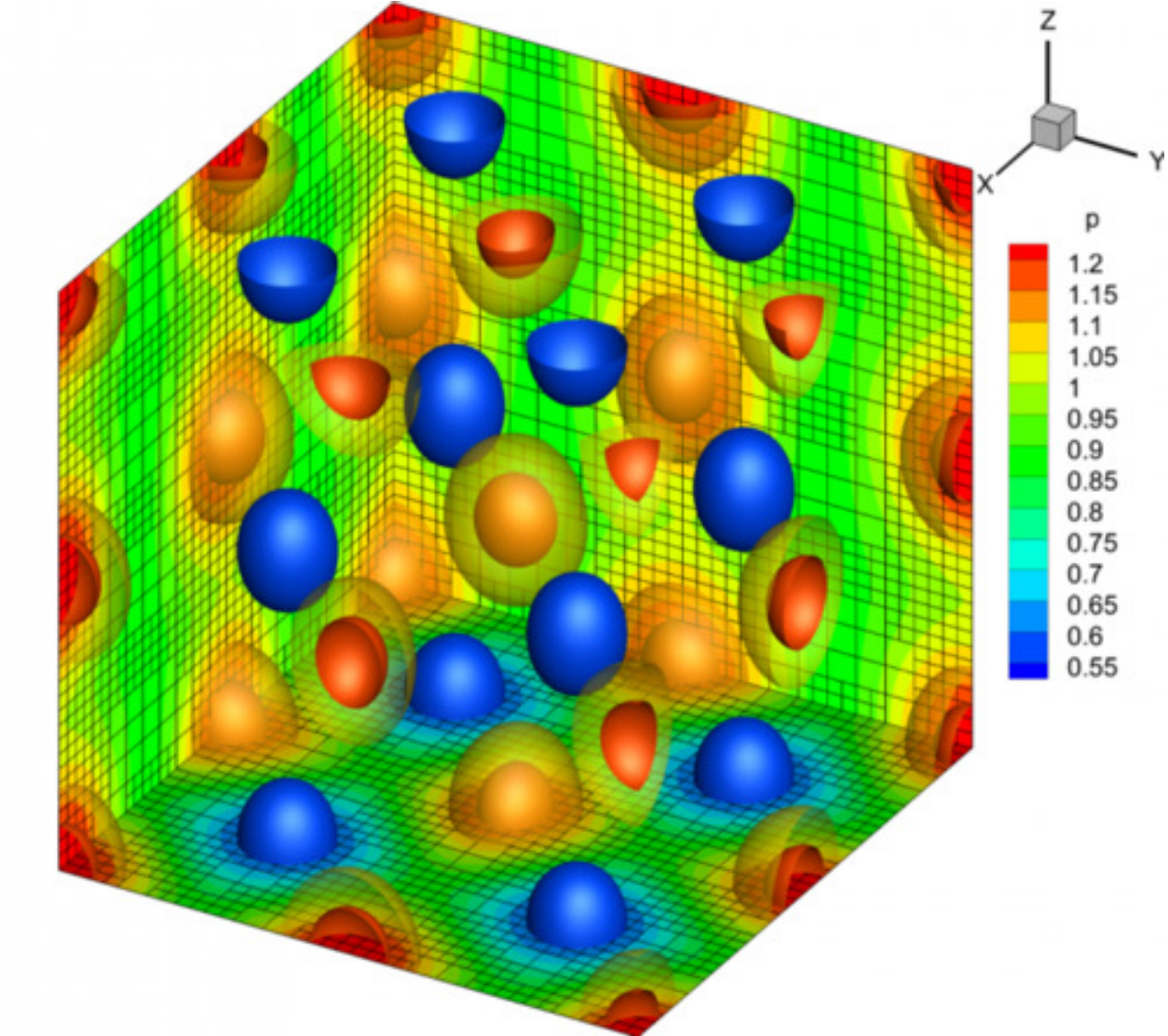}
			\includegraphics[width=0.32\textwidth]{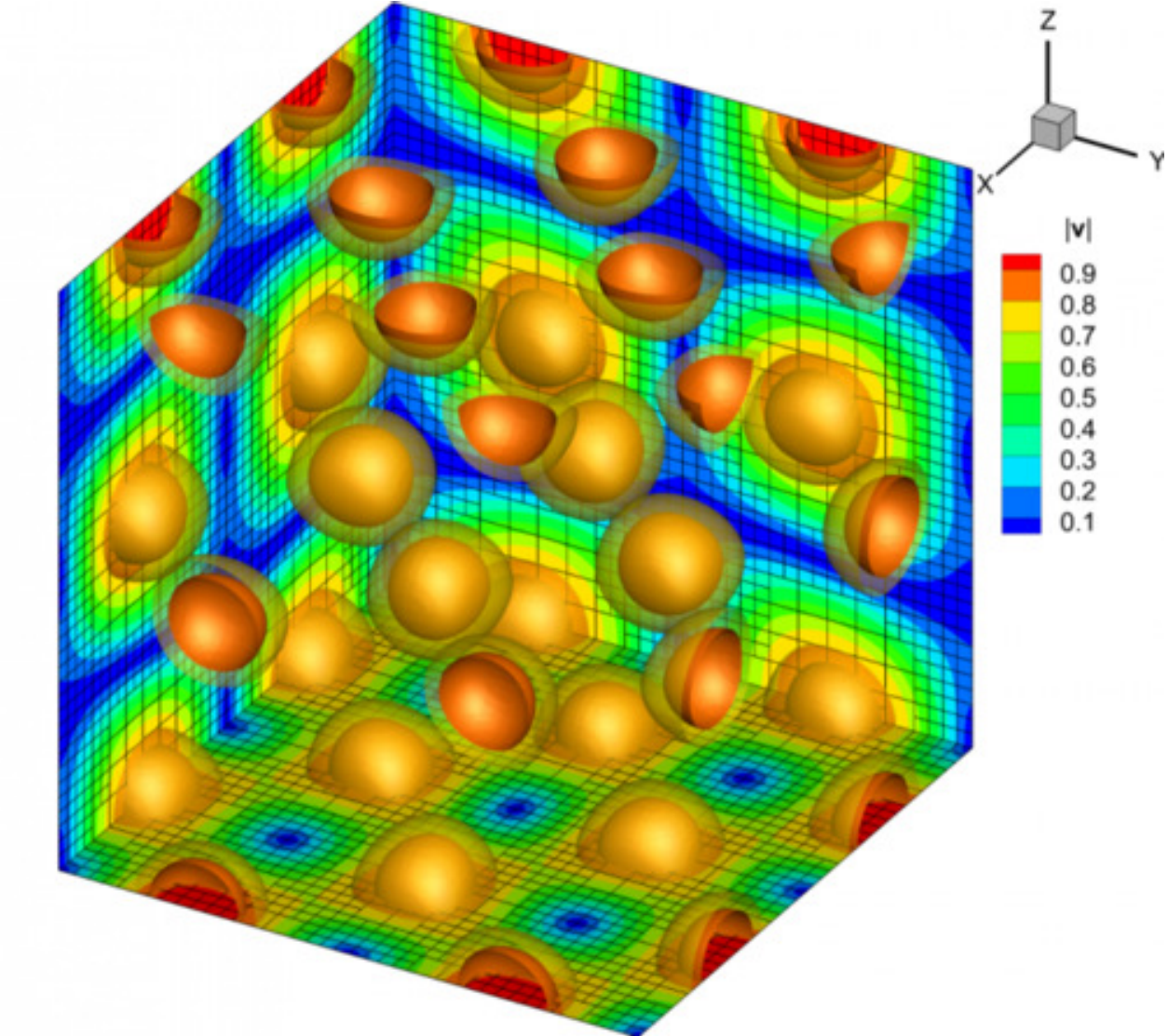}
			\includegraphics[width=0.32\textwidth]{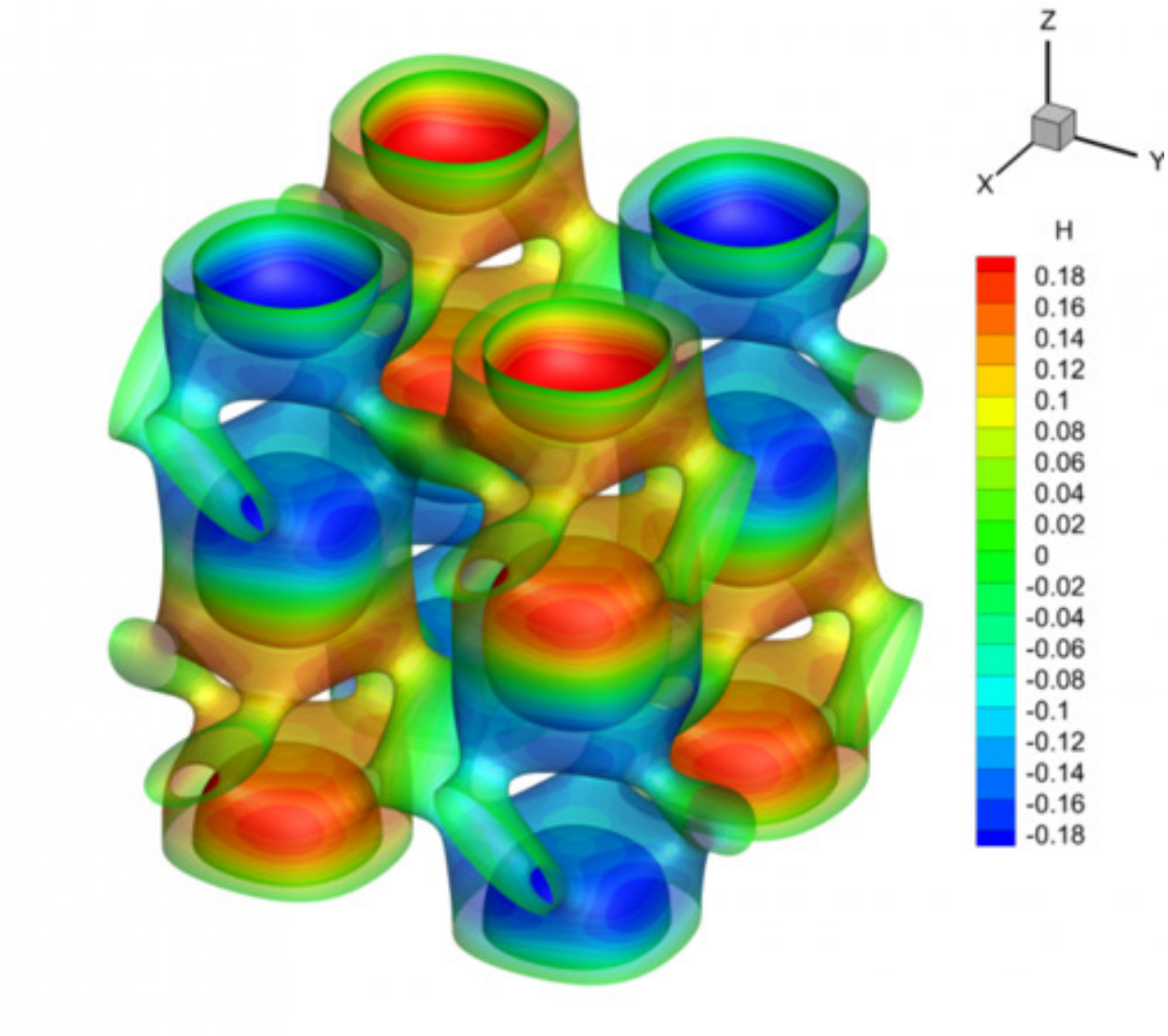}\\
			\includegraphics[width=0.32\textwidth]{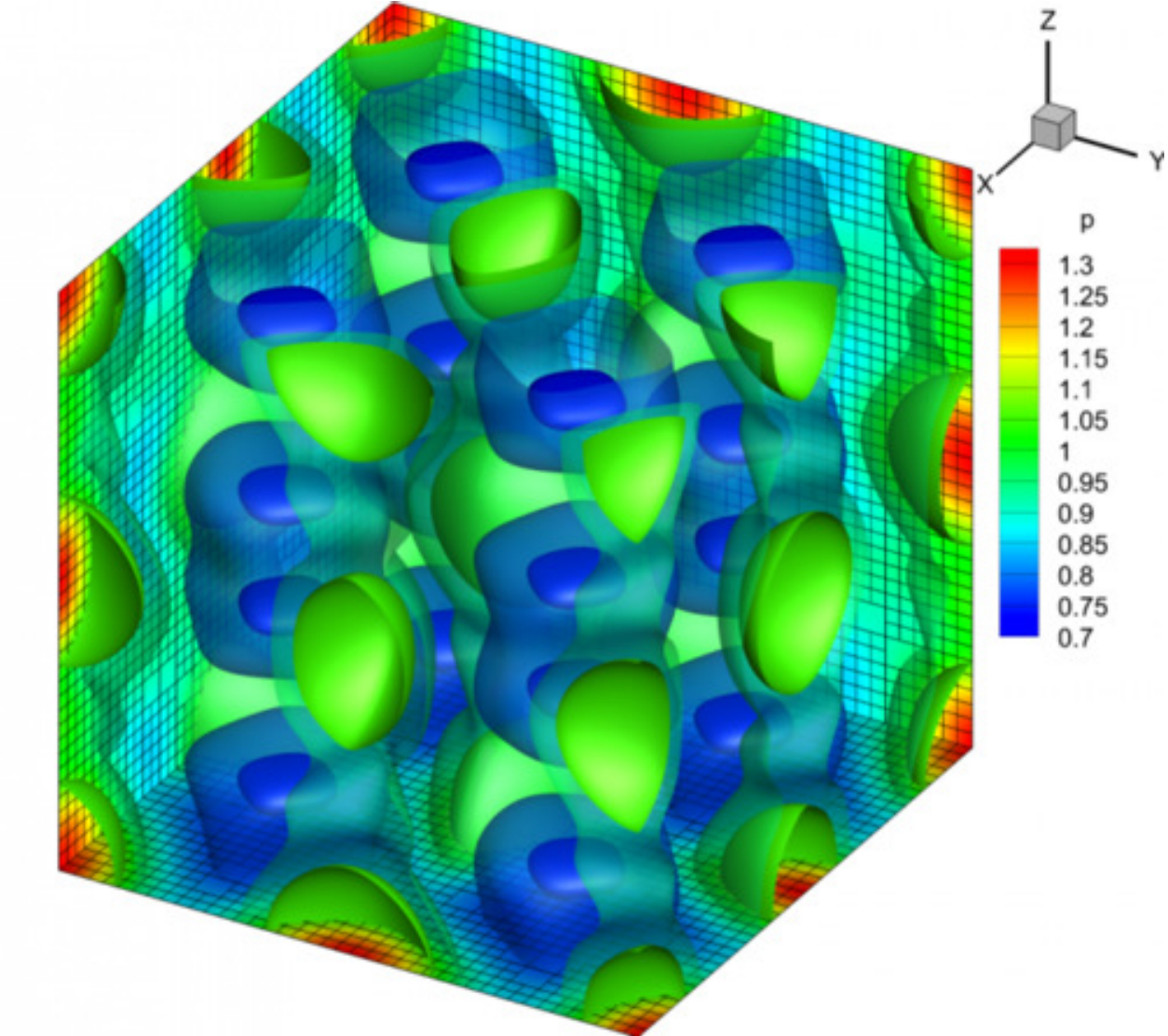}
			\includegraphics[width=0.32\textwidth]{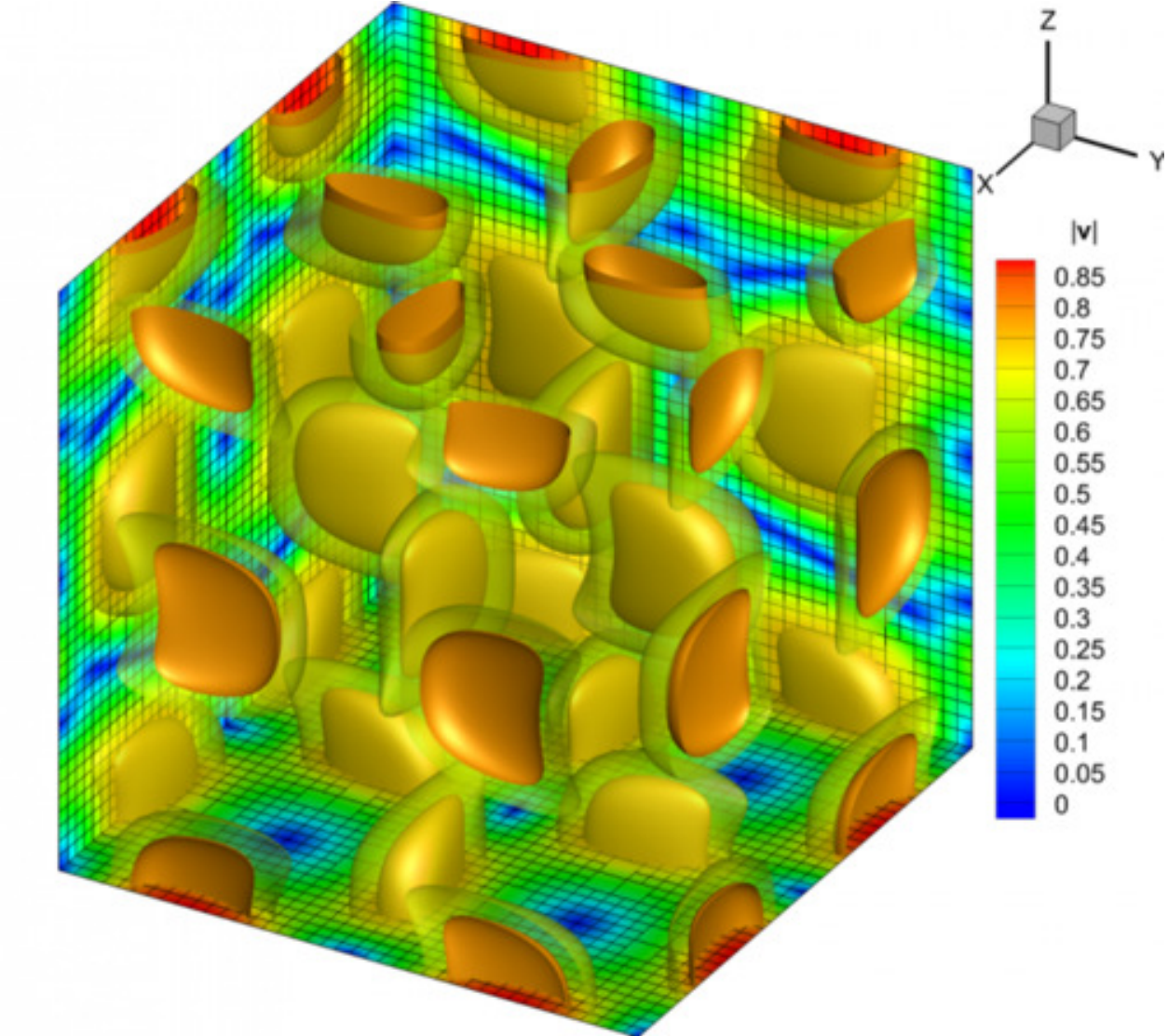}
			\includegraphics[width=0.32\textwidth]{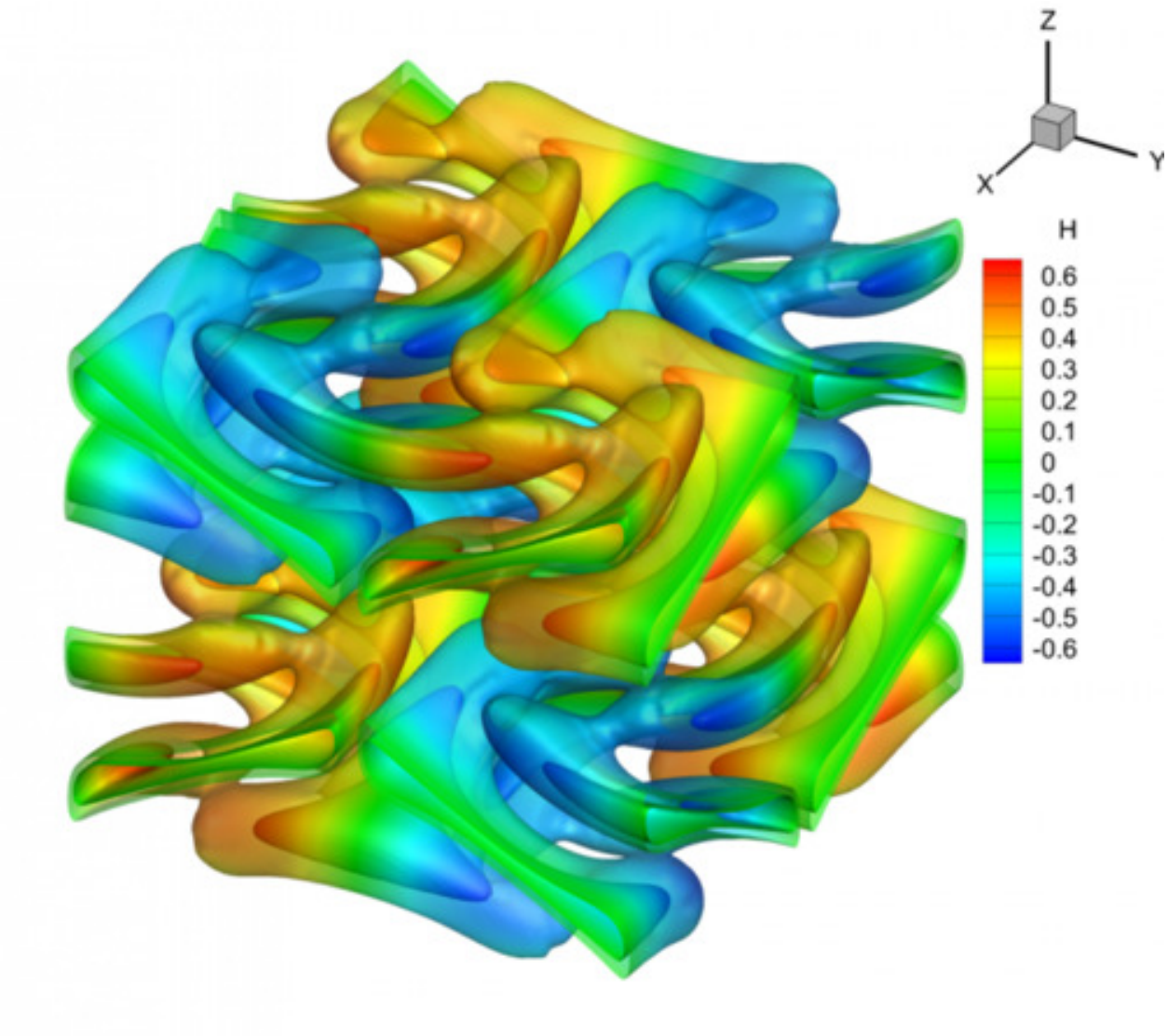}\\
			\includegraphics[width=0.32\textwidth]{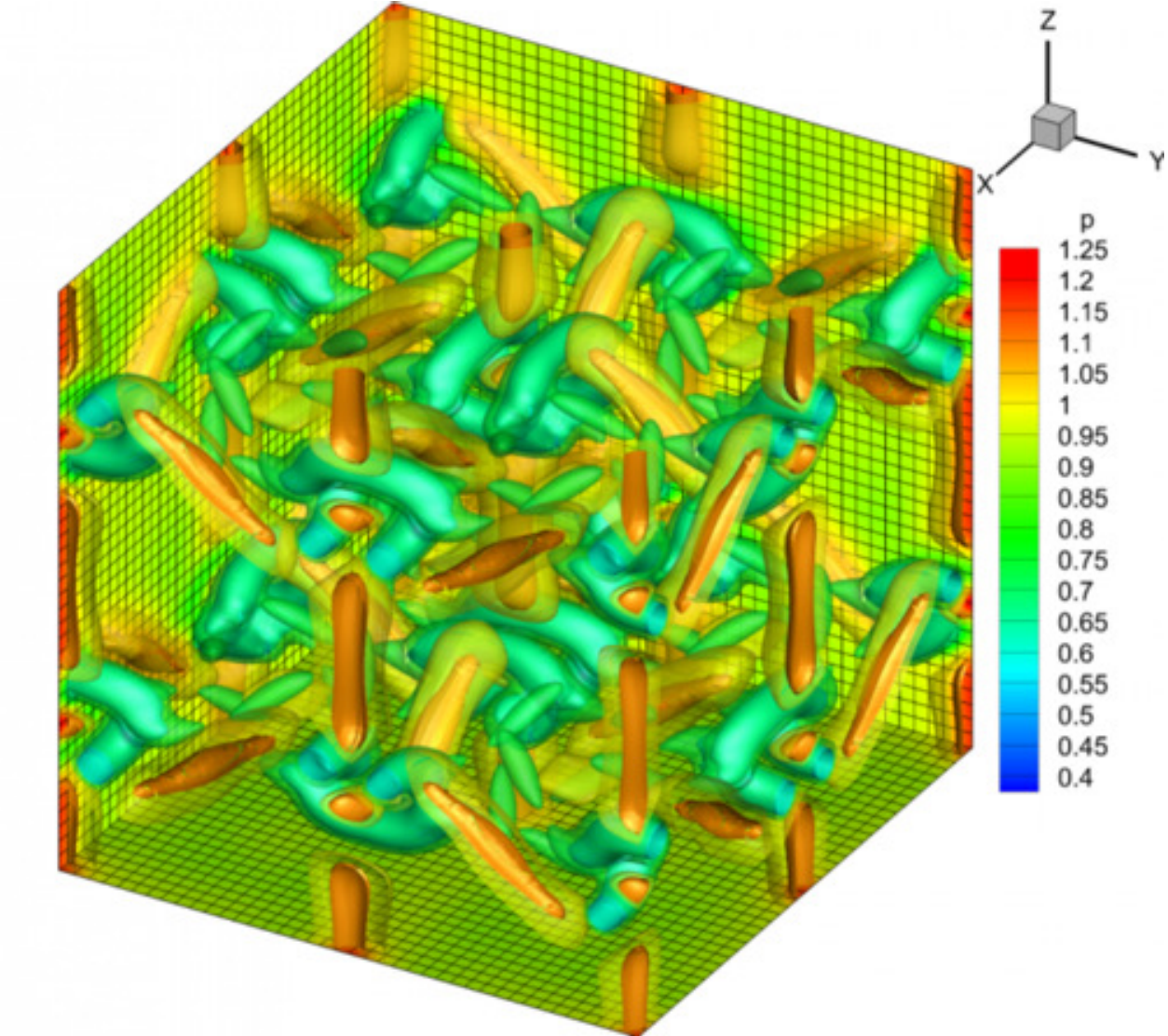}
			\includegraphics[width=0.32\textwidth]{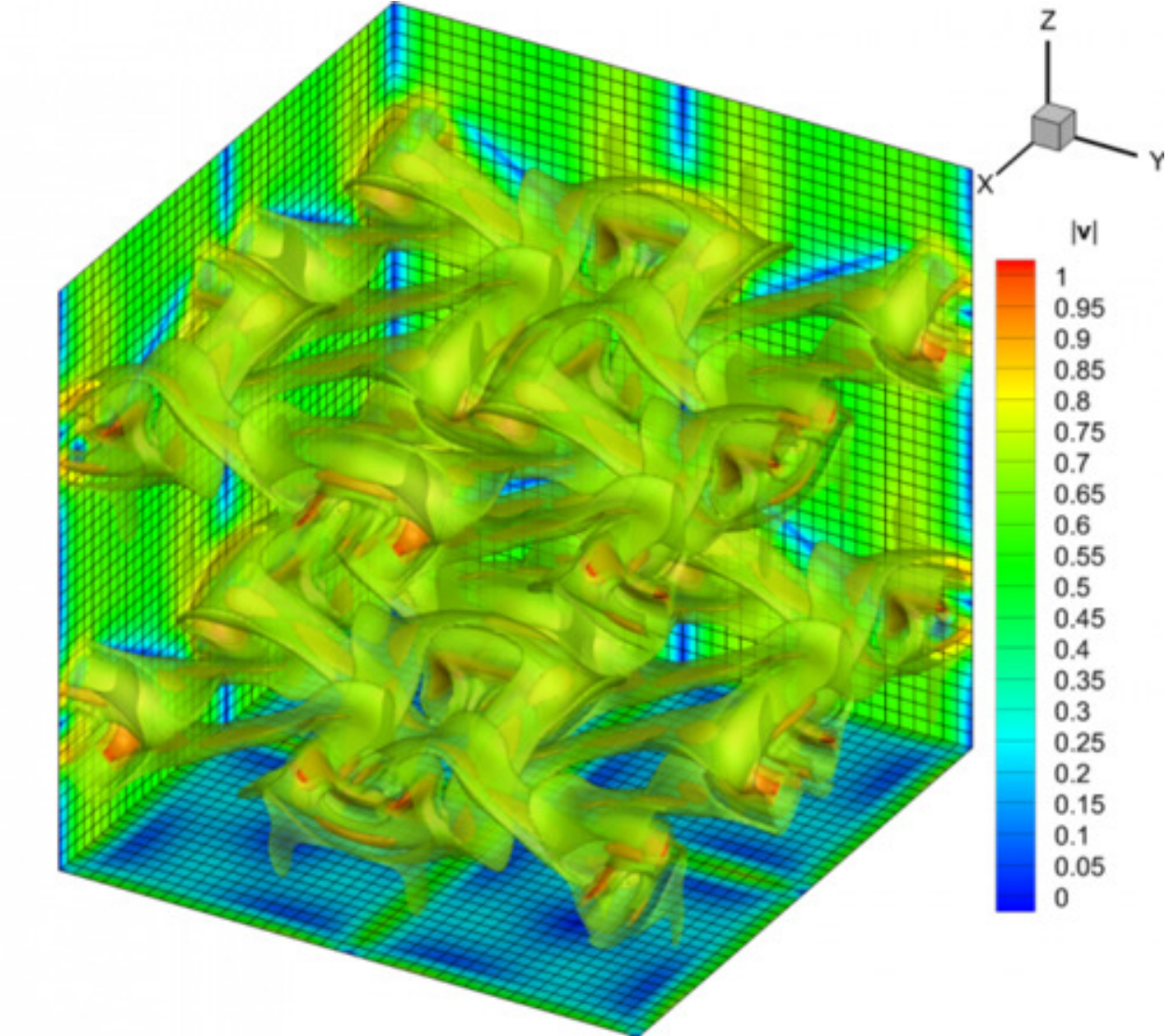}
			\includegraphics[width=0.32\textwidth]{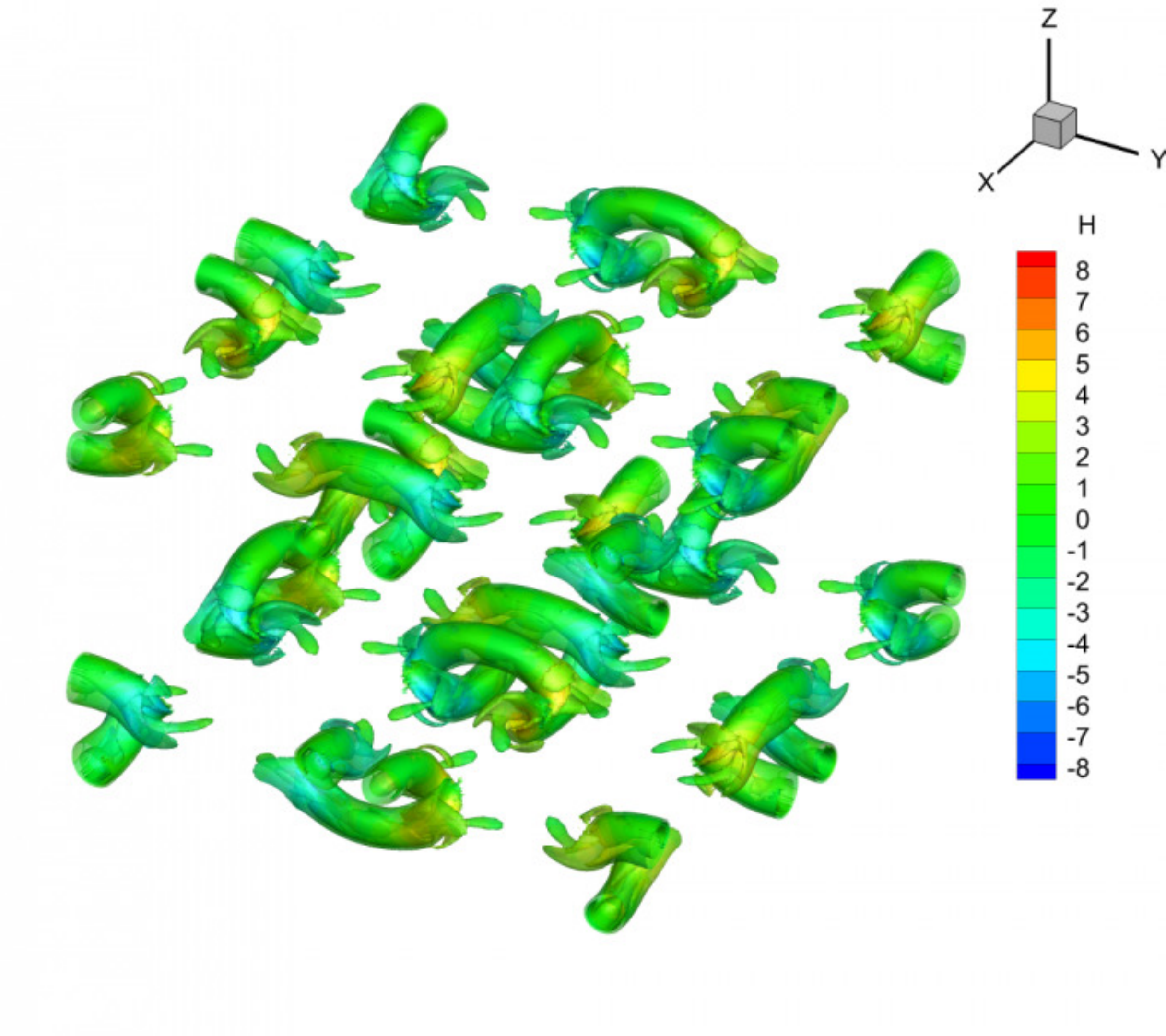}\\
			\includegraphics[width=0.32\textwidth]{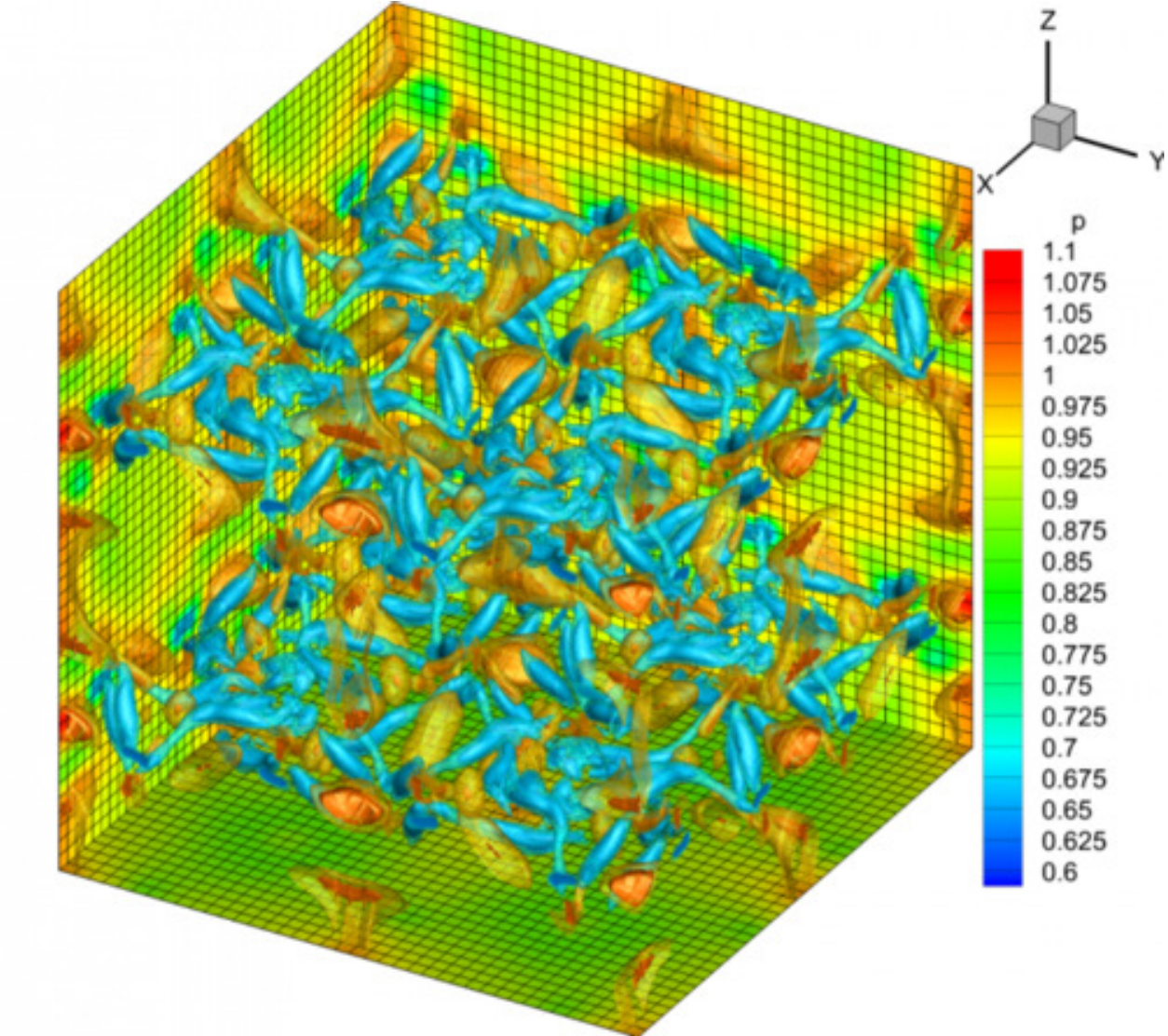}
			\includegraphics[width=0.32\textwidth]{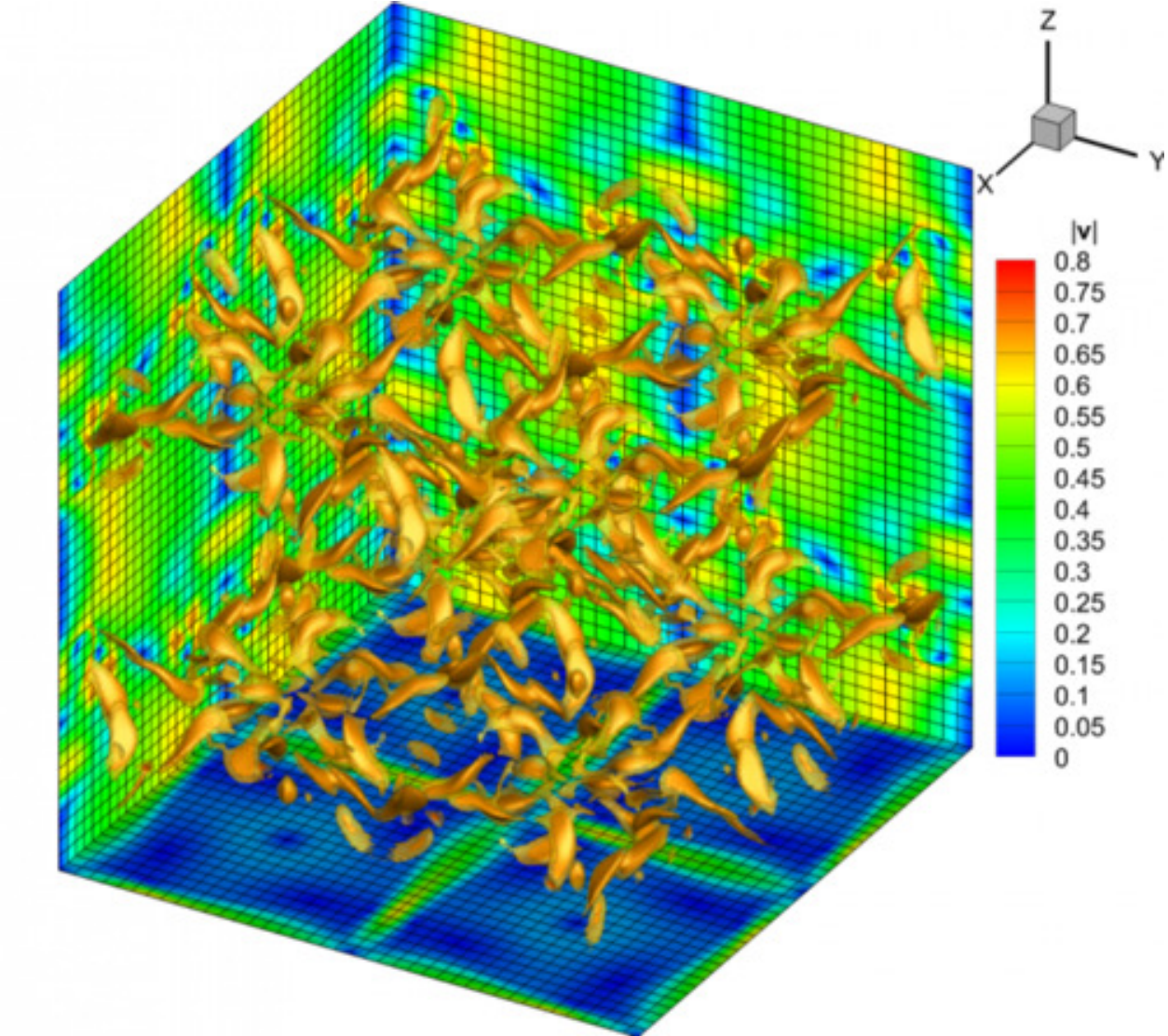}
			\includegraphics[width=0.32\textwidth]{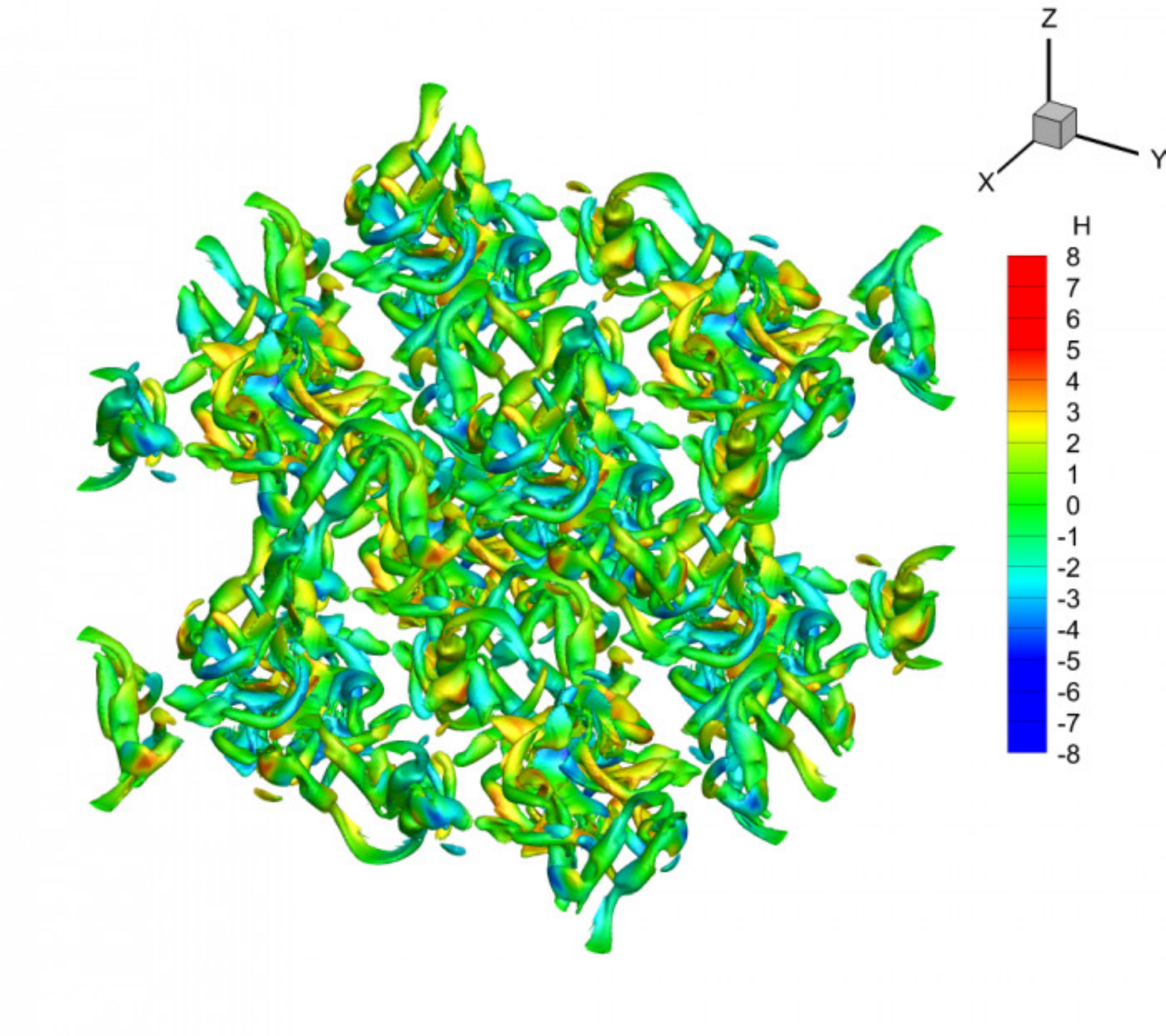}
			\caption{Numerical solution for the three dimensional Taylor-Green vortex flow at $Re=800$ computed with our $\SIDG$-$\p_{4}$ method using $20^3$ elements on the coarsest level $\Omega_h^0$, a refinement factor $\err=2$, $\ell_{\text{max}}=1$ maximum refinement levels.  
			The isosurfaces of the velocity (left), the isosurfaces of the pressure (center) and the isosurfaces of the vorticity colored by the helicity field (right) are plotted at 
			times $t=0.5$, $2.0$, $6.0$ and $10.0$ from the top to the bottom, respectively.}\label{fig:TG3D2}
\end{figure}

		The three-dimensional Taylor-Green vortex is a well-known three-dimensional benchmark problem for testing the capability of the numerical method in reproducing accurately the 
		nonlinear kinetic energy dissipation, where energy is transferred from the large to the small spatial scales. 
		A well resolved DNS reference solution is provided by Brachet et al. in \cite{Brachet1983}, obtained through a very high-resolution direct spectral method, employing up to $256^3$ Fourier
		modes and a power series analysis in time up to the $80$-th order (see \cite{Morf1980}). The physical domain is the box $\Omega = [0,2\pi]^3$ with periodic boundary conditions everywhere. 
		The initial condition consists in 
\begin{align}
&u(x,y,z,0) = \sin(x)\cos(y)\cos(z),\\
&v(x,y,z,0) = -\cos(x)\sin(y)\cos(z),\\
&w(x,y,z,0) = 0,\\
&p(x,y,z,0) = \frac{1}{16}\left(\cos(2x)+ \cos(2y) \right)\left(\cos(2z)+2\right),
\end{align} 
		which provides the initially very smooth vortices in a non-equilibrium state. Then, depending on the kinematic viscosity $\nu$, the fluid flow goes through a highly non-linear decaying process, in which the kinetic energy $E=\mathbf{v}^2/2$ is transported from the lower (large scale processes) to higher modes (smallest scales) until reaching the dissipative viscous regime.
		A quantitative comparison with the reference data of \cite{Brachet1983} is presented in Figure (\ref{fig:TG3D2}) for the time evolution of the kinetic energy \emph{dissipation rate} 
\begin{align}
\epsilon(t) = - \frac{\partial K}{\partial t}= - \frac{1}{\left| \Omega \right|} \frac{\partial}{\partial t}  \int_{\Omega} \frac{1}{2} \mathbf{v}^2  d\x.
\end{align}  
When the dissipation rate $\epsilon(t)$ reaches its peak at a given time $t=t_p$, this means the higher modes reach the highest population, entering into the viscous regime and being soon destroyed by dissipative forces. At this crucial point, the higher resolution provided by the AMR framework is needed for allowing the smallest scales to dissipate properly. Whenever an under-resolved solution is  computed, then the kinetic energy dissipation rate behaves improperly: if the numerical method is non-dissipative, then at higher Reynolds numbers $\epsilon$ is expected to be \emph{under-estimated} leading to possible spurious oscillations, i.e. the higher modes saturate the computational domain for the wave numbers $\mathpzc{F}[\Omega_h]$, far  both from the physical and from the numerical viscous regimes, $\mathpzc{F}$ denoting the Fourier operator; if the numerical method is over-dissipative, then $\epsilon$ is expected to be \emph{over-estimated}, meaning the numerical viscous regime appears sooner in $\mathpzc{F}[\Omega_h]$ with respect to the physical one. From Figure \ref{fig:TG3D} we can conclude that the results obtained with our adaptive semi-implicit staggered DG scheme are in  good agreement with the DNS reference solution.

		Figure \ref{fig:TG3D2} shows the time evolution of pressure, velocity and vorticity magnitude obtained with our $\SIDG$-$\p_4$ at Reynolds $Re=800$, within a main mesh $\Omega_h$ made of only $20^3$ space-elements at the coarsest level $\Omega_h^0$, with a refinement factor $\err=2$ and $\ell_{\text{max}}=1$ refinement levels. 

		\subsection{Vortex ring dynamics} 
		
		For describing the physical state of an incompressible fluid, vorticity and velocity are essentially interchangeable quantities. On the other hand, higher Reynolds numbers and turbulent 
		flow regimes are typically characterized by localized higher vorticity zones. This fact leads to think that a deeper comprehension in the vortex dynamics could allow a deeper comprehension 
		of high Reynolds fluid dynamics. 
		The visualization and measurement of vorticity constitutes a non-trivial problem in experimental studies. Then, over the years, several experiments, theoretical analysis and numerical 
		investigations provided typical test problems that nowadays can be used for testing the ability of a numerical method 
  	to give an accurate description of vortex dynamics. 
	In the following tests the initial conditions are given in terms of vorticity, then the \emph{real} initial condition for the numerical simulations are computed after recycling exactly the same discrete operators depicted in the previous theoretical sections.
	
	Indeed, for incompressible fluids a vector potential $\mathbf{A}$ can be introduced by defining the velocity vector as $\mathbf{v} = \nabla \times \mathbf{A} $. Then the vorticity can be written in terms of  $\mathbf{A} $ as
	\begin{align}
	\omega = \nabla \times \mathbf{v} =  \nabla \times ( \nabla \times \mathbf{A} ) \equiv \nabla  (\nabla \cdot \mathbf{A}) - \nabla ^2 \mathbf{A}.
	\end{align}
	Notice that any divergence-free component of vector $\mathbf{A} $ would not contribute to the corresponding velocity field. 
	Then, by looking for solutions satisfying the divergence-free condition for $\mathbf{A}$, a very simple Poisson equation follows for $\mathbf{A}$, yielding a very familiar system of equations for the potential vector components $A_k$, $k=x$, $y$, $z$, i.e. 
	\begin{align}
	  \mathbf{B} &= \nabla A^{\kapped}\\
	\nabla \cdot \mathbf{B} &= - \omega^{\kapped}
	\end{align}
where vector	$\mathbf{B}$ is an auxiliary variable. Indeed, this system shows to be very similar to our governing equations (\ref{eq:NSmom}-\ref{eq:NSinc}) and a consistent discrete DG-$\p_N$ formulation on the presented staggered grids, $A^{\kapped}_h,\omega^{\kapped}_h\in\p_N(\Omega_h)$ and $B^{\kapped}_h\in \p_N(\Omega_h^{\dualk{k}})$, reads 
	\begin{align}
	\mathbb{H} \tilde{A}^{\kapped}&= \mathbb{M} \tilde{\omega}^{\kapped}
	\end{align}
		after setting $\Delta t=1$, $\tilde{\omega}^{\kapped} $ and $\tilde{A}^{\kapped}$ are respectively the vector of the degrees of freedom for the $k$-th vorticity and potential vector components. After the three equivalent discrete Poisson systems are solved for the potential vector, then the velocity is updated accordingly to the $L_2$ projection over $\p_N$ of the definition $\mathbf{v}=\nabla \times \mathbf{A}$.

		\subsubsection{Vortex ring pair collision} 
	
\begin{figure} 
\centering 
			\includegraphics[width=0.35\textwidth]{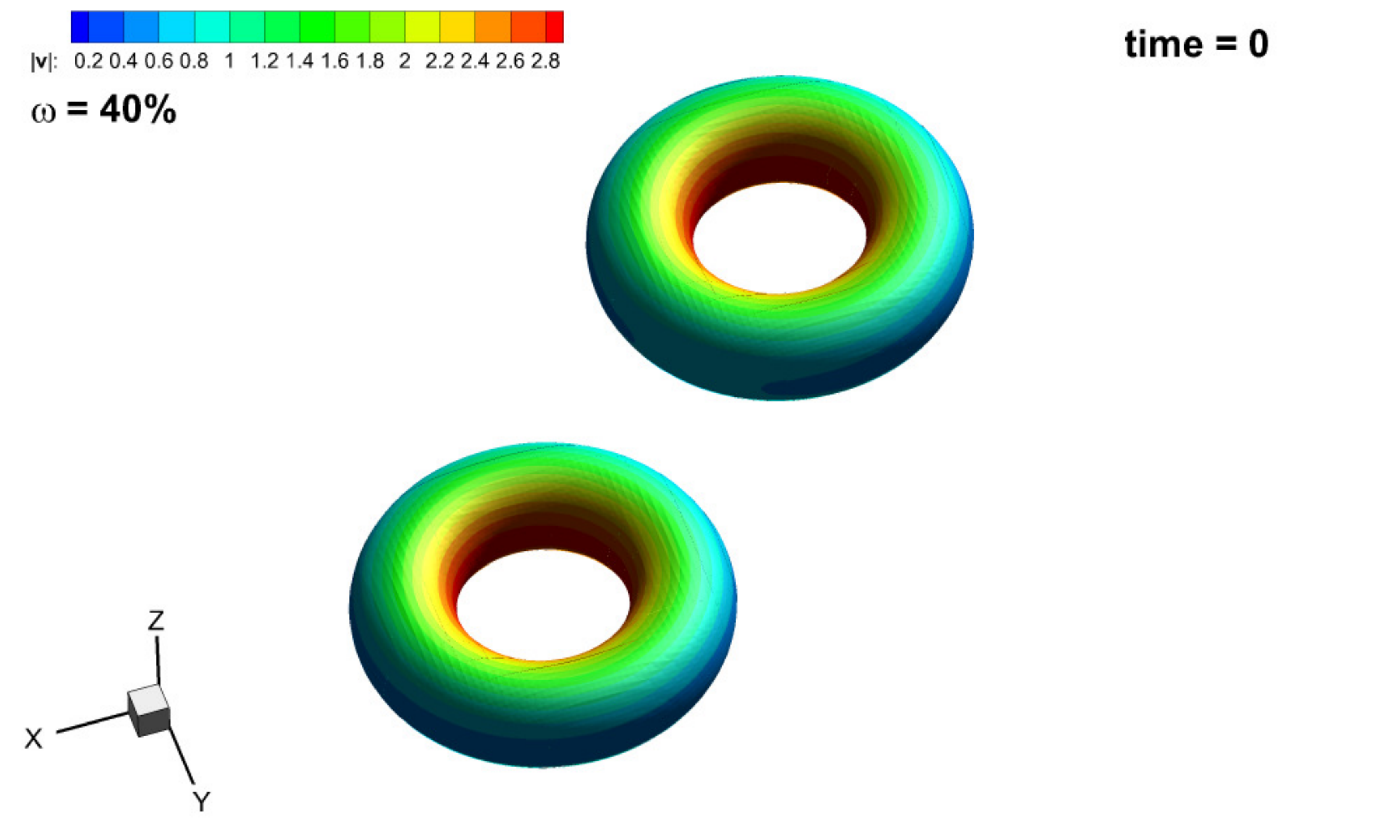}\;\includegraphics[width=0.35\textwidth]{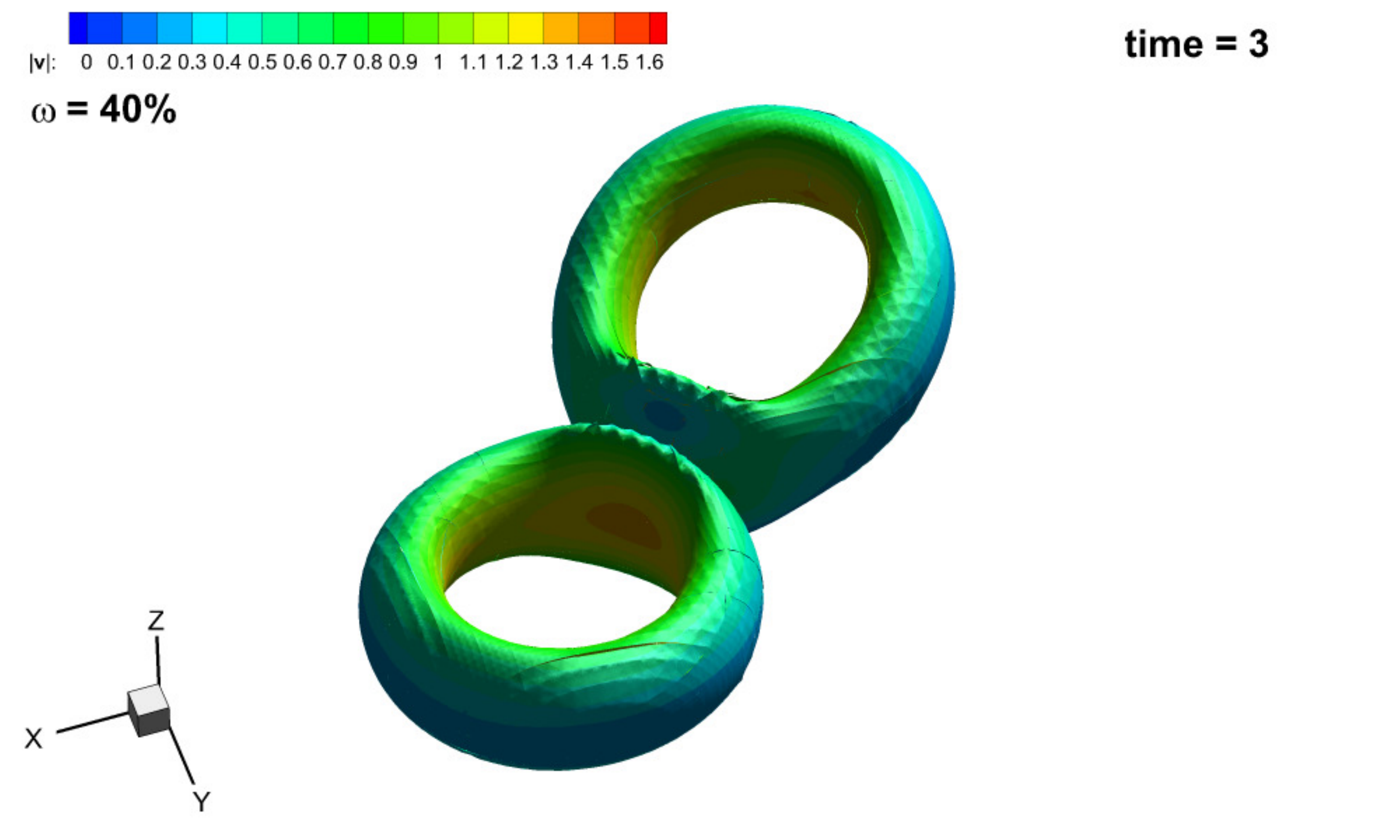}\\
			\includegraphics[width=0.35\textwidth]{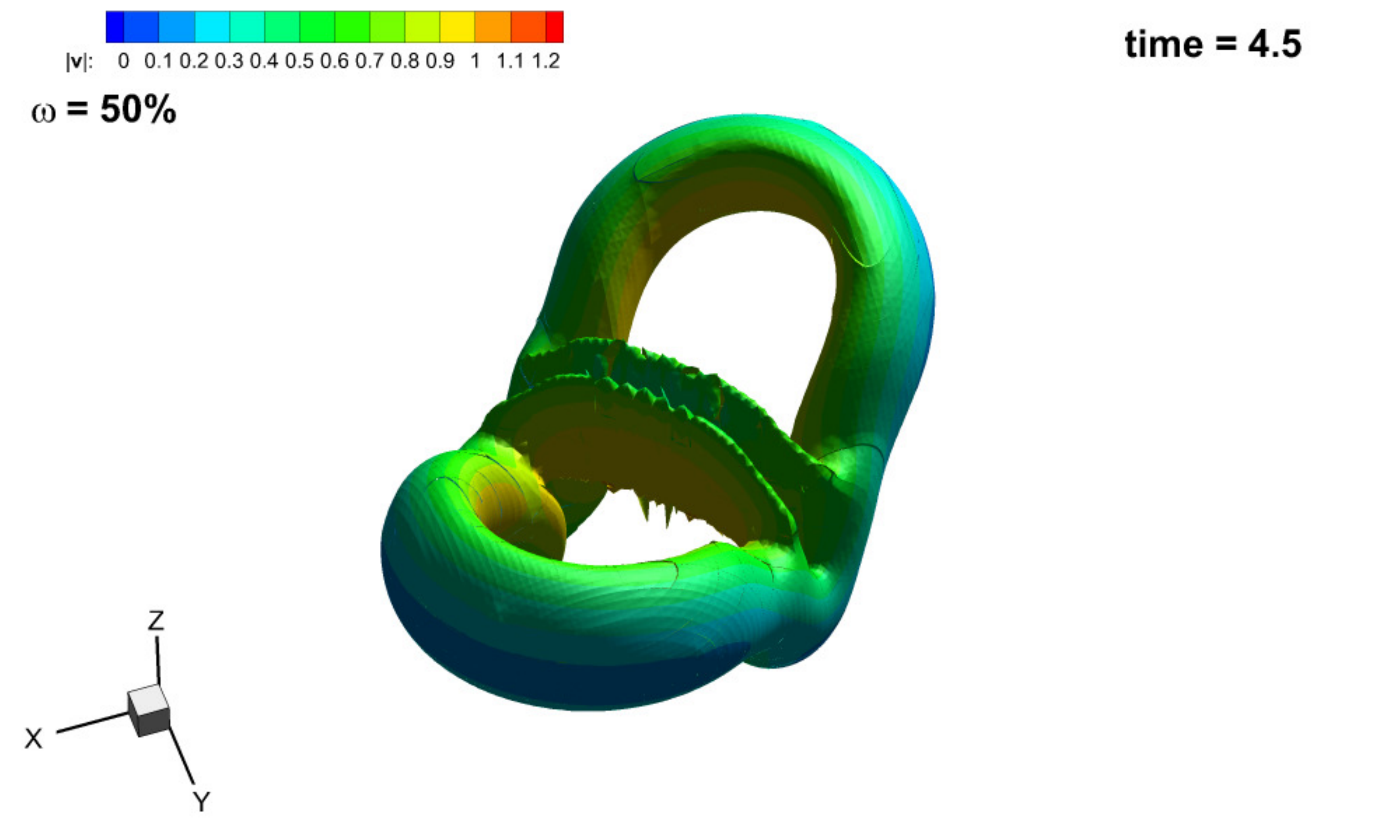}\;\includegraphics[width=0.35\textwidth]{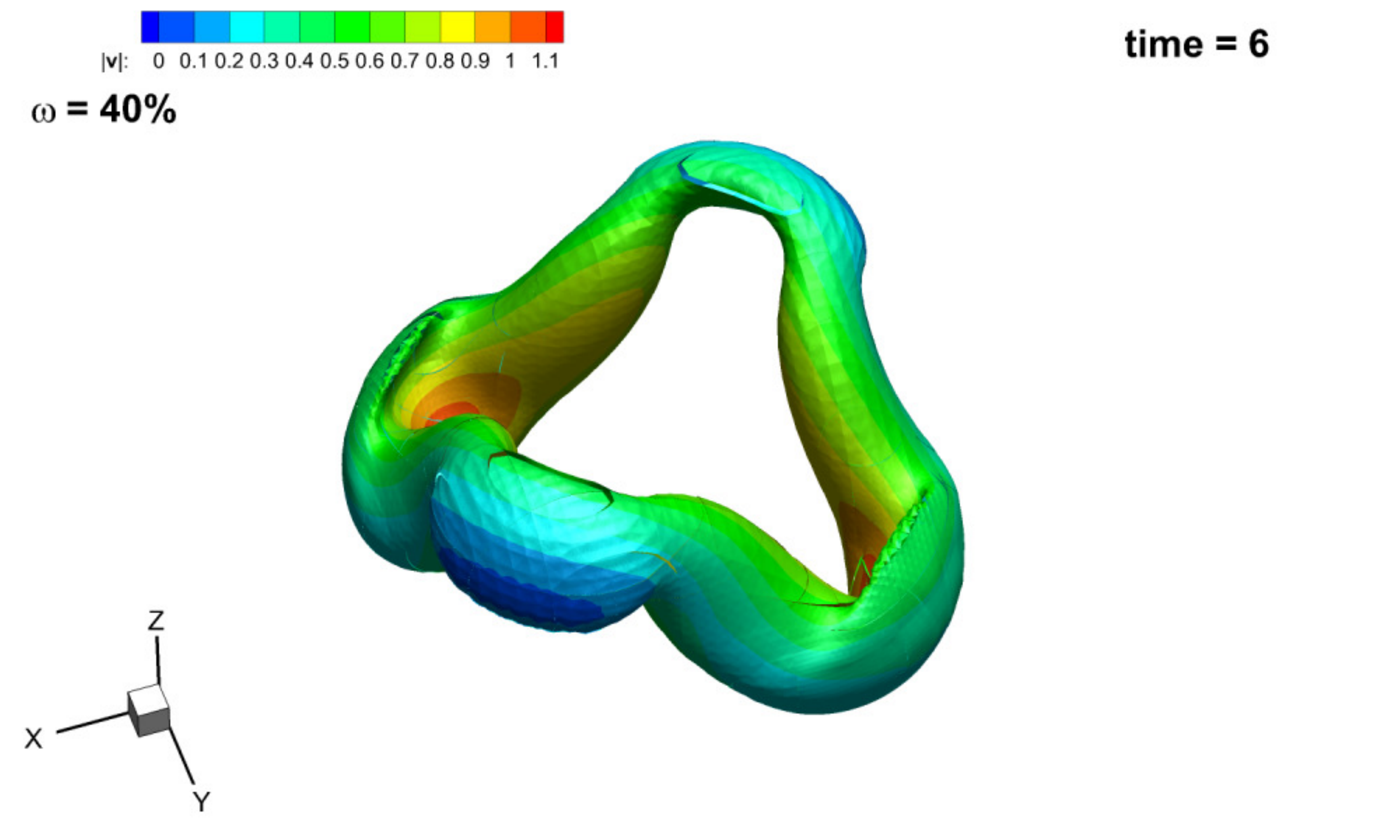}\\
			\includegraphics[width=0.35\textwidth]{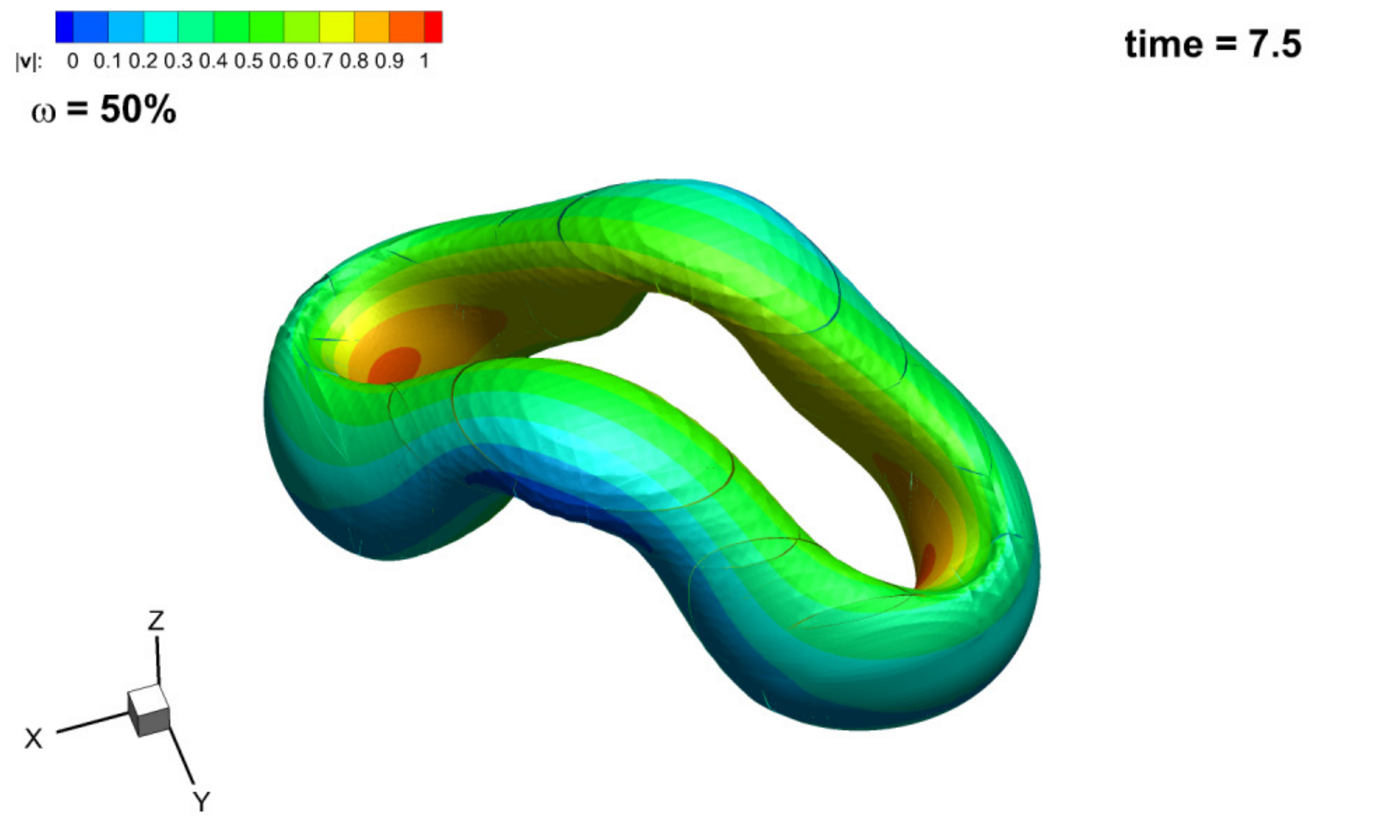}\;\includegraphics[width=0.35\textwidth]{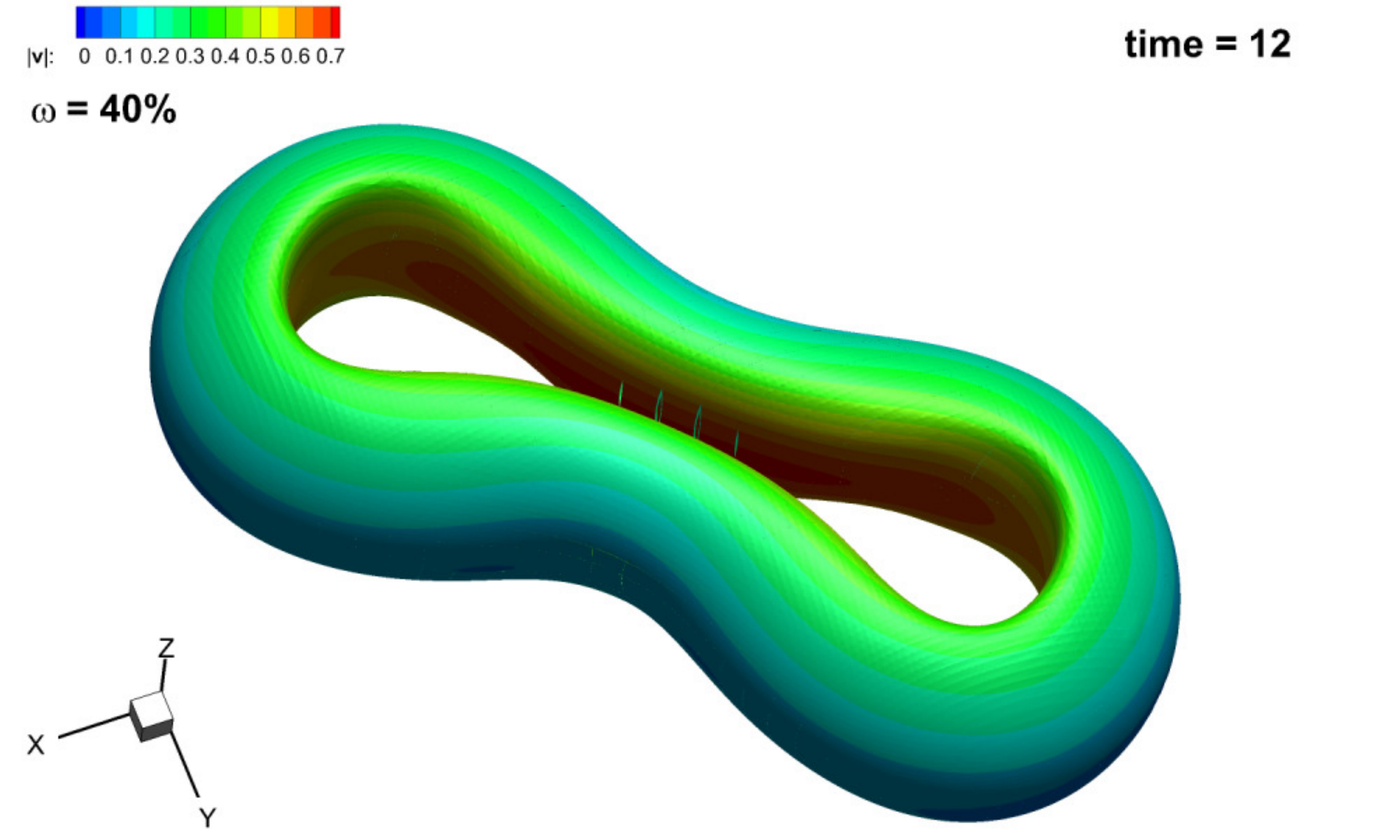}\\
			\includegraphics[width=0.35\textwidth]{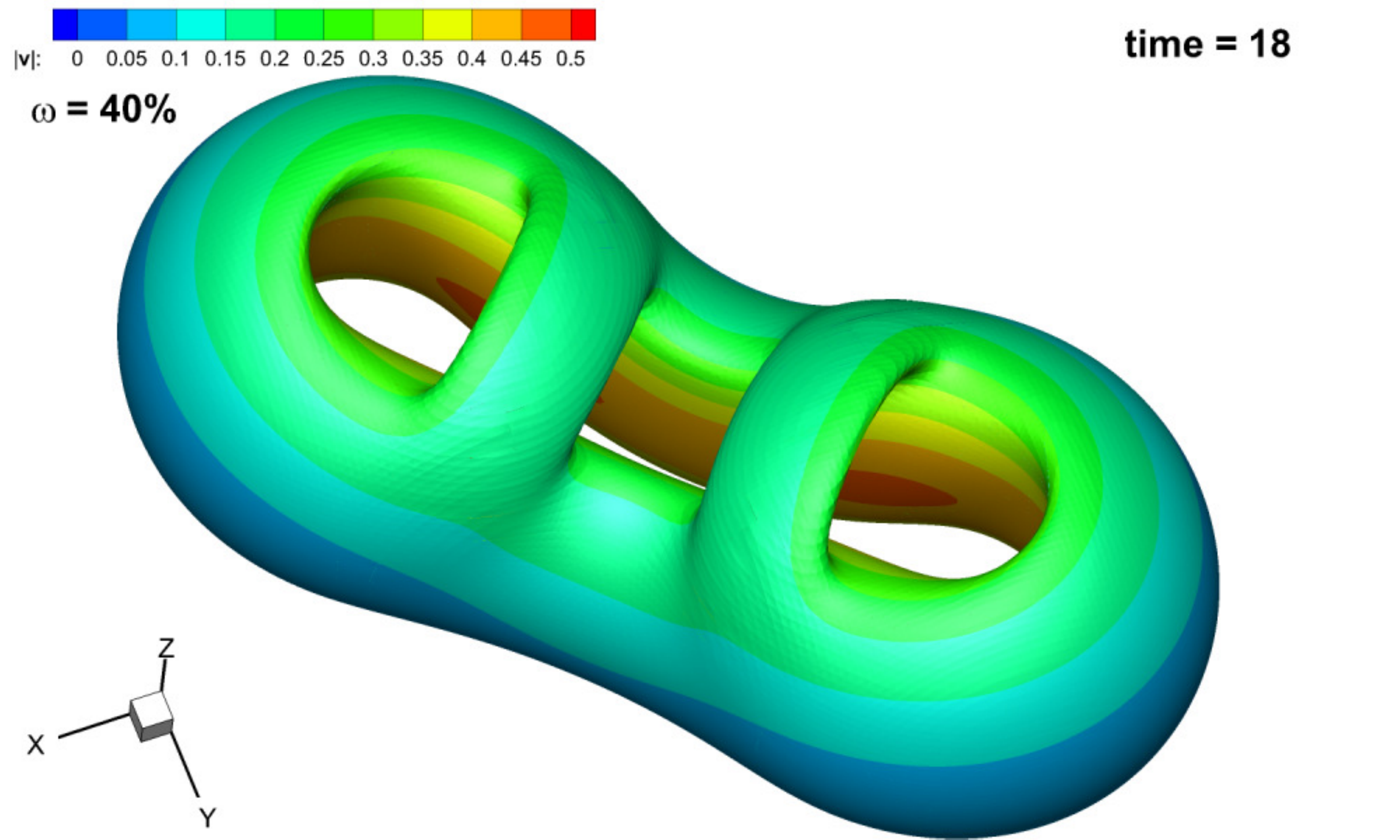} 
\caption{Time evolution  of the iso-surfaces for the vorticity magnitude $|\omega |$ in the three-dimensional vortex-ring pair interaction problem at different times, respectively, from left to right, from the top to the bottom:  t=$0.0$, $3.0$, $4.5$, $6.0$, $7.5$, $12.0$ and $18.0$, using $30^3$ elements on the coarsest grid with periodic boundary conditions; these results are obtained with the $\p_{4}$-version of our \SIDG method. The iso-vorticity values has been chosen accordingly to the percentage referred to the maximum $|\omega\ov{\text{max}}$ accordingly to \cite{Kida1991}.}\label{fig:VPI}
\end{figure}

	In this three-dimensional test, two coplanar vortex rings are initialized within a cubic spatial domain $\Omega=[-\pi,\pi]^3$ with periodic boundary condition everywhere, centered in $C_{1,2}=\pm(D\cos(\pi/4)/2, D  \cos(\pi/4)/2)$ with $D=1.83$, a major radius $R=0.491$, centered along the $x-y$ plane, according to a Gaussian distribution  in the ring core (matches case II of \cite{Kida1991})
	\begin{align}
	|\omega|(r) = \omega_0 \exp\left[-\left(\frac{r}{a}\right)^2\right], \label{eq:vort}
	\end{align} 
where  $r$ is the radial direction, $a=0.196$ is the effective thickness of the core, $\omega_0=23.8$ the vorticity amplitude.  For the present test the Reynolds number is $Re_{\gamma} = \gamma/\nu = 577$, $\gamma= \pi \omega_0 a^2$ being the circulation of the vortex ring, yielding a kinematic viscosity $\nu \sim 5\times 10^{-3}$.
The two ring vortices proceed in the vertical direction because of \emph{self-induction} and, at the same time, the two rings are attracted towards the $x=y$ plane by \emph{mutual-induction}. At around $t=3$ the two vortex rings collide. Notice that vortex lines are anti-parallel  at the contact point, and consequently dissipated by viscous interactions. Then, the resulting dynamics becomes highly non linear, with very complex effects due to both self- and mutual- induced interaction (see a review for vortex reconnection \cite{KidaTakaoka1994} and \cite{Kida1991} for an almost complete overview of the vortex collision process). Figure \ref{fig:VPI} shows the dynamics of the vortex interaction by means of the selected iso-surfaces for the vorticity magnitude $|\omega|$, showing the main evolution phases before the ring collision, the collision phase and the post-collision phase, with a good description of the \emph{bridging mechanism}. The computed results show to be in good agreement with the provided reference solution in literature, e.g. see \cite{Kida1991,Gosh1994,HahnIaccarino}. For this test, the computational domain has been discretized using $30^3$
space-elements on the coarsest grid $\Omega_h^0$ within the space of solutions of our $\SIDG$-$\p_4$ method, a refinement factor $\err=2$ and a maximum number of refinement levels $\ell_{\text{max}}=1$.
The mesh is shown to be automatically refined only next to the vorticity cores, saving much computational effort compared to a uniform fine grid. 
Figure \ref{fig:VPI} shows the time evolution of the vorticity field interpolated along the two orthogonal vertical planes $x=y$ and $x=-y$. At $t=4.5$, some spurious oscillations arise in the vorticity field along the collision plane $x=y$, reflecting the fact that, in the vorticity formulation, a potentially dangerous steep gradient is generated. Indeed, it is a well known fact that vortex reconnection is allowed only for viscous fluids. This means that, in the high Reynolds regime, the vortex rings would approach to be conserved in time, leading, probably, to a highly non-linear inviscid interaction in the collision phase. In this case, very high resolution methods are needed and a limiting strategy would become necessary for a pure DG method for resolving the physics of the fluid flow.

\begin{figure} 
\centering 
			\includegraphics[width=0.5\textwidth]{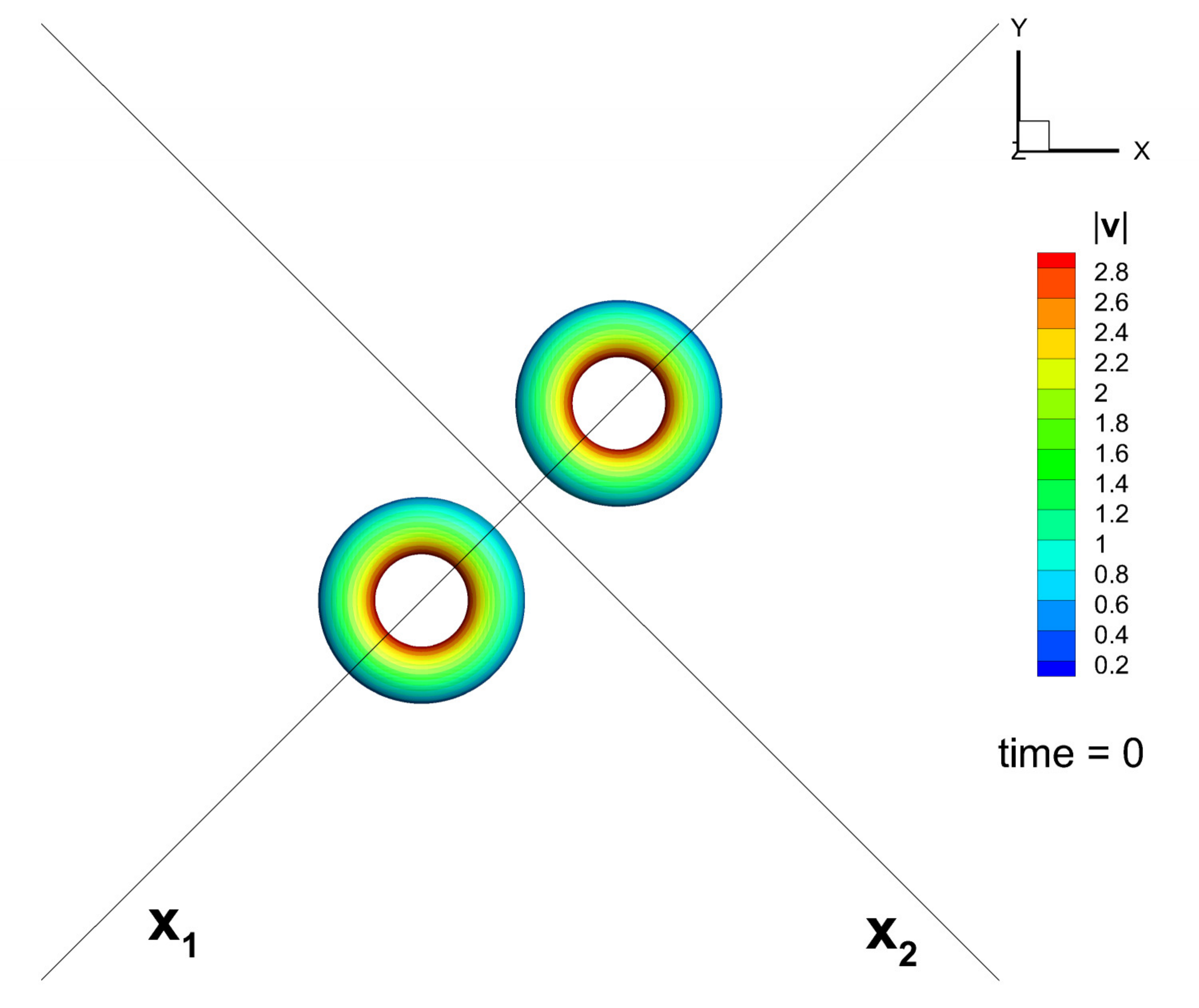}\\
			\includegraphics[width=0.48\textwidth]{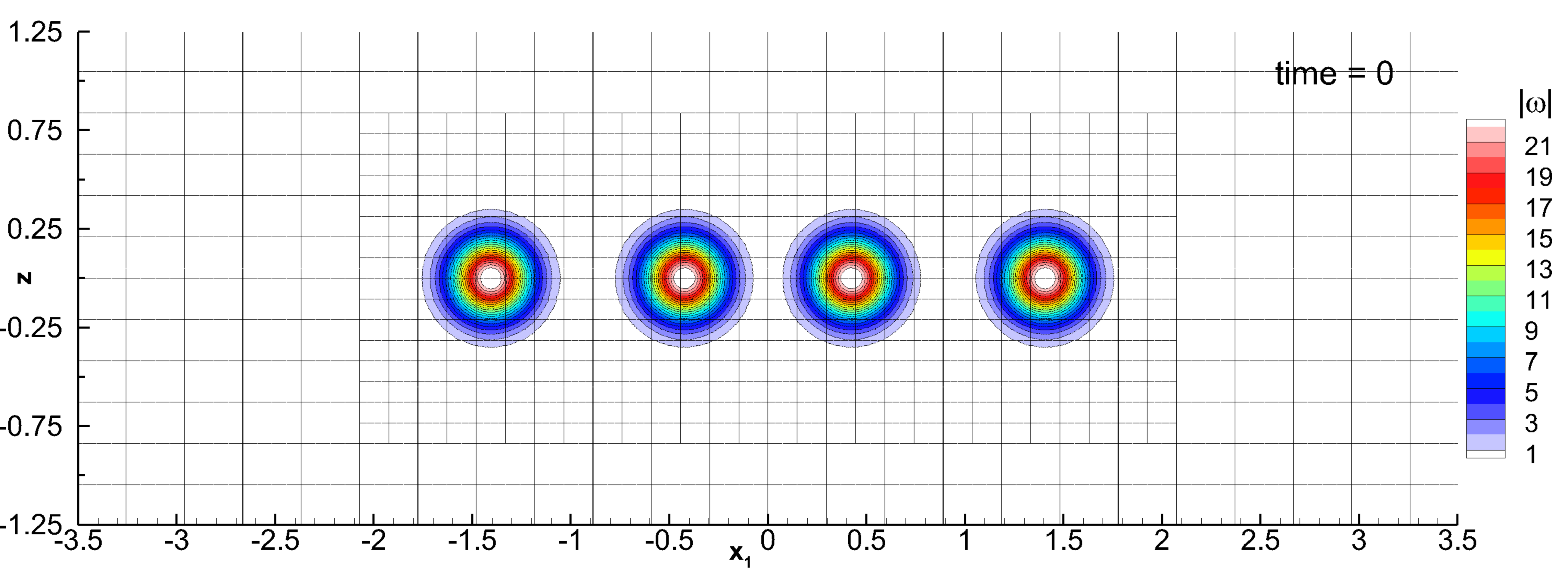}\;\includegraphics[width=0.48\textwidth]{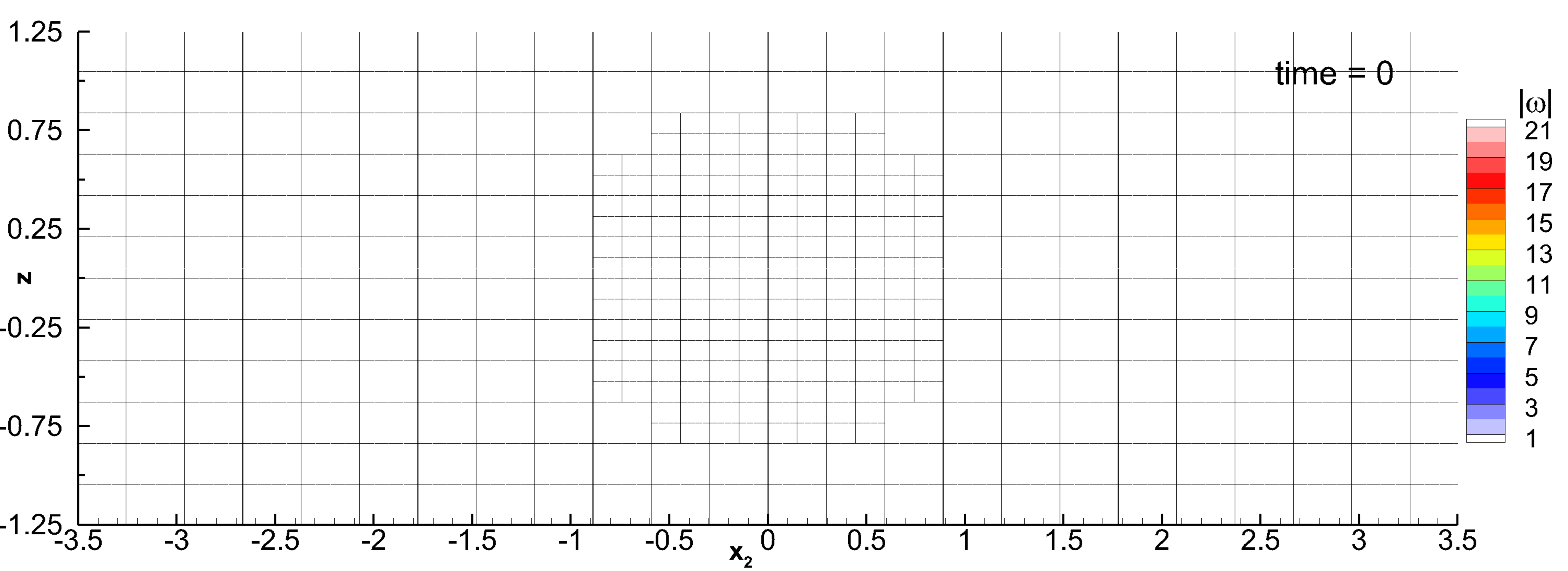}\\
			\includegraphics[width=0.48\textwidth]{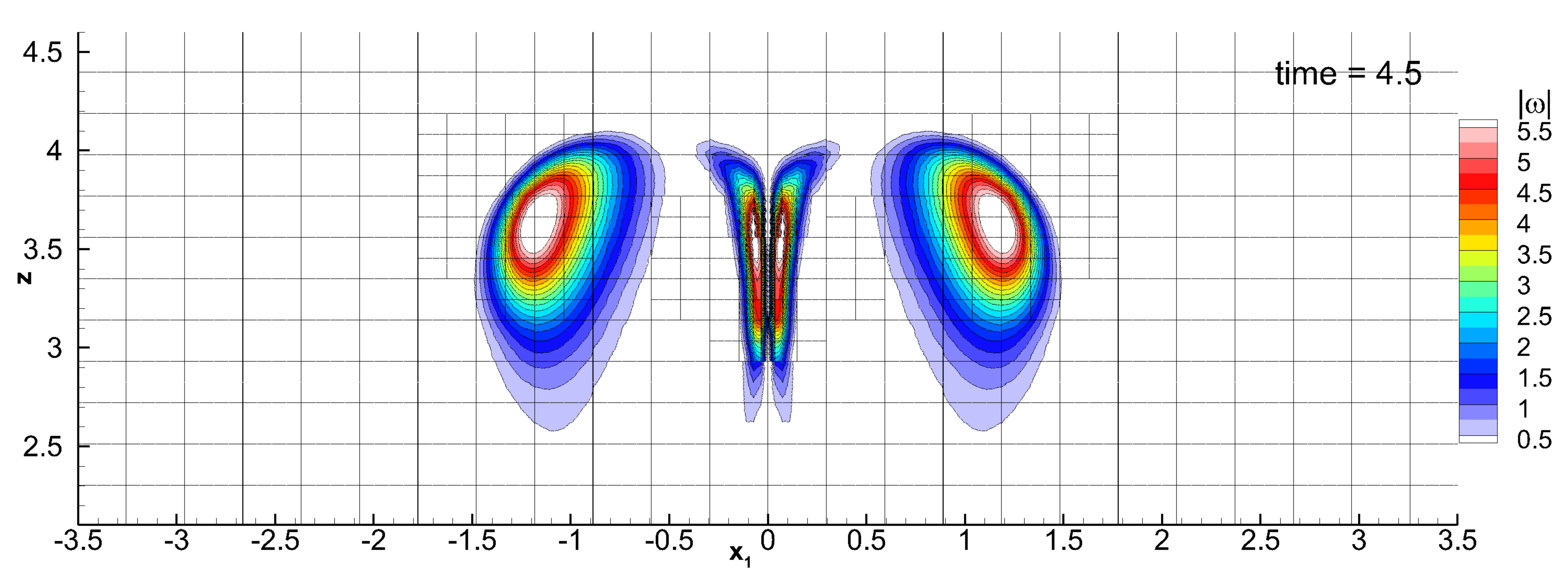}\;\includegraphics[width=0.48\textwidth]{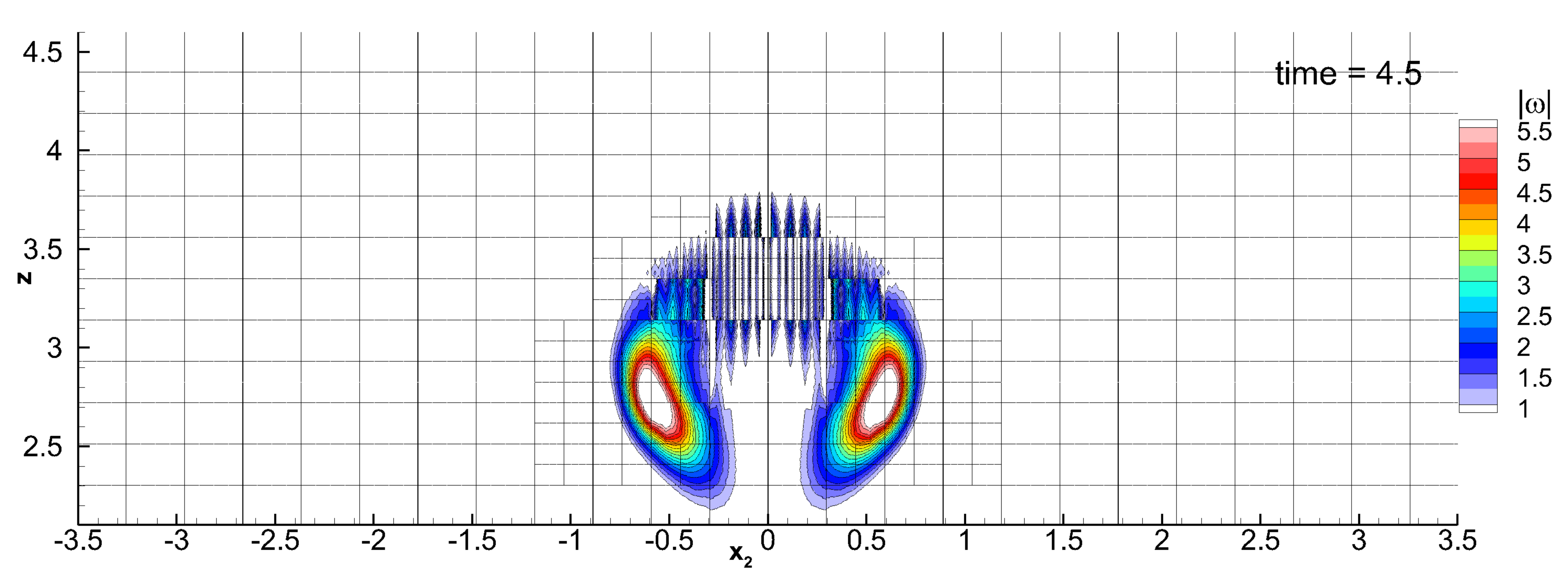}\\
			\includegraphics[width=0.48\textwidth]{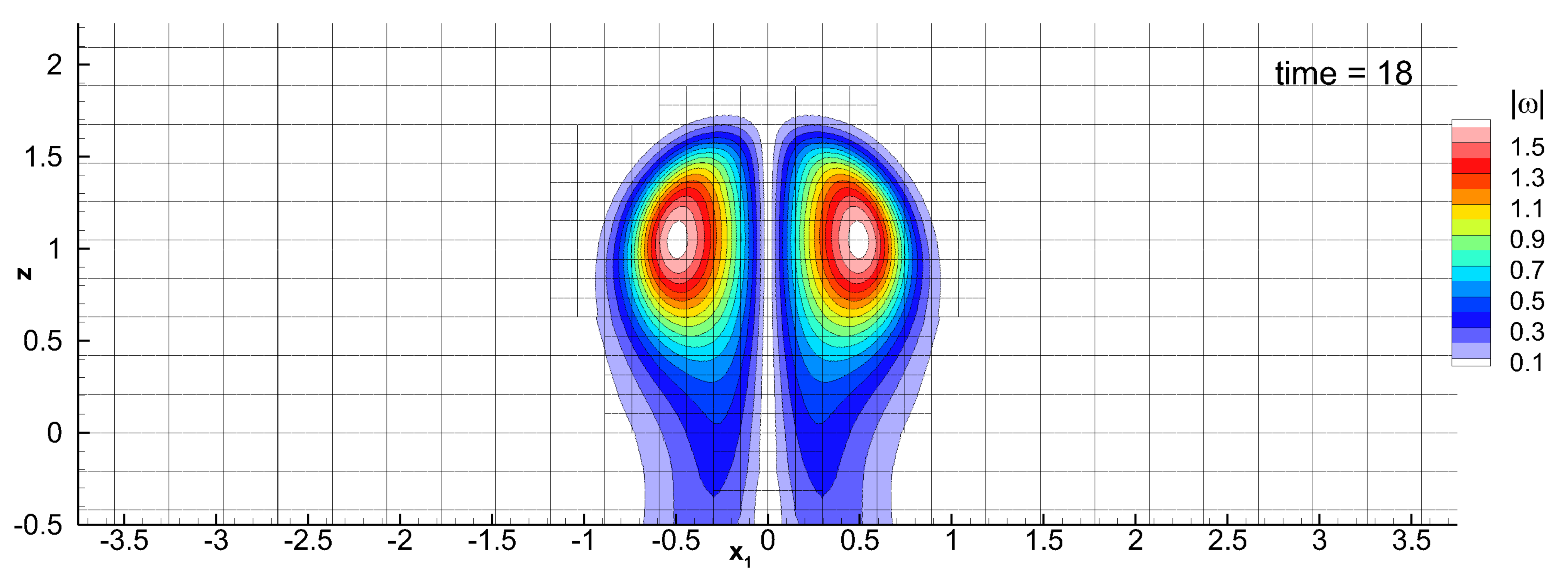}\;\includegraphics[width=0.48\textwidth]{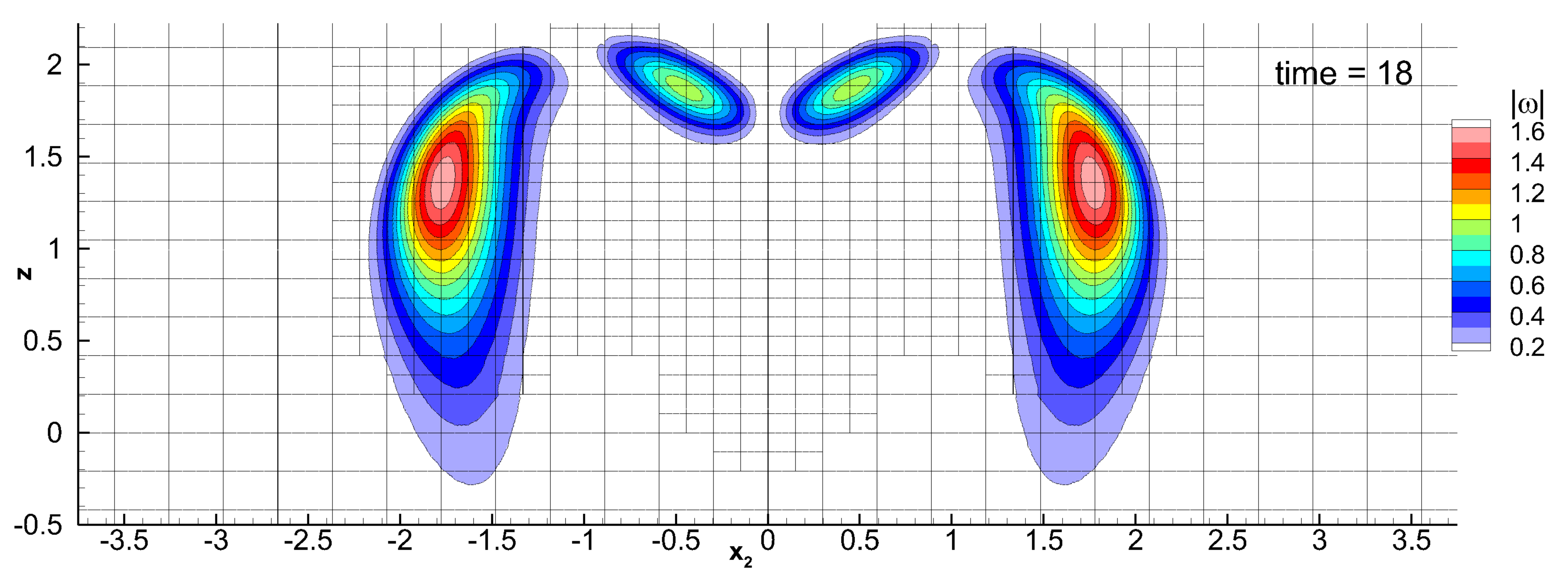}
\caption{Time evolution  of the iso-surfaces for the vorticity magnitude $|\omega|$ in the three-dimensional vortex-ring pair interaction problem at different times; the numerical solution has been interpolated along the two-dimensional orthogonal planes $x_1-z$  (left) and $x_2-z$ (right), passing through the two bisectors $x_1$  and $x_2$ of  $x-y$, respectively (see image at the top);   these results are obtained with the $\p_{4}$-version of our \SIDG  method using $30^3$ elements on the coarsest grid with periodic boundary conditions at times t=$0.0$, $4.5$ and $18.0$.}\label{fig:VPI2}
\end{figure}


		\subsubsection{Vortex ring pair leapfrogging} 
\begin{figure} 
\centering 
			\includegraphics[width=0.45\textwidth]{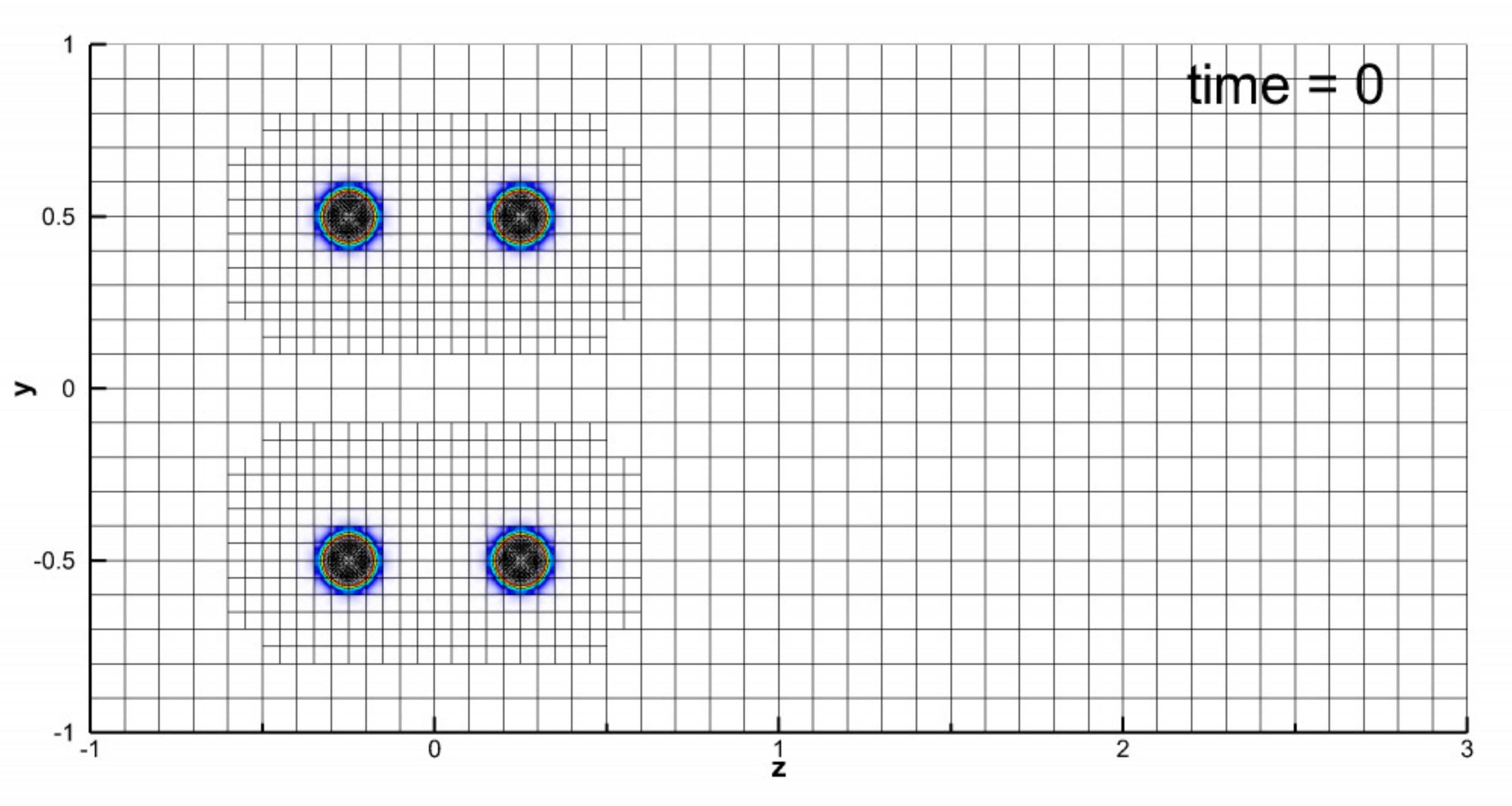}\;\includegraphics[width=0.45\textwidth]{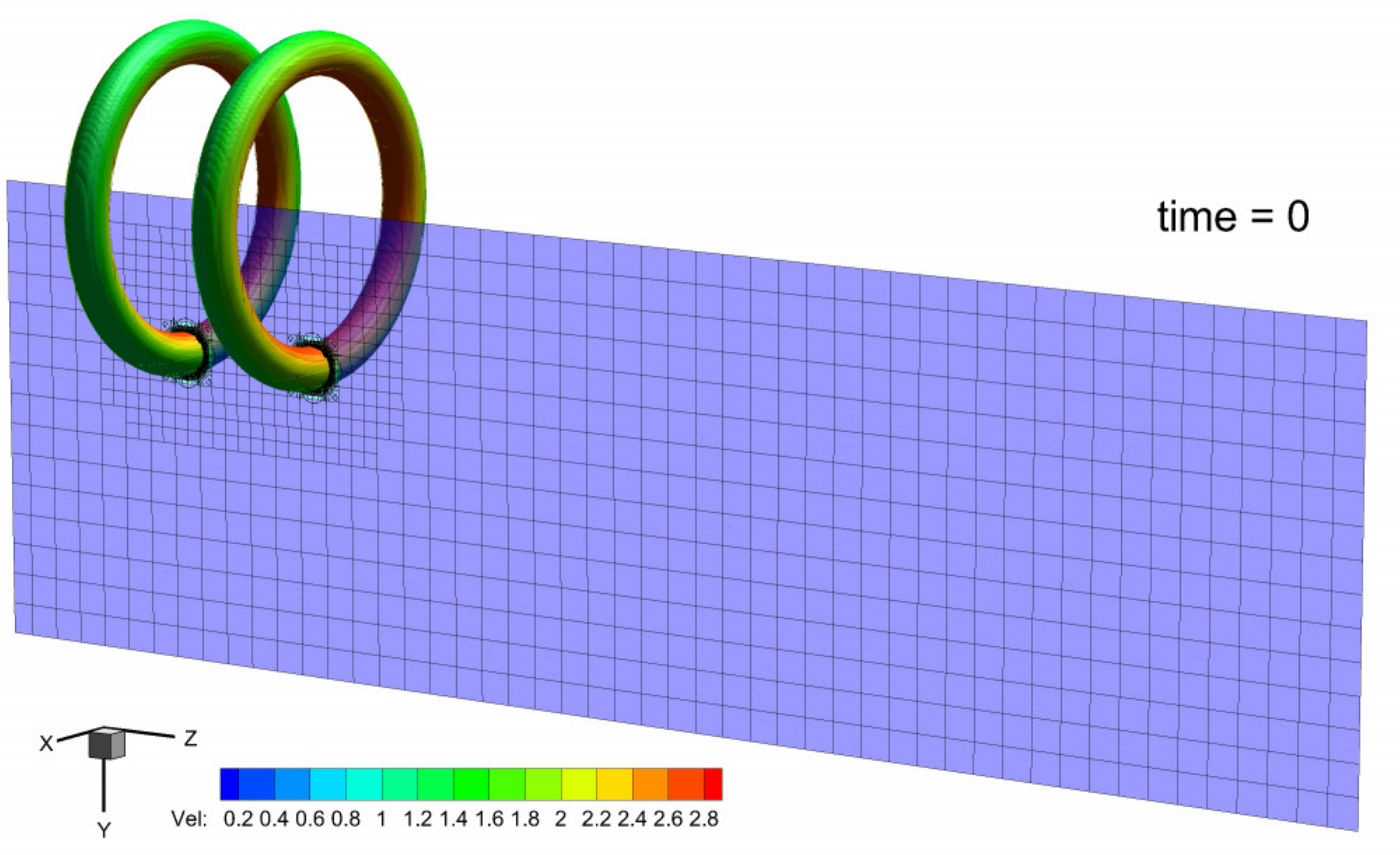}\\
			\includegraphics[width=0.45\textwidth]{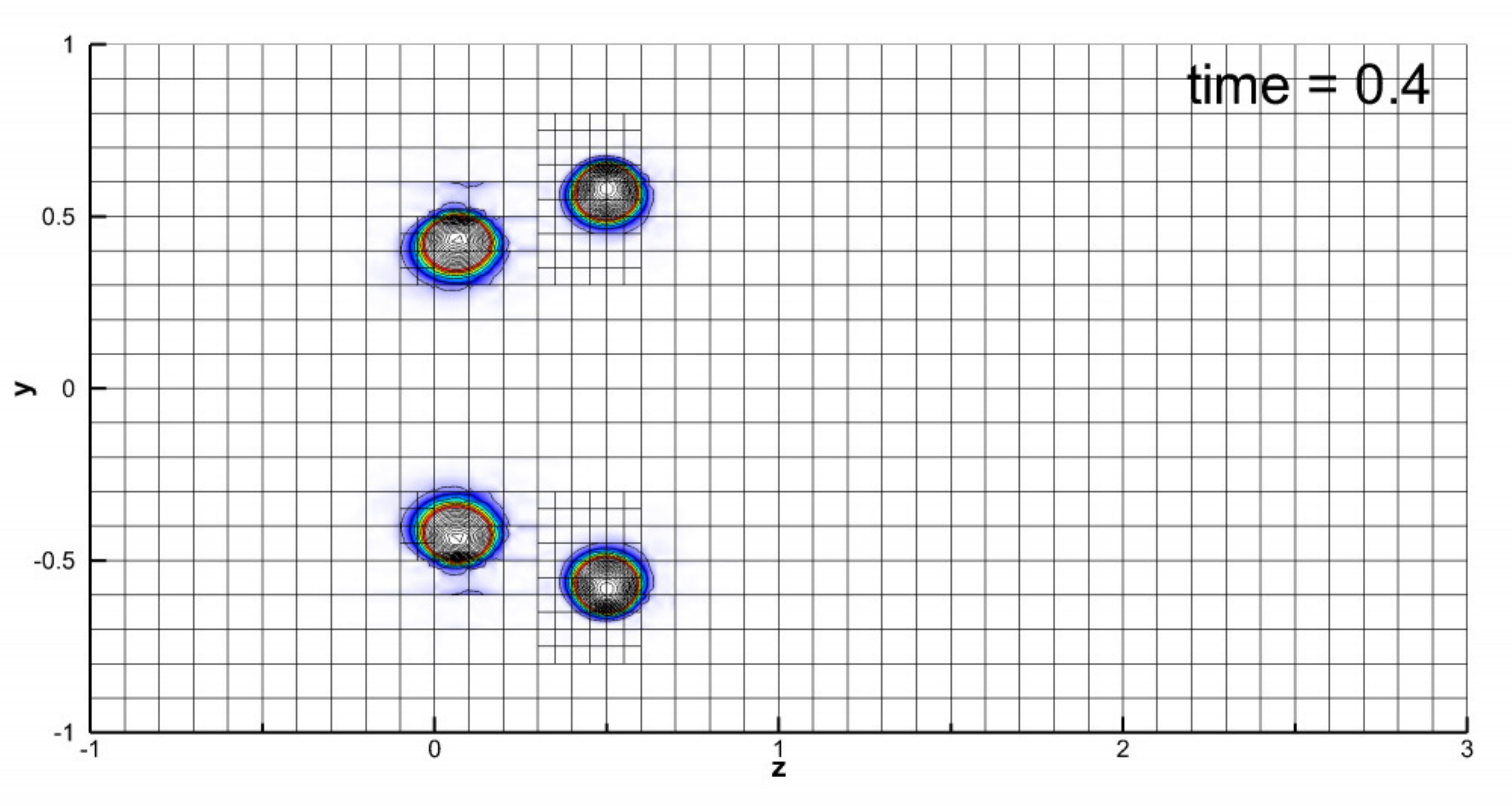}\;\includegraphics[width=0.45\textwidth]{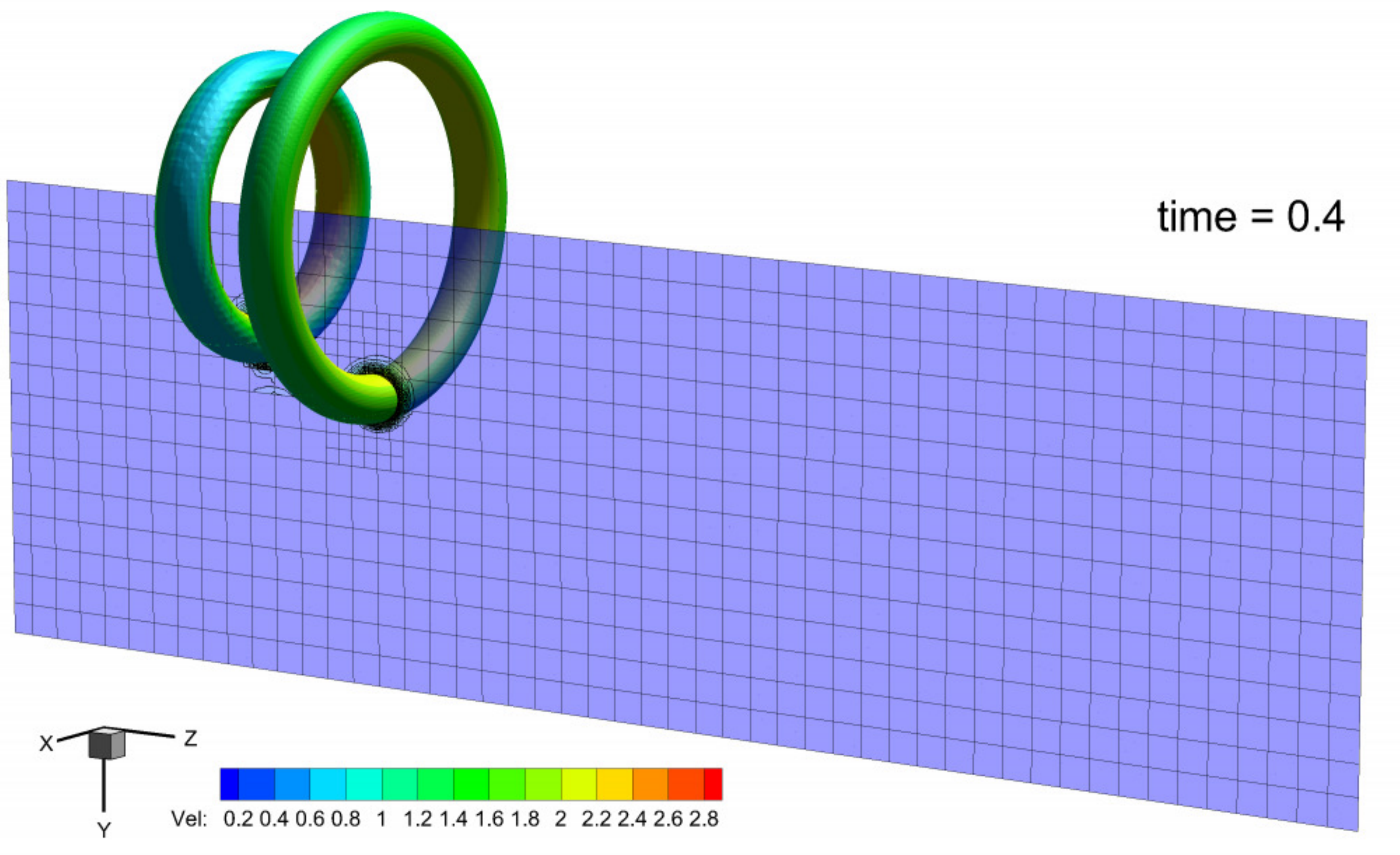}\\
			\includegraphics[width=0.45\textwidth]{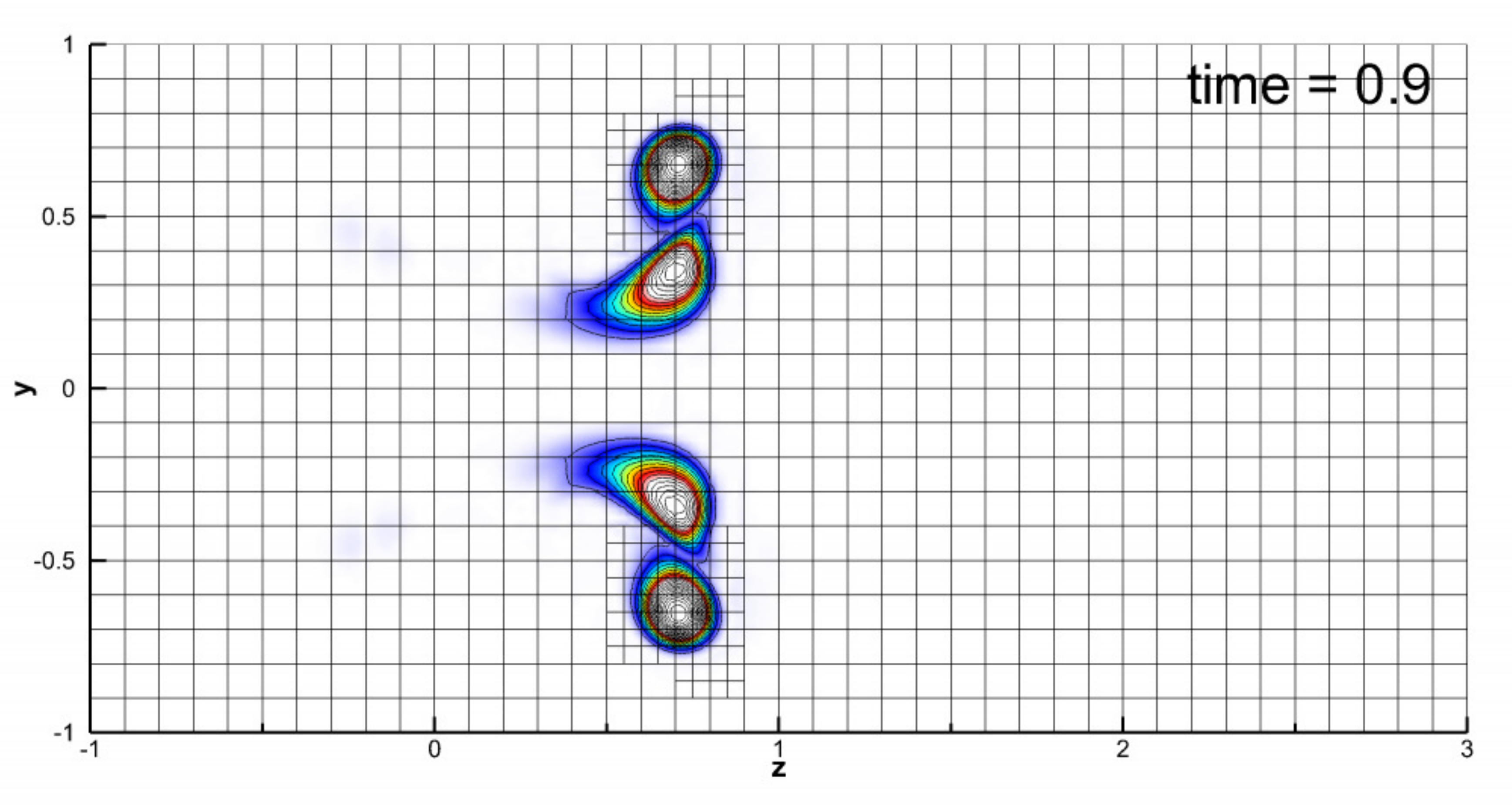}\;\includegraphics[width=0.45\textwidth]{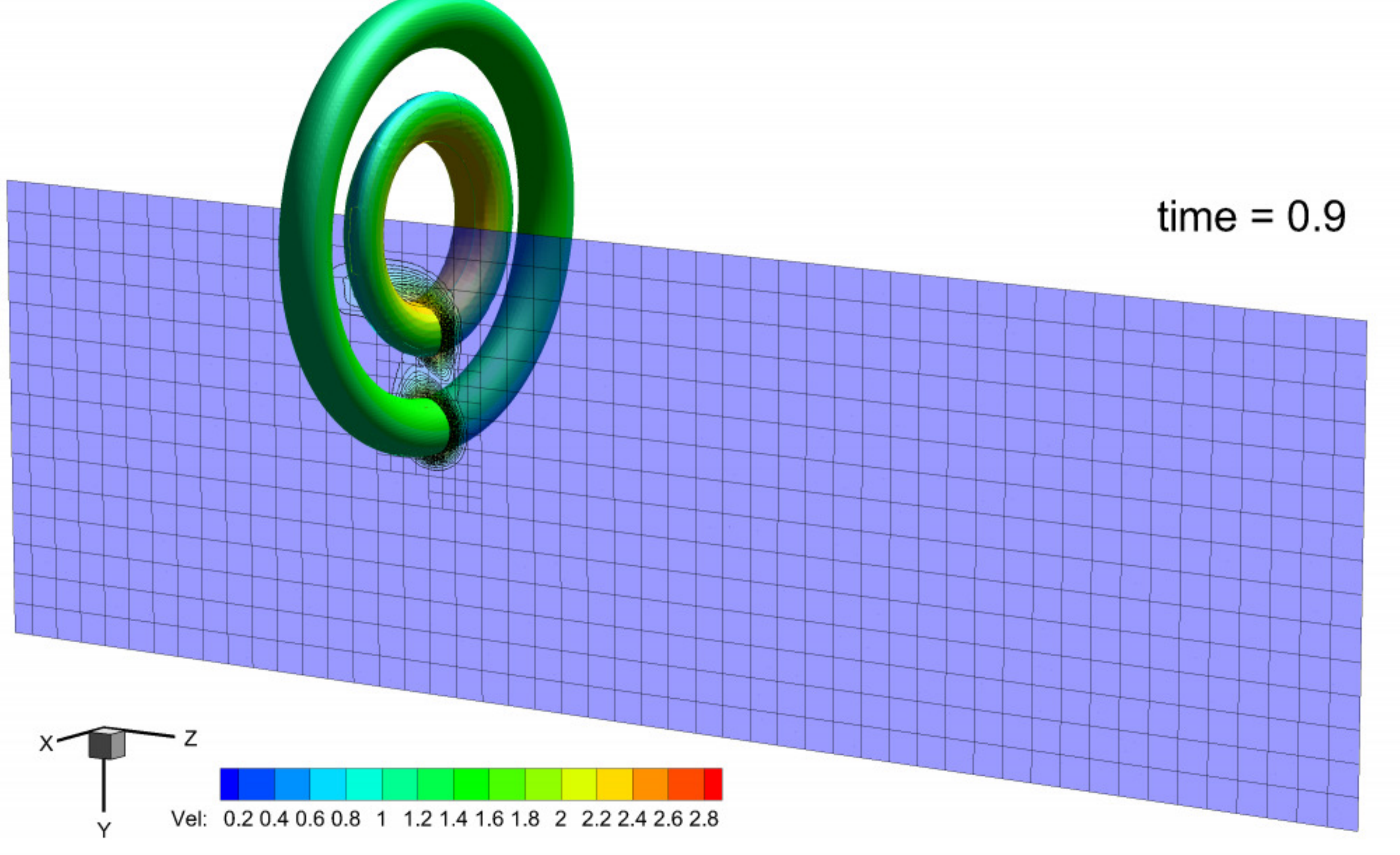} \\
			\includegraphics[width=0.45\textwidth]{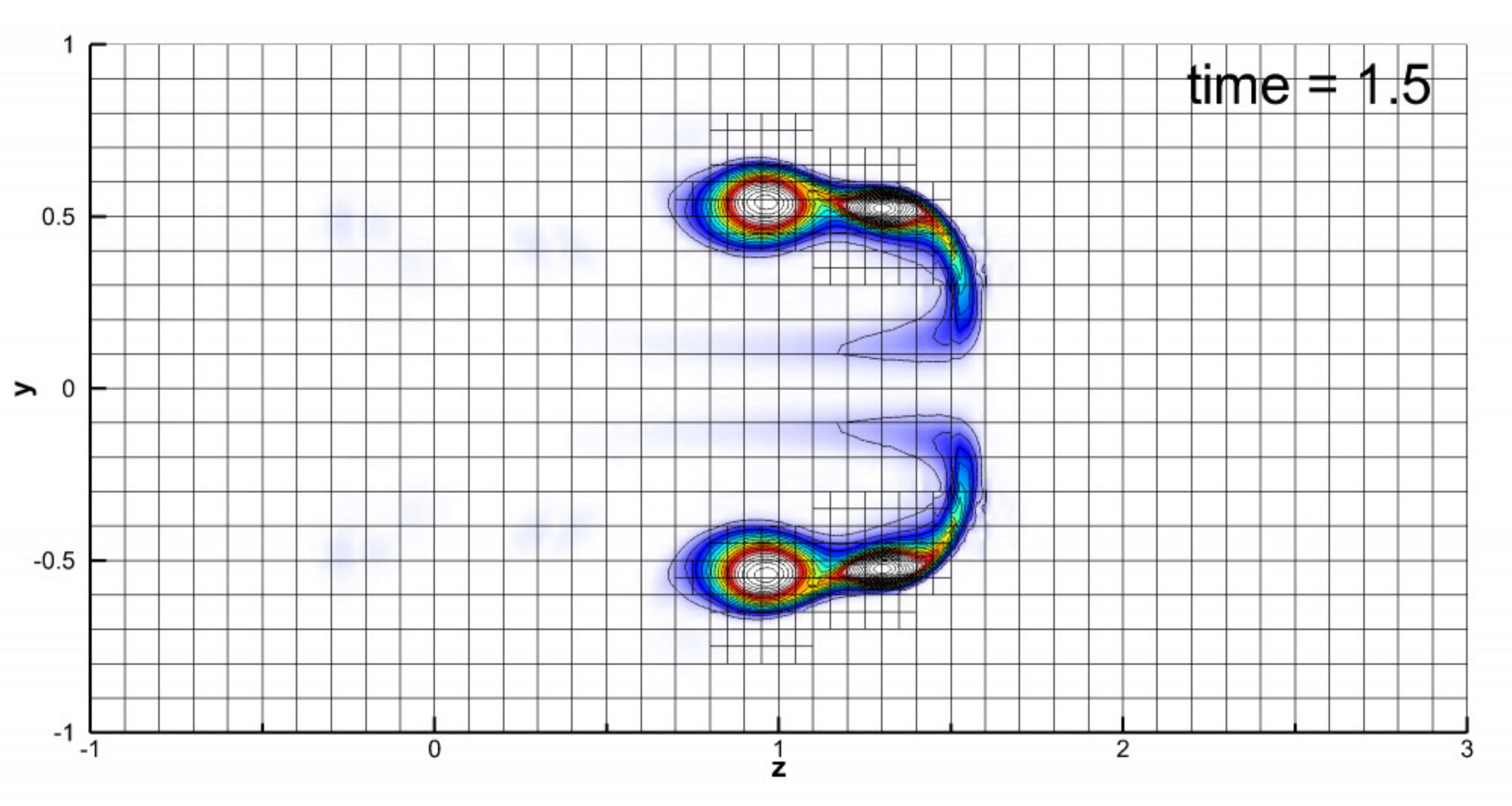}\;\includegraphics[width=0.45\textwidth]{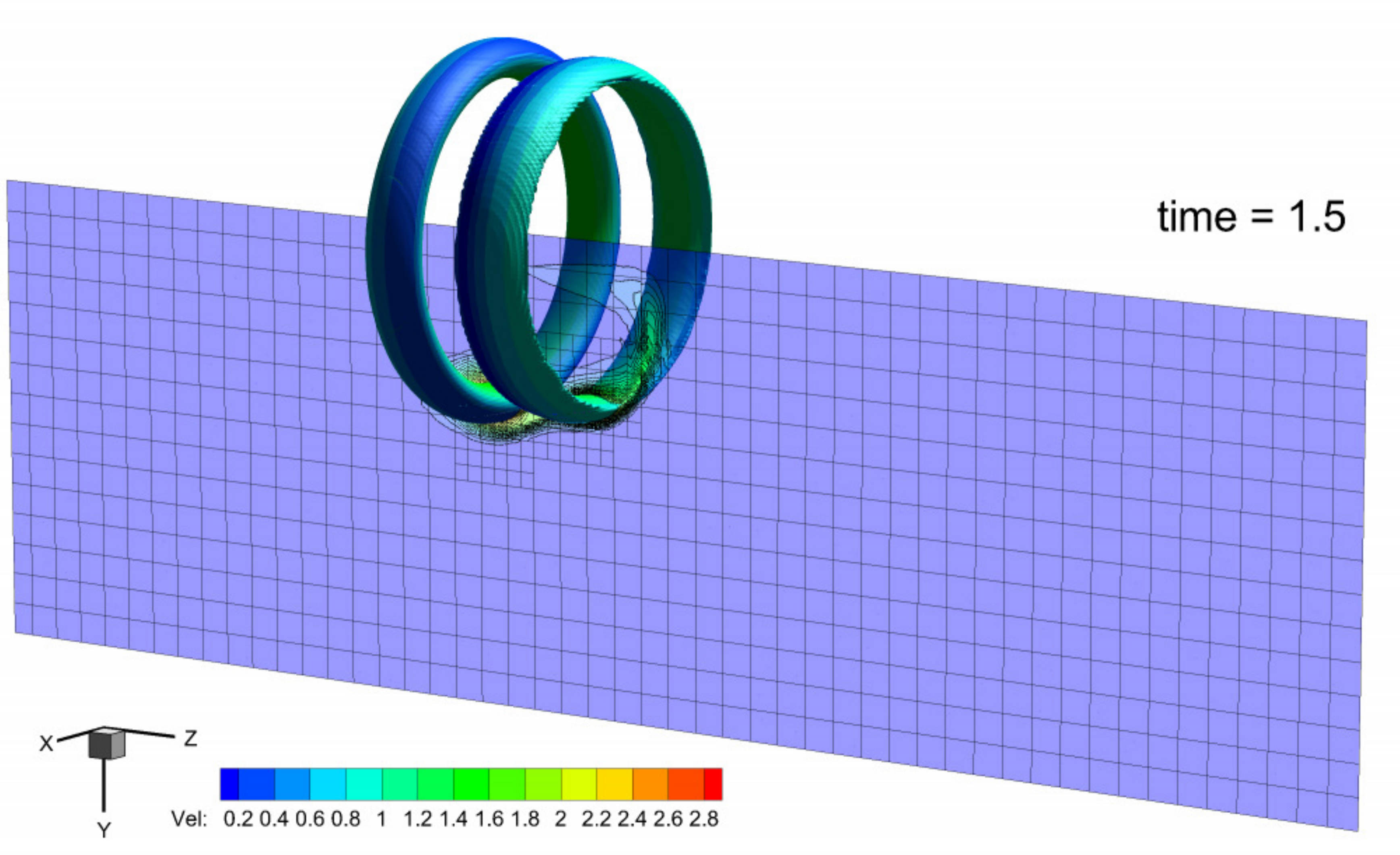} 
\caption{Time evolution of the vorticity field for the three-dimensional vortex-ring pair 'leapfrog' problem, respectively, from the top to the bottom:  t=$0.0$, $0.4$, $0.9$, and $1.5$, using $30\times 30 \times 80$ elements on the coarsest grid with periodic boundary conditions; at the left the numerical solution interpolated along the two-dimensional $y-z$ plane, at the right the three-dimensional view of the isosurfaces of the vorticity magnitude $|\omega|$; these results are obtained with the $\p_{4}$-version of our \SIDG method.}
\label{fig:VPL}
\end{figure}

\begin{figure} 
\centering 
			\includegraphics[width=0.45\textwidth]{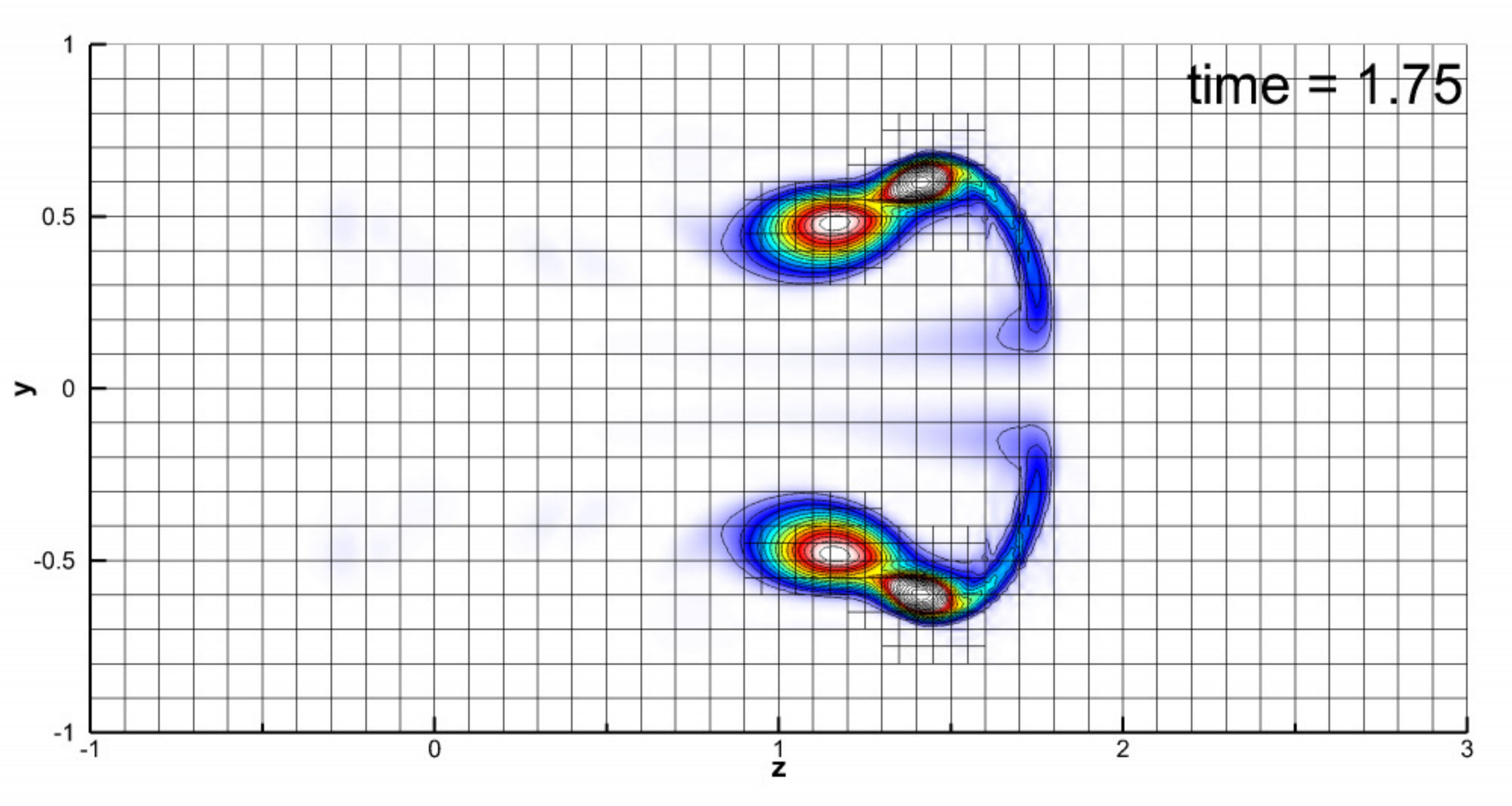}\;\includegraphics[width=0.45\textwidth]{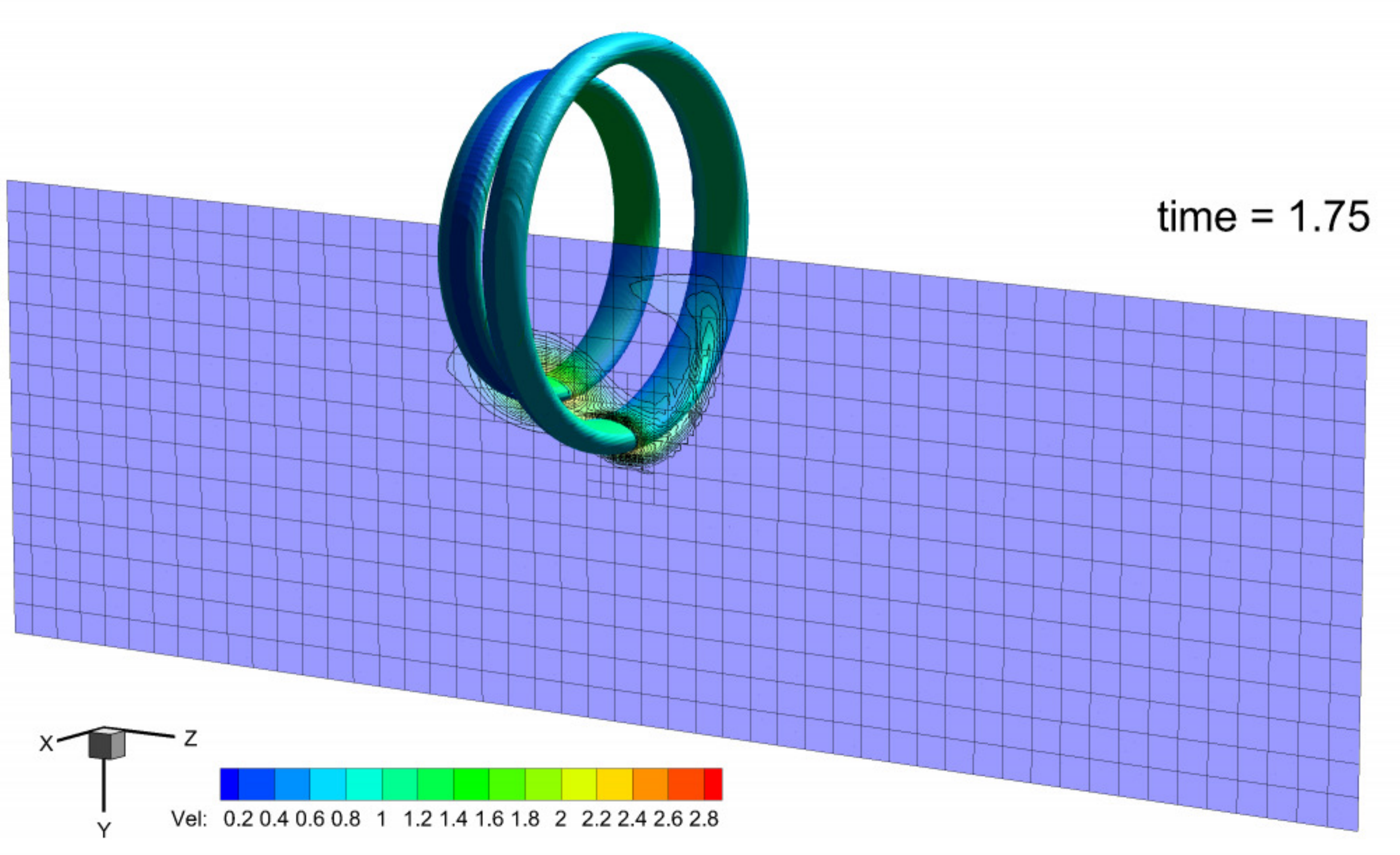}\\
			\includegraphics[width=0.45\textwidth]{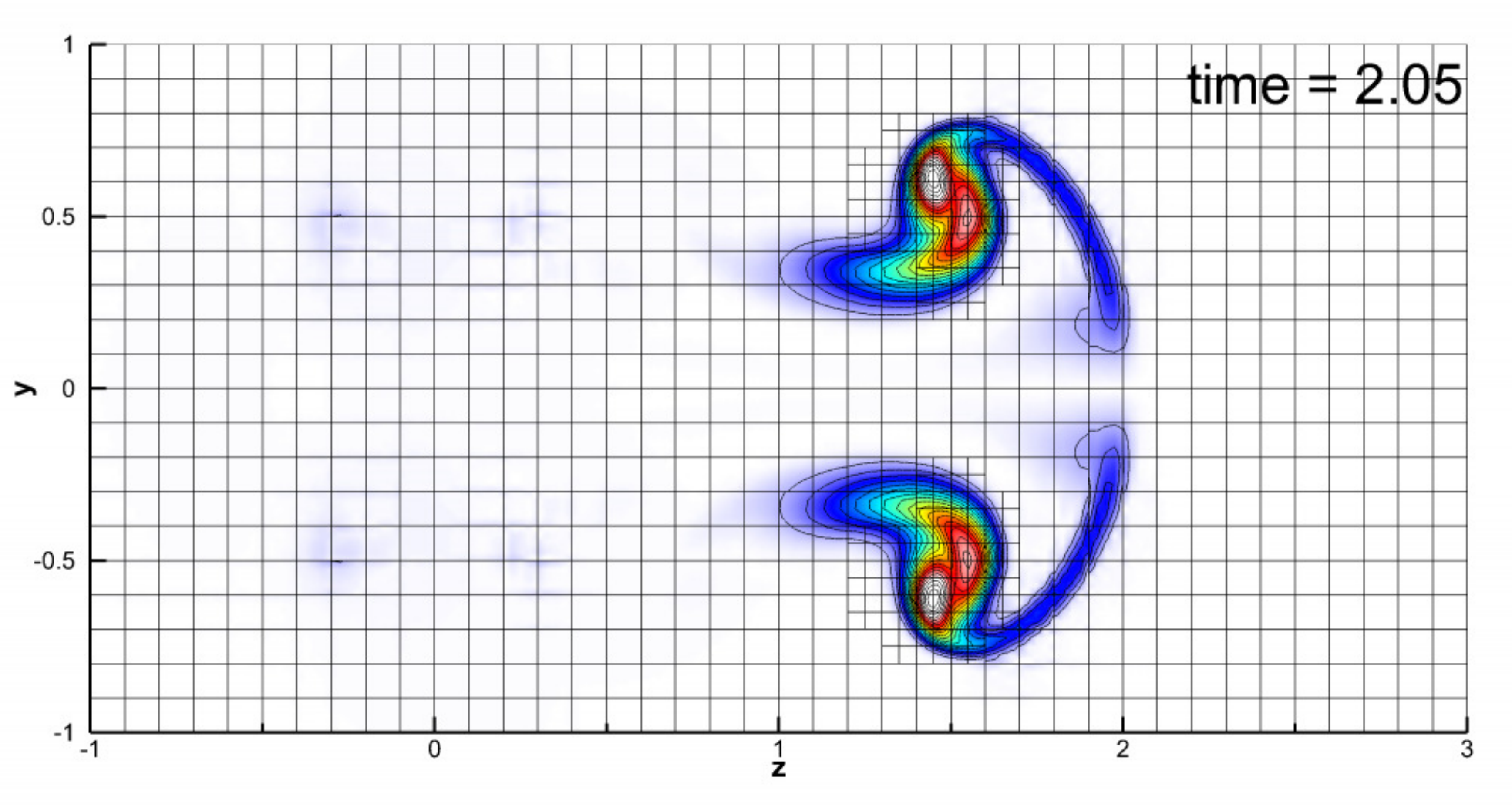}\;\includegraphics[width=0.45\textwidth]{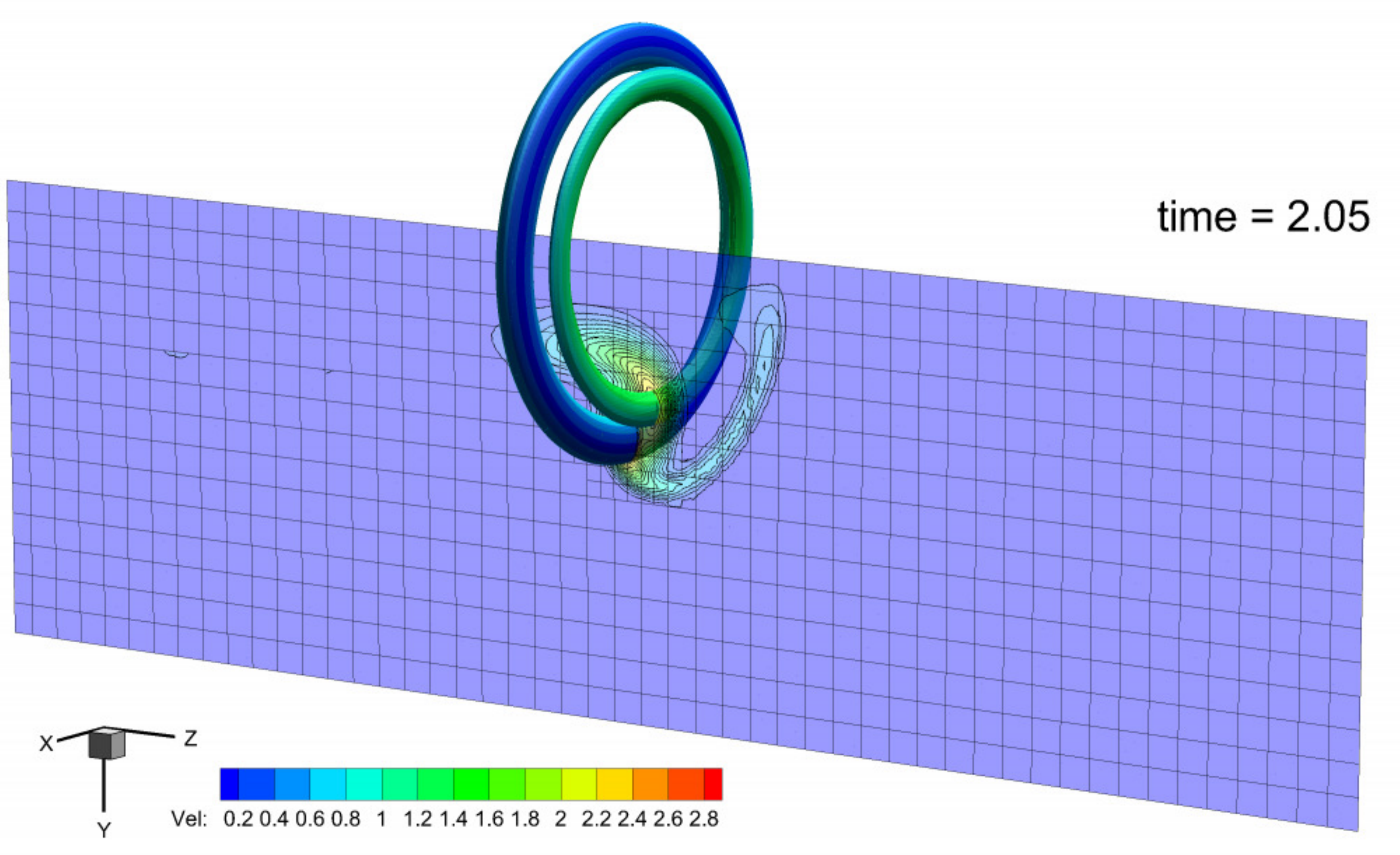}\\
			\includegraphics[width=0.45\textwidth]{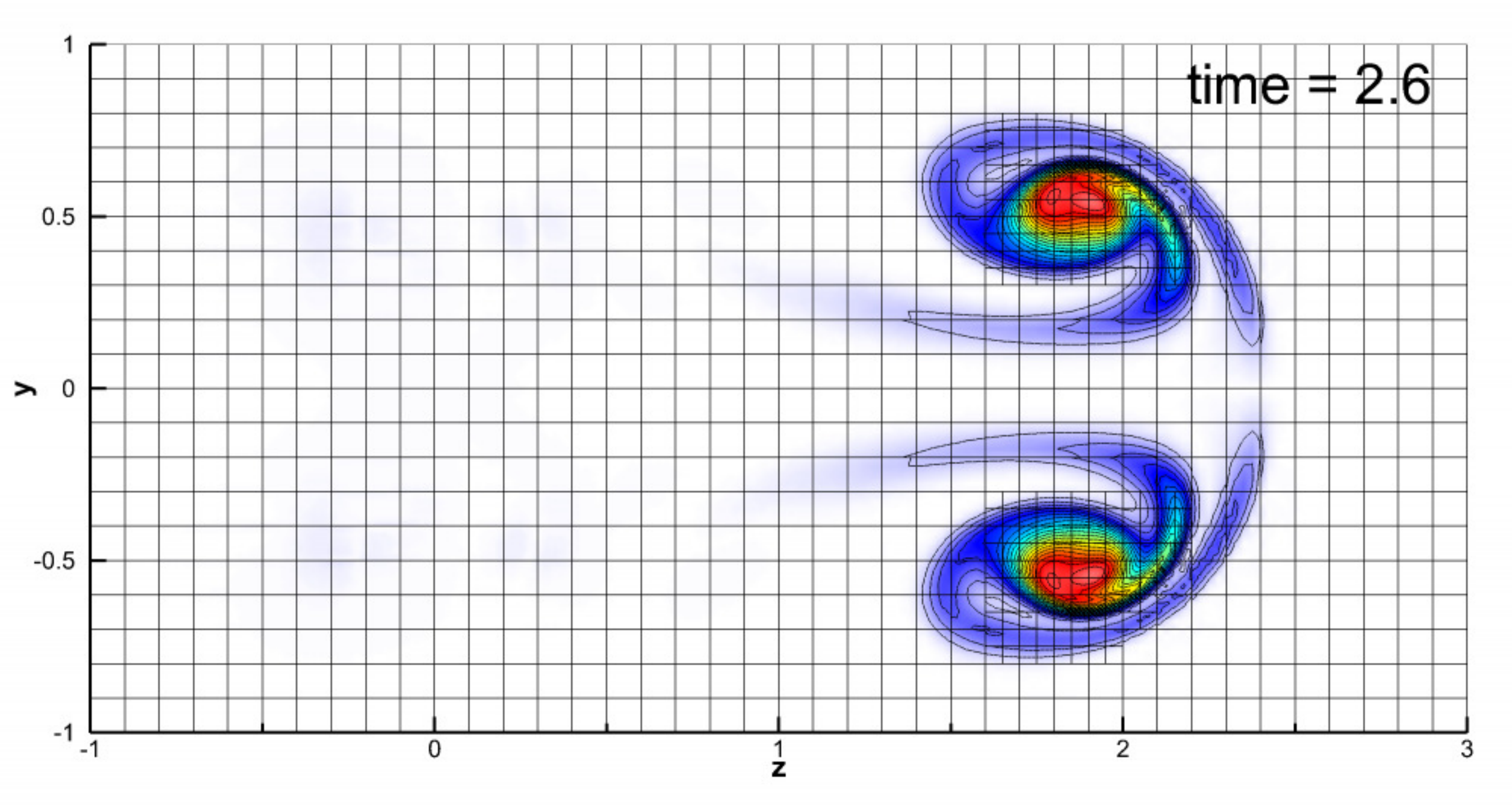}\;\includegraphics[width=0.45\textwidth]{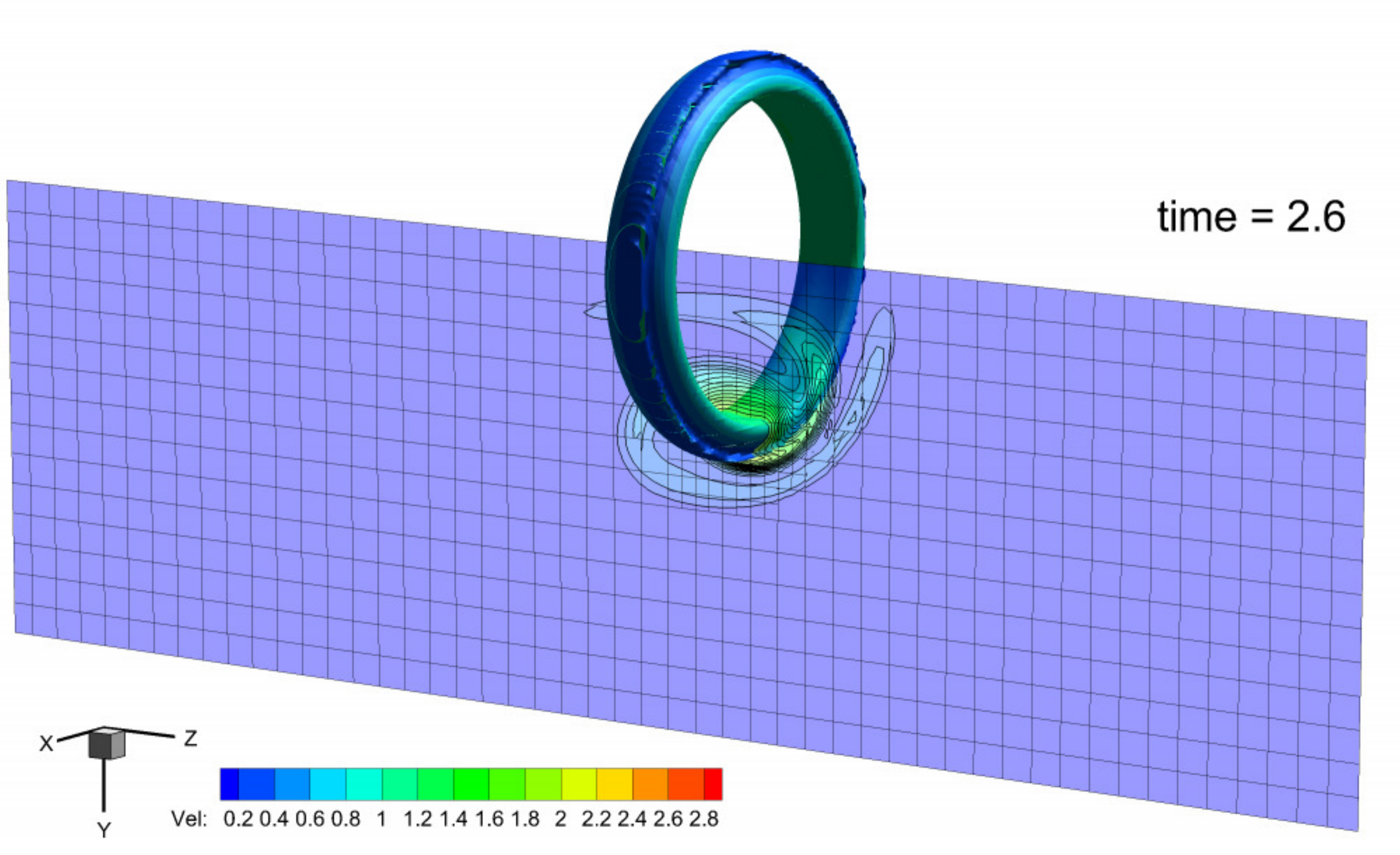} \\
			\includegraphics[width=0.45\textwidth]{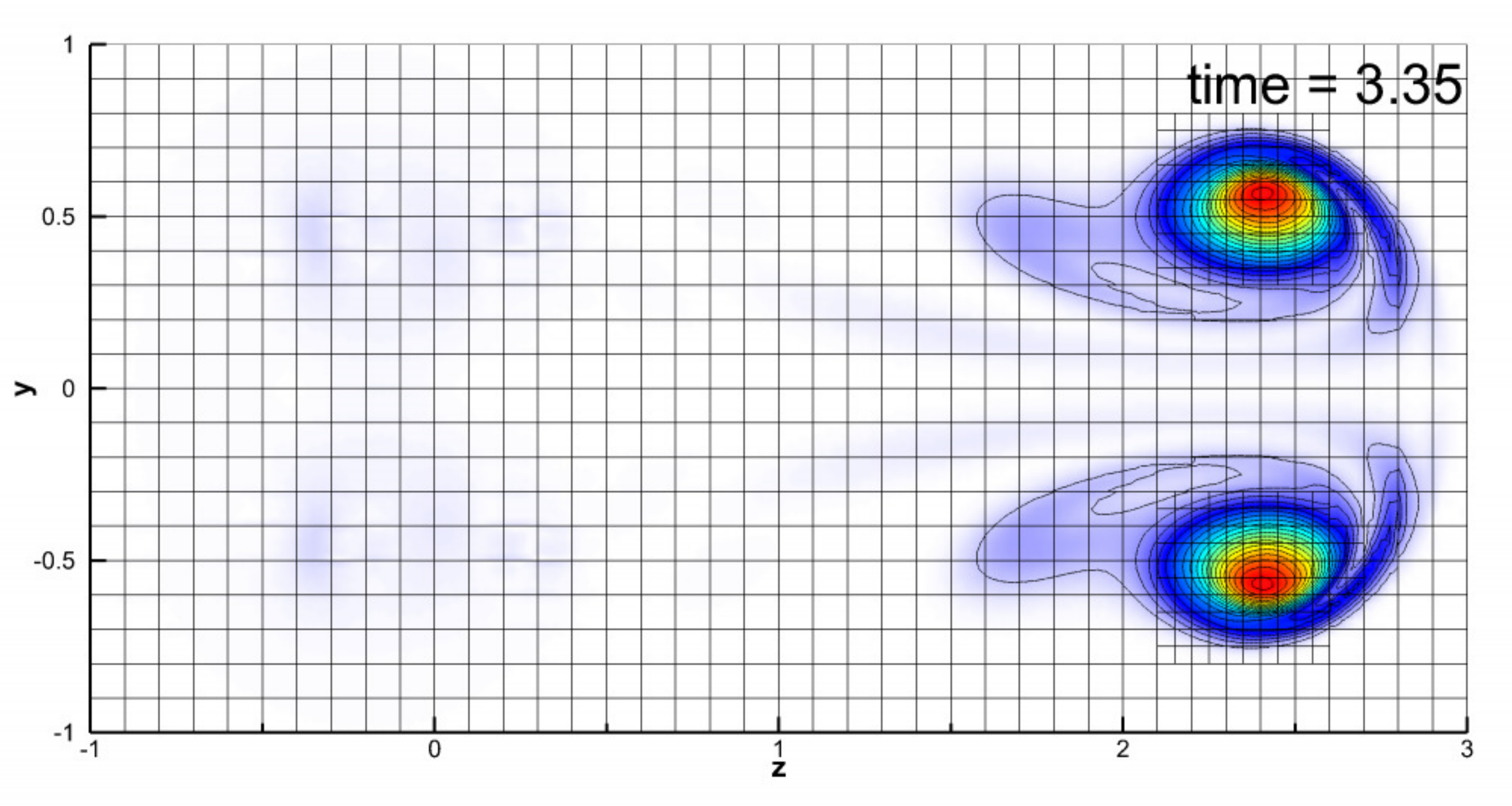}\;\includegraphics[width=0.45\textwidth]{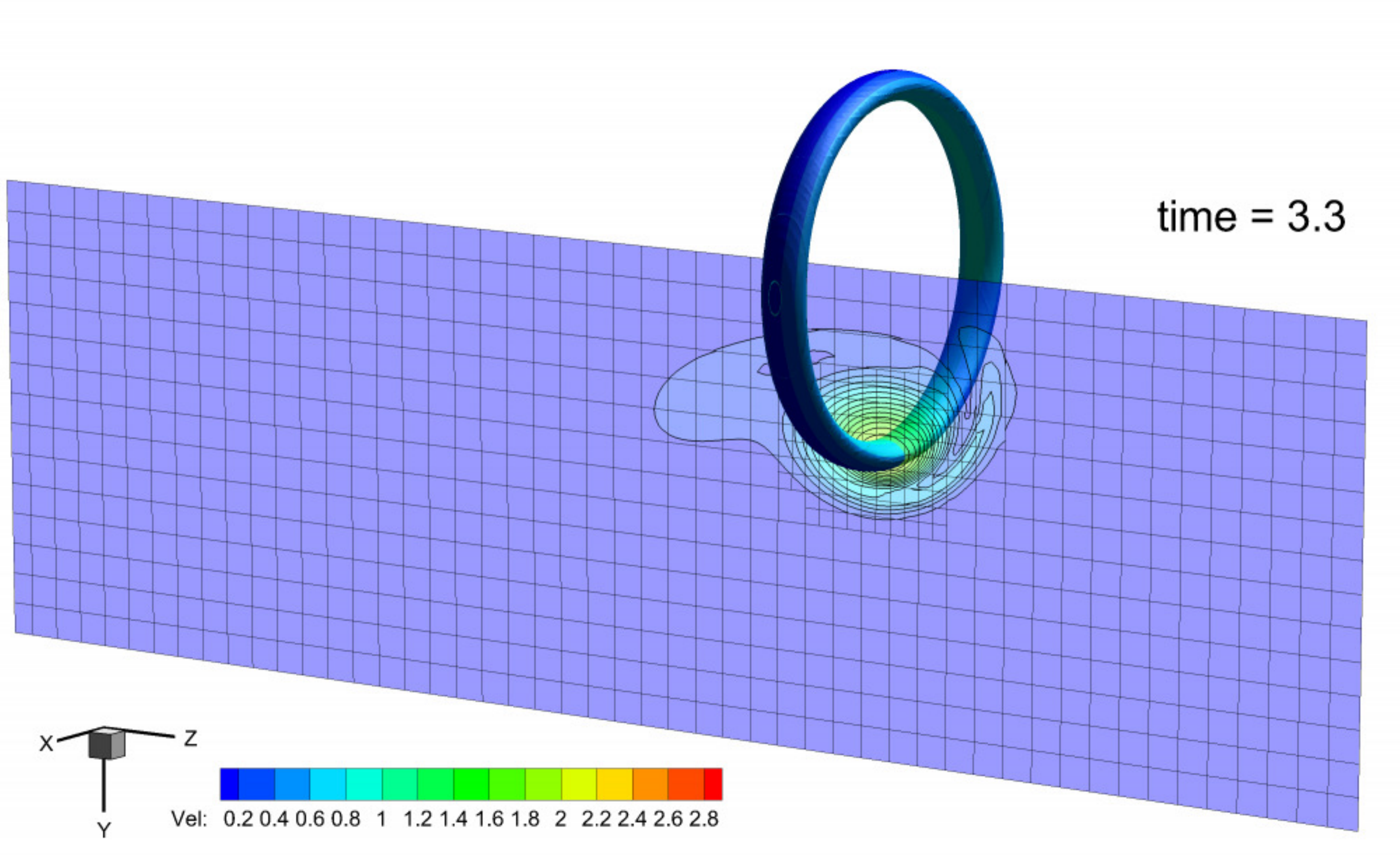} 
\caption{Time evolution of the vorticity field for the three-dimensional vortex-ring pair 'leapfrog' problem, respectively, from the top to the bottom:  t=$1.75$, $2.05$, $2.60$, and $3.35$, using $30\times 30\times 80$ elements on the coarsest grid with periodic boundary conditions; at the left the numerical solution interpolated along the two-dimensional $y-z$ plane, at the right the three-dimensional view of the isosurfaces of the vorticity magnitude $|\omega|$; these results are obtained with the $\p_{4}$-version of our \SIDG method.}\label{fig:VPL2}
\end{figure}

In this section the interaction between two coaxial vortex rings is simulated within the spatial domain $\Omega = [-1.5,1.5]^2\times[-4.4]$ with periodic boundary conditions everywhere. Two vortex rings are initialized according to the Gaussian distribution (\ref{eq:vort}), following \cite{Stanaway1988}, with a Reynolds number $Re_{\Gamma}=1000$, a ratio of major and minor radius of $R/r=10$, the ring  centers being \emph{one radius} $R$ away from each other. In particular the chosen parameters are $R_0=0.5$, $a=0.05$, vorticity amplitude $\omega_0=1$, and kinematic viscosity $\nu=10^{-3}$. In this case, the self-induction leads the ring pair to move vertically in the $z$ direction, the mutual-induction leads the last (backward) vortex to accelerate, being scaled down and going past the second vortex through the inner orifice, i.e. '\emph{leapfrogging}'. The time evolution for the computed vorticity magnitude interpolated along the two arbitrary (central symmetry holds) vertical and 
orthogonal planes is shown in Figures \ref{fig:VPL} and \ref{fig:VPL2}. Also in this case, good agreement with the provided reference solution of \cite{Stanaway1988,ChengLouLim2015} is verified. 
The physical domain $\Omega$ has been discretized within a mesh of $30\times 30\times 80$ space elements on the coarsest grid $\Omega_h^0$ within the space of solutions of our $\SIDG$-$\p_4$ method, refinement factor $\err=2$, up to one single refinement level $\ell_{\text{max}}=1$. The chosen AMR grid corresponds to a maximum number of  $N_{\text{dof}}^{\text{max}}=\Nel^{\text{max}} \times (N+1)^3 = 72\ 000 \ 000$ of degrees of freedom per physical variable if a uniform fine grid was used. Notice that, thanks to the AMR framework the real total number of degrees of freedom per physical variable is  reduced approximately by a factor of $\err^d=8$. Indeed, the mesh is dynamically refined only close to the vortex rings.

\section{Conclusion}
\label{sec:conclusions}

In this paper the family of staggered spectral semi-implicit DG methods for the incompressible Navier-Stokes equations recently proposed in \cite{FambriDumbser} for uniform Cartesian grids 
has been successfully extended to \emph{staggered AMR meshes} in two and three space dimensions. 
A similar formulation for staggered DG schemes on conforming \textit{unstructured} simplex meshes has been recently outlined in \cite{TavelliDumbser2014,TavelliDumbser2014b,TavelliDumbser2015}. 

It has been shown that the main advantage of edge/face-based staggered grid methods is to improve substantially the sparsity of the main linear system compared to traditional DG schemes on 
collocated grids. 
In particular the discrete Laplace operator $\mathbb{H}$ is shown to be symmetric and positive-definite and has maximum sparsity, since it only involves the element
and its direct edge/face neighbors. On uniform meshes, it is block penta-diagonal (in 2D) or block hepta-diagonal (in 3D). 
This has to be compared to classical DG methods on collocated grids: those are characterized by discrete operators with either larger computational stencils, 
if a discrete pressure Poisson equation is solved, since in this case the discrete pressure will also depend on neighbors of neighbors; or, one solves directly
the saddle point problem associated with the PDE system and thus keeps the small stencil, but in that case there are much more unknowns, namely the scalar pressure
and the components of the velocity vector. 
Moreover,  in our scheme $\mathbb{H}$ is shown to be the Schur complement of the discrete saddle point system of the incompressible Navier-Stokes equations. From another point of view,  
$\mathbb{H} \equiv \mathcal{D} ^T \mathbb{M}^{-1}\mathcal{D}$ accounts for the jumps of the piecewise polynomials in the discrete gradients $\mathcal{D}$ resembling a 
Bassi-Rebay-type lifting operator \emph{evaluated on the dual grid} $\Omega_h^{*}$, see equation (\ref{eq:D}). 
The same well-conditioned coefficient matrix $\mathbb{H}$ appears both in the pressure Poisson equation, but also in the linear system arising from the implicit discretization 
of the viscous terms, with only one additional symmetric positive definite block diagonal term coming from the element mass matrix associated with the time derivative of the 
velocity. As a consequence, the pressure system and the viscous systems can be solved very efficiently by means of a classical matrix-free conjugate gradient method without 
recurring, in this paper, to any preconditioner. The spectral properties of the system matrix have been theoretically analyzed in \cite{SIDG_analysis2016} for the uniform Cartesian 
grid case. 

The final algorithm is verified to be high order accurate in space with spectral convergence property after completing a thorough numerical convergence test that accounts also for the  
AMR grid refinement. The presented $\SIDG$-$\p_N$ method has been observed to be stable, robust and very accurate. The performance of the method has been shown on a large set of 
non-trivial test cases in two and three space-dimensions.

Future research will concern the extension of the present scheme on AMR meshes to a staggered space-time DG method in order to achieve also high order of accuracy in time, see 
\cite{FambriDumbser,TavelliDumbser2016}. 
%
%
%
Further work will also concern the extension to the compressible Euler and Navier-Stokes equations, following the ideas put forward in \cite{DumbserCasulli2016} for pressure-based 
semi-implicit finite volume schemes on staggered grids. In this case, a proper limiter will be needed, not only for the treatment of shock waves, but also in order to avoid spurious 
oscillations in the presence of steep gradients or under-resolved flow features, for example following the novel ideas on \textit{a posteriori} subcell finite volume limiters
outlined in \cite{Dumbser2014,Zanotti2015c,Zanotti2015d}. 

\section*{Acknowledgments}
The research presented in this paper was funded by the European Research Council (ERC) under the European Union's Seventh Framework
Programme (FP7/2007-2013) within the research project \textit{STiMulUs}, ERC Grant agreement no. 278267.

The authors also acknowledge the Leibniz Rechenzentrum (LRZ) in Munich, Germany, for awarding access to the \textit{SuperMUC} 
supercomputer, as well as the support of the HLRS in Stuttgart, Germany, for awarding access to the \textit{Hazel Hen} supercomputer.


\bibliographystyle{plain} 
\bibliography{./references7} 


\end{document}